\documentclass[1p]{elsarticle}

\usepackage[utf8]{inputenc}
\usepackage{hyperref}
\usepackage{amsfonts}
\usepackage{amsmath,amssymb}
\usepackage{bbm}
\usepackage{bigints}
\usepackage{multirow,booktabs}
\usepackage{mathtools}
\usepackage{cleveref}
\usepackage[caption = false]{subfig}
\usepackage[table]{xcolor}

\newcommand{\R}{\mathbbm{R}}
\newcommand{\N}{\mathbbm{N}}

\newcommand{\cG}{\Pi}
\newcommand{\cP}{\mathcal{P}}
\newcommand{\hh}{\widehat{h}}
\newcommand{\bh}{\overline{\boldsymbol{h}}}
\newcommand{\hU}{\widehat{U}}
\newcommand{\bU}{\overline{\boldsymbol{U}}}
\newcommand{\hu}{\widehat{u}}
\newcommand{\hv}{\widehat{v}}
\newcommand{\hqx}{\widehat{q^x}}
\newcommand{\hqy}{\widehat{q^y}}
\newcommand{\bqx}{\overline{\boldsymbol{q}^x}}
\newcommand{\bqy}{\overline{\boldsymbol{q}^y}}
\newcommand{\co}{\text{cor}}
\newcommand{\rone}[1]{{\leavevmode\color{black}{#1}}}

\newtheorem{thm}{Theorem}[section]
\newtheorem{lem}{Lemma}[section]
\newdefinition{rmk}{Remark}[section]
\newproof{pf}{Proof}[section]
\newtheorem{corollary}{Corollary}[thm]
\newtheorem{proposition}{Proposition}[section]

\journal{}
\begin{document}

\begin{frontmatter}
\title{Hyperbolicity-Preserving and Well-Balanced Stochastic Galerkin Method for Two-Dimensional Shallow Water Equations}
\author{Dihan Dai\corref{cor1}%
\fnref{fn1}}
\ead{dai@math.utah.edu}
\author{Yekaterina Epshteyn\fnref{fn1}}
\ead{epshteyn@math.utah.edu}
\author{Akil Narayan\fnref{fn1,fn2}}
\ead{akil@sci.utah.edu}

\cortext[cor1]{Corresponding author}
\address[fn1]{Department of Mathematics, University of Utah, Salt Lake City, UT 84112, USA}
\address[fn2]{Scientific Computing and Imaging (SCI) Institute, University of Utah, Salt Lake City, UT 84112, USA}
\begin{abstract}
Stochastic Galerkin formulations of the two-dimensional shallow water systems parameterized with random variables may lose hyperbolicity, and hence change the nature of the original model. In this work, we present a hyperbolicity-preserving stochastic Galerkin formulation by carefully selecting the polynomial chaos approximations to the nonlinear terms in the shallow water equations. We derive a sufficient condition to preserve the hyperbolicity of the stochastic Galerkin system which requires only a finite collection of positivity conditions on the stochastic water height at selected quadrature points in parameter space. Based on our theoretical results for the stochastic Galerkin formulation, we develop a corresponding well-balanced hyperbolicity-preserving central-upwind scheme. We demonstrate the accuracy and the robustness of the new scheme on several challenging numerical tests.  
\end{abstract}
\begin{keyword}
   finite volume method, stochastic Galerkin method, shallow water equations, hyperbolic systems of conservation and balance laws\\\vspace{1em}
\end{keyword}
\end{frontmatter}

\section{Introduction}
The system of shallow water equations and closely related models have found tremendous success in modelling flows whose vertical scales are significantly smaller than their horizontal scales. Such flows include, but are not limited to, the water flows in rivers, lakes, and coastal areas. In addition, uncertainties can enter models, for instance, through measuring or empirically approximating the bottom topography or initial conditions. This makes deterministic-model-based predictions less reliable and therefore it is also important to model this uncertainty in the governing shallow water systems. In this work, we will consider the \emph{parameterized} stochastic shallow water equations (SWE), where the parameter is a random variable that models the uncertainty.
 
For discretizing parameterized partial differential equations, the polynomial chaos expansion (PCE) technique is a commonly-used approach. It seeks to represent a second-order random field using an expansion in an orthogonal polynomial basis. It was first introduced by Wiener \cite{wiener1938homogeneous} to represent the Gaussian processes using Hermite polynomials. The idea has been later extended to the representations of random outputs and fields with more general distributions using more general families of orthogonal polynomials \cite{ghanem1991stochastic,xiu2002wiener,le2004uncertainty,wan2005adaptive}. 
With PCE the strategy for encoding uncertainty, this gives rise to parameterized partial differential equations. Two major classes of methods are used to solve the resulting system and compute the statistics of the outputs. One class of \emph{non-intrusive} methods is based on sampling (see, e.g.,\cite{xiu2005high, nobile2008sparse, mishra2012multilevel}). The implementations and the parallelism can be easily achieved for non-intrusive methods. However, such approximations can be less accurate than the solutions of intrusive methods.  

The class of \emph{intrusive methods}, in the PCE context, usually refers to stochastic Galerkin (SG) method. In the SG framework, the truncated PCE representations are projected by a Galerkin method in stochastic (parameter) space to obtain a new system of deterministic partial differential equations (PDEs) that describes the evolution of the coefficients in the PCE representation. The new system typically requires new solvers and a substantial rewrite of legacy codes. However, it is expected that the SG method will deliver a more accurate solution than, for example, some alternative non-intrusive strategies. In particular, it is known that intrusive methods have near-optimal accuracy in an $L^2$ sense for static problems  \cite{babuska2004galerkin,le2010spectral}. In addition, for SG formulations, one can generally prove important properties of the system or of their subsequent numerical discretizations. Because of this, applications of SG methods have found success in modeling uncertainty in diffusion equations, e.g. \cite{xiu2009efficient, eigel2014adaptive}, kinetic equations, e.g. \cite{hu2016stochastic,shu2017stochastic}, and for conservation and balance laws with symmetric Jacobian matrices, e.g. \cite{tryoen2010intrusive}.

However, when the SG method is applied to a general hyperbolic system of conservation/balance laws such as the SWE, the associated SG system may lose hyperbolicity \cite{despres2013robust, jin2019study,gerster2019hyperbolic}. The loss of hyperbolicity may lead to unphysical solutions and prevent the use of robust numerical schemes. Efforts have been made to overcome this difficulty. In \cite{wu2017stochastic}, a strategy to regain hyperbolicity for the quasilinear form of hyperbolic systems is proposed. Hyperbolicity can be regained by multiplying the SG formulation of the system by the left eigenvector matrix of the flux Jacobian matrix, but this approach is limited to quasilinear forms, and the numerical schemes designed for conservative forms cannot be applied directly. Operator-splitting approaches for one-dimensional Euler equations \cite{chertock2015operator} and one-dimensional shallow water equations \cite{chertock2015welluq} have also been investigated. The original systems are split into linear hyperbolic and scalar subsystems, whose SG formulations remain hyperbolic. However, it has been shown in \cite{schlachter2018hyperbolicity} that the splitting approach may lead to an overall non-hyperbolic system. Strategies that introduce a PCE of auxiliary variables have also been explored. For example, the SG system of a balance/conservation laws in terms of the entropic variables can be shown to be hyperbolic \cite{poette2009uncertainty,poette2019contribution}. An optimization-based intrusive polynomial moment method (IPMM) was proposed to compute the PCE of entropic variables given the PCE of the conserved variables \cite{despres2013robust,poette2009uncertainty,poette2019contribution}. But the IPMM method can be computationally expensive since one must solve an optimization problem for each cell at each discrete time level. The PCE using a Roe variable formulation is introduced in \cite{pettersson2014stochastic,gerster2019hyperbolic,gerster2020entropies}, where the flux of the SG system is constructed using PCE of the Roe variables. Thus, the conservative form of the system is preserved. It has been shown that, for the Wiener-Haar expansion, both the SG Roe formulations of the Euler equation \cite{pettersson2014stochastic} and of the 1D SWE \cite{gerster2019hyperbolic} are hyperbolic. For the system of isothermal Euler equations, its SG Roe formulation is hyperbolic under a positive definiteness condition of a matrix depending on the positive square root of the water height, regardless of the choices of distributions of the random parameters \cite{gerster2019hyperbolic}. However, the Roe formulation may still be expensive to implement since the PCE of the Roe variables need to be calculated by solving both a nonlinear equation and a linear equation. More discussion about the SG Roe formulation can be found in \cite{gerster2020entropies}. Recently we developed a hyperbolicity-preserving SG formulation for the one-dimensional shallow water equations \cite{doi:10.1137/20M1360736}, where we carefully tailored PCE truncation operators for the nonlinear non-polynomial term in the system using only the conserved variables. For arbitrary choices of distributions of the random parameters, the resulting SG system is hyperbolic under only a positive-definiteness condition on a matrix which depends on the PCE of the water height. The condition is a stochastic variant of a positivity condition for the water height for the deterministic system. The current paper is an extension of the work \cite{doi:10.1137/20M1360736} to the two dimensional shallow water equations.

In this work, the numerical scheme for the SG formulation in two-dimensional physical space is based on the extension of the finite-volume central-upwind scheme. The semi-discrete central-upwind schemes \cite{kurganov2001semidiscrete} are a class of robust Godunov-type Riemann-solver-free projection-evolution methods for hyperbolic systems. They are the generalizations of the central Nessyahu-Tadmor scheme \cite{nessyahu1990non} and Kurganov-Tadmor high resolution central scheme \cite{kurganov2000new}. The central Nessyahu-Tadmor scheme, the Kurganov-Tadmor central scheme, and the family of the central-upwind schemes have been successfully applied for the solution of the general multi-dimensional hyperbolic systems of conservation/balance laws, including for the solution of the deterministic SWE and related models \cite{kurganov2018finite}. A second-order central-upwind scheme for the solution of SWE was first proposed in \cite{kurganov2002central}, but it was not designed to preserve the positivity and the well-balanced property simultaneously. The scheme was further improved in \cite{kurganov2007second} to capture the ``lake-at-rest'' steady state and to preserve the positivity of the water height. See \cite{kurganov2007adaptive,kurganov2007reduction,M2AN_2011__45_3_423_0,chertock2015well,liu2018well,kurganov2018finite} for more examples and details of other closely related work. The scheme in this paper is mainly based on further extensions to stochastic SWE of the schemes proposed in \cite{kurganov2001semidiscrete, kurganov2007second}.  We note that, although the developed framework in \cite{doi:10.1137/20M1360736} and in the current paper is not limited to the choice of central-upwind schemes, the Riemann-solver-free nature of these schemes allows us to design an efficient deterministic solver as a part of the numerical discretization of the random SWEs.

\subsection{Contributions of This Paper}
In this paper we extend the ideas that we developed in \cite{doi:10.1137/20M1360736} to stochastic shallow water equations in two-dimensional physical space. In order to develop a hyperbolicity-preserving stochastic Galerkin formulation, we use specially tailored PCE truncation operators for the SWE nonlinear term $(q^xq^y)/h$ that occurs in both the $x$- and the $y$- component of the flux. For any probability distributions of the random parameters, our SG formulation is conservative and we prove it is also hyperbolic under a positive-definiteness condition on the matrix associated with the PCE of the water height (\Cref{thm:hyperbolicity}). This result is similar to our one-dimensional physical space result in \cite{doi:10.1137/20M1360736} and is a generalization to the results in \cite{gerster2019hyperbolic} for matrices related to Roe variables. This positive-definiteness condition is a stochastic variant of the positivity condition on the water height for the deterministic model. We show that the positive-definiteness condition can be guaranteed by enforcing a finite collection of positivity conditions on the stochastic water height at the stochastic quadrature nodes. Based on our SG formulation, we also develop a corresponding well-balanced hyperbolicity-preserving central-upwind scheme. 

This paper is structured as follows. In \Cref{sec:model}, we introduce the parameterized stochastic shallow water model in 2D physical space, as well as the stochastic Galerkin discretization of the system using a particular choice of the PCE truncation for the nonlinear terms $(q^x)^2/h,(q^y)^2/h$, and $(q^xq^y)/h$. In \Cref{sec:hyperbolicity}, we discuss the hyperbolicity of the resulting stochastic Galerkin system and present a sufficient condition to guarantee hyperbolicity of the stochastic Galerkin formulation of the 2D shallow water system. In \Cref{sec:scheme} we develop a well-balanced central-upwind scheme for the resulting model and obtain a hyperbolicity-preserving CFL-type condition. In addition, we use a filter to avoid the possible loss of hyperbolicity at pointwise reconstructed values. In \Cref{sec:results} we demonstrate the accuracy and the robustness of the developed numerical scheme on several challenging numerical tests. Conclusions and future directions are discussed in \Cref{sec:conclusion}.

\section{Modeling Stochastic Shallow Water Equations}\label{sec:model}
\subsection{Stochastic Modelling of the Two-Dimensional Shallow Water Equations}
The two-dimensional system of shallow water equations read 
\begin{equation}\label{eq:swesg1}
    \begin{aligned}
          &\dfrac{\partial h}{\partial t}+\dfrac{\partial q^x}{\partial x}+\dfrac{\partial q^y}{\partial y} = 0,\\[10pt]
          &\dfrac{\partial q^x}{\partial t}+\dfrac{\partial }{\partial x}\left(\frac{(q^x)^2}{h}+\frac{gh^2}{2}\right)+\dfrac{\partial }{\partial y}\left(\frac{q^xq^y}{h}\right)= -gh\dfrac{\partial B}{\partial x},\\[10pt]
          &\dfrac{\partial q^y}{\partial t}+\dfrac{\partial}{\partial x}\left(\frac{q^xq^y}{h}\right)+\dfrac{\partial}{\partial y}\left(\frac{(q^y)^2}{h}+\frac{gh^2}{2}\right) = -gh\dfrac{\partial B}{\partial y},\\
\end{aligned}
\end{equation}
where $h$ is the height of the water, $g$ is the gravitational constant, $B$ is the time-independent bottom topography function, and $q^x$ and $q^y$ are the discharges in the $x$- and $y$-directions, respectively. The quantities $h$, $q^x$, and $q^y$ depend on the spatial variables $(x,y)$ and time $t$, and $B$ depends only on $(x,y)$. We are interested in modeling a random/uncertain variant of this system.

Let $(\Omega, \mathcal{F}, P)$ be a complete probability space with event space $\Omega$, $\sigma$-algebra $\mathcal{F}$, and probability measure $P$ defined on $\mathcal{F}$. Uncertainty can enter the deterministic system, for example, via the bottom topography $B$. For example, the stochastic bottom can be modeled by a $d$-dimensional random field  
\begin{align*}
    B = B(x,y,\xi(\omega)) = B_0(x,y) + \sum_{k=1}^d B_k(x,y) \xi_k,
\end{align*}
where $\xi(\omega) = (\xi_1(\omega), \ldots, \xi_d(\omega))$ is a $d$-dimensional random variable for $\omega\in \Omega$. Such a stochastic bottom function will lead to a $\xi$-parameterized version of the 2D shallow water equations \eqref{eq:swesg1}, so that now $h = h(x, y, t, \xi)$, $q^x = q^x(x, y, t, \xi)$, $q^y = q^x(x,y,t,\xi)$, and $B = B(x, y, \xi)$. Our results can be generalized to other models of uncertainty, e.g., the uncertainty from the measurements of the initial water displacements $\eta\coloneqq h+B$. In the following context, we will refer to the quantities associated only with $(x,y)$ as ``physical'' quantities and the quantities associated only with $\xi$ as ``stochastic'' quantities.

\subsubsection{Polynomial Chaos Expansion - A Review}
In this section, we  briefly review the results and notation for polynomial chaos expansions (PCE), which are similar to what is used in \cite{doi:10.1137/20M1360736}. We include this section here in order to be self-contained; interested readers can find a more technical discussion of PCE in \cite{debusschere2004numerical,xiu2010numerical,sullivan2015introduction}.

Let $\rho$ be the Lebesgue density function for the random variable $\xi\in\mathbbm{R}^d$. The strategy of PCE is to approximate the dependence of a square-integrable random field by a polynomial of $\xi$. Let $\nu = (\nu_1, \ldots, \nu_d) \in \N_0^d$ be a multi-index. For a $d$-dimensional vector $\zeta \in \R^d$, the standard notation for the monomials with multi-index $\nu$ reads
\begin{align*}
  \zeta^\nu &\coloneqq \prod_{j=1}^d \zeta_j^{\nu_j}, & \zeta^0 = \zeta^{(0,0,\ldots,0)} &= 1.
\end{align*}
We consider the PCE approximations in the polynomial subspace defined by $\Lambda$:
\begin{align*}
  P_\Lambda &= \mathrm{span} \{ \zeta^\nu\;\; \big|\;\; \nu \in \Lambda\}, & \dim P_\Lambda &= K \coloneqq |\Lambda|,
\end{align*}
where $\Lambda\subset \N_0^d$ denotes any non-empty, size-$K$ finite set of multi-indices with $0 = (0,0,\ldots,0) \in \Lambda$.
We will also need ``powers" of this set, defined by $r$-fold products of $P_\Lambda$ elements:
{\footnotesize
\begin{align}\label{eq:PL}
 P_\Lambda^r &\coloneqq \mathrm{span} \left\{ \prod_{j=1}^r p_j \;\; \big|\;\; p_j \in P_\Lambda, \; j = 1,\ldots, r \right\}, & \dim P_\Lambda^r &\leq \left(\!\!\left(\begin{array}{c} K \\ r \end{array}\right)\!\!\right) = \left(\begin{array}{c} K+r-1 \\ r \end{array}\right),
\end{align}
}%
where the dimension bound results from a combinatorial argument. Note that since $0 \in \Lambda$, then $P^r_\Lambda \subseteq P^s_\Lambda$ for any $r \leq s$. We will later exercise the notation above for $r = 3$ (\Cref{thm:h-positivity}). Define the function space
\begin{align*}
    L^2_\rho(\R^d) \coloneqq \left\{f: \R^d\to\R\middle\vert\left(\int_{\R^d}f^2(s)\rho(s)ds\right)^{\frac{1}{2}}<+\infty\right\}.
\end{align*}
Assuming finite polynomial moments of all orders for $\rho$, 
there exists an $L^2_\rho(\R^d)$-orthonormal polynomial basis $\{\phi_k\}_{k=1}^\infty$ of $P_\Lambda$, i.e., 
\begin{align}\label{eq:orthocond}
  \langle\phi_k,\phi_\ell\rangle_{\rho} &\coloneqq \int_{\R^d} \phi_k(s) \phi_\ell(s) {\rho}(s)d s  = \delta_{k, \ell}, & \phi_1(\xi) &\equiv 1,
\end{align}
for all $k, \ell \in \N$, where $\delta_{k,\ell}$ is the Kronecker delta.

If $z(x,y,t,\cdot) \in L^2_\rho(\R^d)$, then under mild conditions on the probability measure $\rho$ (see \cite{ernst_convergence_2012}) there exists a convergent expansion of $z$ in these basis functions,
\begin{align*}
  z(x,y,t,\xi) &\stackrel{L^2_\rho}{=} \sum_{k=1}^\infty \widehat{z}_k(x,y,t) \phi_k(\xi), 
\end{align*}
where $\widehat{z}_{k}(x,y,t)$ are (deterministic) Fourier coefficients in the basis $\{\phi_k\}_{k \in \N}$. If $P_\Lambda = \mathrm{span}\{\phi_k\;\;\big|\;\; 1 \leq k \leq K\}$, then we define the $K$-term PCE \textit{approximation} of $z$ on $P_\Lambda$ to be the truncation of the first $K$ terms in the summation above,
\begin{align}\label{eq:PCEex} 
    z(x,y,t,\xi) \approx\sum_{k=1}^K \widehat{z}_k(x,y,t)\phi_k(\xi)\eqqcolon\Pi_\Lambda[z](x,y,t,\xi).
\end{align}
The operator $\Pi_\Lambda:L^2_\rho(\R^d)\to P_\Lambda$ is the linear projection operator onto $P_\Lambda$.

Orthogonality of the basis functions allows us to compute statistics of the the truncated PCE by manipulating the coefficients in a straightforward way. For example, the expectation and the variance of $\Pi_\Lambda[z](x,y,t,\xi)$ are given by 
\begin{equation}\label{eq:expvar}
  \mathbb{E}[\Pi_\Lambda[z](x,y,t,\xi)] = \widehat{z}_1(x,y,t),\quad \text{Var}[\Pi_\Lambda[z](x,y,t,\xi)] = \sum_{k=2}^{K}\widehat{z}_k^2(x,y,t),
\end{equation}
respectively.

We introduce the \emph{pseudo-spectral product} as the approximation to the $K$-term PCE of the product of two random fields $a(x,y,t,\xi)$ and $b(x,y,t,\xi)$,
\begin{align}\label{eq:mult-assump}
\Pi_\Lambda[a,b]
  &\coloneqq \Pi_\Lambda \left[ \Pi_\Lambda [a]\; \Pi_\Lambda[b] \right] 
  = \sum_{m=1}^{K}\left(\sum_{k,\ell=1}^K\widehat{a}_k\widehat{b}_\ell\langle\phi_k\phi_\ell,\phi_m\rangle_{\rho}\right)\phi_m(\xi).
\end{align}
This expression is an exact projection onto $P_\Lambda$ of the product of two $P_\Lambda$ projections. The approximation \eqref{eq:mult-assump} is a widely used strategy for computing PCE expansion of products (e.g. \cite{debusschere2004numerical,gerster2019hyperbolic}). 
The choice above is not the only approximation to $\Pi_\Lambda[ab]$, but is one that renders the 1D SG formulation to be hyperbolicity-preserving \cite{doi:10.1137/20M1360736}.
This choice is also a natural one since $\Pi_\Lambda\left[a,b\right] = \Pi_\Lambda\left[ab\right]$ if $(ab)\in P_\Lambda$. For general $a, b$, then $\Pi_\Lambda\left[a, b\right]$ and $\Pi_\Lambda\left[a b\right]$ are not equal since \eqref{eq:mult-assump} truncates the polynomial terms of $a$ and $b$ whose multi-indices are not in $\Lambda$. Such an operation can be cast in linear algebraic terms by considering vectors comprised of the PCE expansion coefficients. Let $\widehat{z} \in \R^K$ denote the vector of the $\phi_k$-expansion coefficients of $z \in P_\Lambda$. Define the linear operator $\mathcal{P}: \R^{K} \rightarrow \R^{K\times K}$ as
\begin{align}\label{eq:pmatrix}
  \mathcal{P}(\widehat{z}) &\coloneqq \sum_{k=1}^K\widehat{z}_k\mathcal{M}_k, &
  \mathcal{M}_k &\in \R^{K \times K}, &
  (\mathcal{M}_k)_{\ell m} &= \langle\phi_k,\phi_\ell\phi_m\rangle_{\rho},
\end{align}
where $\mathcal{M}_k$ is a symmetric matrix for each $k$.
The following properties can be verified:
\begin{align}\label{eq:pmatrixproperty}
  \mathcal{P}(\widehat{a}) &= \begin{pmatrix}\mathcal{M}_1\widehat{a}|\mathcal{M}_2\widehat{a}|\ldots |\mathcal{M}_K\widehat{a} \end{pmatrix}, & 
  \mathcal{P}(\widehat{a})\widehat{b} &= \mathcal{P}(\widehat{b})\widehat{a}, & 
  \widehat{\cG_\Lambda[a, b]} = \mathcal{P}(\widehat{a})\widehat{b}, 
\end{align}
where the last property is due to \eqref{eq:mult-assump}. If $\mathcal{P}(\widehat{a})$ is invertible, we further have 
\begin{equation}\label{eq:Pderivative}
\begin{aligned}
    &\frac{\partial (\mathcal{P}(\widehat{a})\mathcal{P}^{-1}(\widehat{a}))}{\partial \widehat{a}_\ell} = \frac{\partial I}{\partial \widehat{a}_\ell}= 0,\\
    \Rightarrow\;\;&\frac{\partial \mathcal{P}(\widehat{a})}{\partial \widehat{a}_\ell}\mathcal{P}^{-1}(\widehat{a}) 
    +  \mathcal{P}(\widehat{a})\frac{\partial \mathcal{P}^{-1}(\widehat{a})}{\partial \widehat{a}_\ell}= 0,\\
    \overset{\eqref{eq:pmatrix}\eqref{eq:pmatrixproperty}}{\Longrightarrow}\;\;& \frac{\partial \mathcal{P}^{-1}(\widehat{a})}{\partial \widehat{a}_\ell} = -\mathcal{P}^{-1}(\widehat{a})\mathcal{M}_\ell\mathcal{P}^{-1}(\widehat{a}).
\end{aligned}
\end{equation}
Equation \eqref{eq:Pderivative} will be exercised later for computing the Jacobian matrix of \eqref{eq:swesg4}.

We will also need to compute $P_\Lambda$ truncations of ratios of random fields (when for each $(x,y,t)$ the denominator is a single-signed variable with probability 1). For a single-signed random field $a$, we have
\begin{equation}\label{eq:prod-assump}
    \Pi_\Lambda\left[a\,\frac{b}{a}\right](x,y,t,\xi) =\Pi_\Lambda[b](x,y,t,\xi).
\end{equation}
Using \eqref{eq:prod-assump} as motivation, we make the assumption that the PCE of $b$ is equal to the pseudo-spectral product of $a$ and $b/a$,
\begin{align}\label{eq:div-assump}
    \Pi_\Lambda\left[a,\frac{b}{a}\right] = \Pi_\Lambda[b]
    \hskip 5pt \stackrel{\eqref{eq:pmatrixproperty}}{\Longleftrightarrow} \hskip 5pt
    \mathcal{P}(\widehat{a}) \widehat{\left(\frac{b}{a}\right)} = \widehat{b}.
\end{align}
In particular, if $b\in P_\Lambda$ (which is that case that we consider later), then implication \eqref{eq:div-assump} is true. Equation \eqref{eq:div-assump} then motivates the definition for a new operator for the $K$-term PCE of the ratios,
\begin{equation}\label{eq:div-pce}
    \Pi^{\dagger}_\Lambda\left[\frac{b}{a}\right](\xi) \coloneqq \sum_{k=1}^{K}c_k\phi_k(\xi),
\end{equation}
  where $c_k$ is the $k$th element of $\widehat{\left(\frac{b}{a}\right)}$ defined by \eqref{eq:div-assump}, 
  assuming $\mathcal{P}(\widehat{a})$ is invertible. \Cref{eq:div-pce} is a frequently used strategy for computing the PCE of the ratio of two polynomials \cite{debusschere2004numerical}. In the following section we will exercise \eqref{eq:div-pce} for the case where $b\in P_\Lambda$.
\subsection{Stochastic Galerkin Formulation for 2D SWE}
\begin{figure}[htbp]
    \centering
    \includegraphics[width=0.6\textwidth]{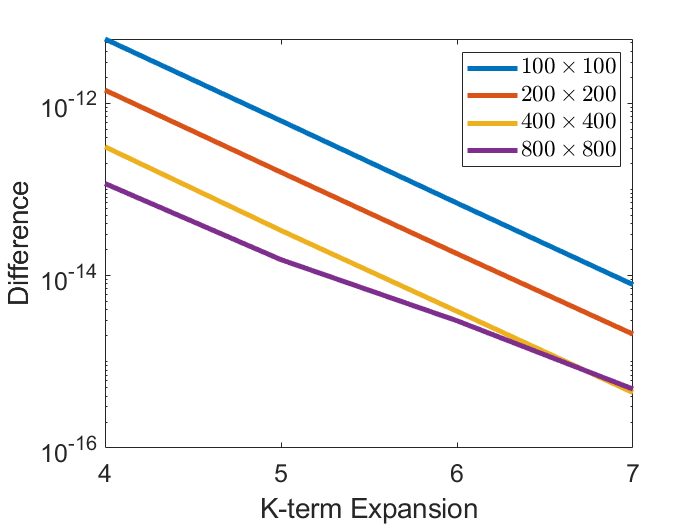}
    \caption{Numerical verification: the two approximations \eqref{eq:sg-assump-2}-\eqref{eq:sg-assump-3} are in general not the same. The difference is computed as the maximum of $L^2_\rho$ norm of the difference between the two PCE approximations of $q^xq^y/h$ in corresponding cells over the whole computational domain. The data are taken from the numerical solutions to \eqref{eq:initex1-1uq-conv}-\eqref{eq:bottomex1-1uq-conv} at $t = 0.07$ for different 2D grid sizes.}
    \label{fig:qxqyhdiff}    
\end{figure}

Following a standard Galerkin procedure in stochastic ($\xi$) space, we first reduce the problem to a finite-dimensional one by making the ansatz that the solutions to $h, q^x,$ and $q^y$ lie in the polynomial space $P_\Lambda$,
\begin{equation}\label{eq:sg-ansatz}
\begin{aligned}
  &h \simeq h_\Lambda \;\coloneqq \sum_{k=1}^K \hh_k(x,y,t) \phi_k(\xi), \\
  &q^x \simeq q^x_\Lambda \coloneqq \sum_{k=1}^K (\hqx)_k(x,y,t) \phi_k(\xi), \\
  &q^y \simeq q^y_\Lambda \coloneqq \sum_{k=1}^K (\hqy)_k(x,y,t) \phi_k(\xi),
\end{aligned}
\end{equation}
respectively, and $B$ is replaced by $\cG_\Lambda[B]$. We then apply the projection operator $\cG_\Lambda$ to both sides of \eqref{eq:swesg1} and enforce equality. Similar to \cite{doi:10.1137/20M1360736}, we make the assumption that,
\begin{equation}\label{eq:sg-assump-1}
\begin{aligned}
  \frac{(q^x)^2}{h} = \frac{q^x}{h}\; q^x \hskip 10pt \longrightarrow \hskip 10pt 
  \cG_\Lambda\left[ \frac{(q^x_\Lambda)^2}{h_\Lambda}\right] = \cG_\Lambda \left[  q^x_\Lambda\; \cG^\dagger_\Lambda\left[\frac{q^x_\Lambda}{h_\Lambda}\right]\right],\\
  \frac{(q^y)^2}{h} = \frac{q^y}{h}\; q^y \hskip 10pt \longrightarrow \hskip 10pt 
  \cG_\Lambda\left[ \frac{(q^y_\Lambda)^2}{h_\Lambda}\right] = \cG_\Lambda \left[  q^y_\Lambda\; \cG^\dagger_\Lambda\left[\frac{q^y_\Lambda}{h_\Lambda}\right]\right].\\
\end{aligned}
\end{equation}
In addition, we make another crucial assumption on the approximation of the nonlinear term $(q^xq^y)/h$ that occurs in the fluxes of both $x$- and $y$-directions. Specifically, for the nonlinear term $q^xq^y/h$ in $(q^xq^y/h)_x$, we use the approximation,
\begin{equation}\label{eq:sg-assump-2}
  \frac{q^xq^y}{h} = (q^x)\frac{q^y}{h}\; \hskip 10pt \longrightarrow \hskip 10pt 
  \cG_\Lambda\left[ \frac{q^x_\Lambda q^y_\Lambda}{h_\Lambda}\right] = \cG_\Lambda \left[  q^x_\Lambda\; \cG^\dagger_\Lambda\left[\frac{q^y_\Lambda}{h_\Lambda}\right]\right],
\end{equation}
while for the same nonlinear term in $(q^xq^y/h)_y$, we adopt the approximation,
\begin{equation}\label{eq:sg-assump-3}
  \frac{q^xq^y}{h} = (q^y)\frac{q^x}{h}\; \hskip 10pt \longrightarrow \hskip 10pt 
  \cG_\Lambda\left[ \frac{q^x_\Lambda q^y_\Lambda}{h_\Lambda}\right] = \cG_\Lambda \left[  q^y_\Lambda\; \cG^\dagger_\Lambda\left[\frac{q^x_\Lambda}{h_\Lambda}\right]\right].
\end{equation}
The two different approximations \eqref{eq:sg-assump-2}-\eqref{eq:sg-assump-3} turn out to be the crucial ingredient of our hyperbolicity-preserving SG formulation (\Cref{thm:hyperbolicity}). The two approximations are in general not equal. But as $K$ increases, we observe empirically that the discrepancy decreases exponentially (see \Cref{fig:qxqyhdiff}).

The steps introduced above lead to a new system of balance laws for the PCE vector of the stochastic state variables, 
\begin{equation}\label{eq:swesg4}
    \dfrac{\partial }{\partial t}(\hU)+\dfrac{\partial }{\partial x}(\widehat{F}(\hU))+\dfrac{\partial }{\partial y}(\widehat{G}(\hU)) = \widehat{S}(\hU,\widehat{B}).
\end{equation}
Here $\hU \coloneqq (\hh^\top, \hqx^\top, \hqy^\top)^\top$, where $\hh, \hqx$, and $\hqy$ are each length-$K$ vectors whose entries are the coefficients introduced in \eqref{eq:sg-ansatz}. The flux terms are, 
{\small
\begin{align}\label{eq:sgfluxes}
    &\widehat{F}(\hU) = \begin{pmatrix}
     \hqx\\\cP(\hqx)\cP^{-1}(\hh)\hqx+\frac{1}{2}g\cP(\hh)\hh\\\cP(\hqx)\cP^{-1}(\hh)\hqy
 \end{pmatrix},&&
 \widehat{G}(\hU) =  \begin{pmatrix}
     \hqy\\\cP(\hqy)\cP^{-1}(\hh)\hqx\\\cP(\hqy)\cP^{-1}(\hh)\hqy+\frac{1}{2}g\cP(\hh)\hh
 \end{pmatrix},
\end{align}
}
and the source term is
\begin{align}\label{eq:sgsource}
\widehat{S}(\hU,\widehat{B}) = \begin{pmatrix}0\\-g\cP(\hh)\widehat{B_x}\\-g\cP(\hh)\widehat{B_y}\end{pmatrix}.
\end{align}

Using \Cref{eq:pmatrix}-\eqref{eq:Pderivative} and following a similar calculation as in \cite{jin2019study}, we obtain the flux Jacobians written in terms of $K \times K$ blocks,
 \begin{align}\label{eq:x-jacobian}
   \frac{\partial \widehat{F}}{\partial \hU} &= \begin{pmatrix}O&I&O\\g\cP(\hh)-\cP(\hqx)\cP^{-1}(\hh)\cP(\hu)&\cP(\hqx)\cP^{-1}(\hh)+\cP(\hu)&O\\-\cP(\hqx)\cP^{-1}(\hh)\cP(\hv)&\cP(\hv)&\cP(\hqx)\cP^{-1}(\hh)\end{pmatrix}, \\
% \end{align}
% and 
% \begin{align}
     \label{eq:y-jacobian}
   \frac{\partial \widehat{G}}{\partial \hU} &= \begin{pmatrix}O&O&I\\-\cP(\hqy)\cP^{-1}(\hh)\cP(\hu)&\cP(\hqy)\cP^{-1}(\hh)&\cP(\hu)\\g\cP(\hh)-\cP(\hqy)\cP^{-1}(\hh)\cP(\hv)&O&\cP(\hqy)\cP^{-1}(\hh)+\cP(\hv)\end{pmatrix},
 \end{align}
 where we have introduced 
 \begin{align}\label{eq:uPCE}
    &\hu = \cP^{-1}(\hh)\hqx,&& \hv = \cP^{-1}(\hh)\hqy
\end{align} 
which can be interpreted as PCE coefficients of the $x$- and $y$-velocity, respectively.
\begin{rmk}
     In the numerical section \Cref{sec:scheme}, our scheme requires us to compute the maximum and minimum eigenvalues of both $\frac{\partial \widehat{F}}{\partial \hU}$ and $\frac{\partial \widehat{G}}{\partial \hU}$. Some algebraic manipulations can accelerate the computation. For example, the flux Jacobian $\frac{\partial \widehat{F}}{\partial \hU}$ has the same eigenvalues as the block diagonal matrix
     $$
        \begin{pmatrix}O&I&O\\g\cP(\hh)-\cP(\hqx)\cP^{-1}(\hh)\cP(\hu)&\cP(\hqx)\cP^{-1}(\hh)+\cP(\hu)&O\\O&O&\cP(\hqx)\cP^{-1}(\hh)\end{pmatrix}.
     $$
     One can first compute the maximum and minimum eigenvalues of the first $2K\times 2K$ submatrix and the $K\times K$ submatrix on the right corner, and then determine the maximum and minimum eigenvalues of the original matrix by comparing the eigenvalues of the submatrices. Similar calculations can be applied to $\frac{\partial \widehat{G}}{\partial \hU}$, where we can replace the $-\cP(\hqy)\cP^{-1}(\hh)\cP(\hu)$ and $\cP(\hu)$ blocks by zero matrices when computing eigenvalues.
\end{rmk}
\section{Hyperbolicity of the SG System}\label{sec:hyperbolicity}
In this section we will present a sufficient condition to ensure hyperbolicity of the system \eqref{eq:swesg4} as well as a computationally tractable way to preserve hyperbolicity of the system.
\begin{thm}\label{thm:hyperbolicity}
If the matrix $\mathcal{P}(\hh)$ is strictly positive definite for every point $(x,y,t)$ in the computational spatial-temporal domain, then the SG formulation \eqref{eq:swesg4} is hyperbolic.
\end{thm}
\begin{pf}
To show that the system \eqref{eq:swesg4} is hyperbolic, we need to show that, for any unit vector $\boldsymbol{n} = (n_x, n_y)^\top$, the matrix
\begin{align}
    n_x\frac{\partial \widehat{F}}{\partial \hU}+n_y\frac{\partial \widehat{G}}{\partial \hU}
\end{align}
is real diagonalizable. Define,
\begin{align}\label{eq:matrixnotations}
 &A \coloneqq g\mathcal{P}(\widehat{q^x}),
 &&B \coloneqq \mathcal{P}(\widehat{u}),
 &C \coloneqq g\mathcal{P}(\widehat{q^y}),
 &&D \coloneqq \mathcal{P}(\widehat{v}),
 &&E \coloneqq \sqrt{g\mathcal{P}(\hh)},
 \end{align}
where $\sqrt{M}$ is the (unique) symmetric positive definite square root of a symmetric positive definite matrix $M$. Define
 \begin{equation*}
     P=
     \begin{pmatrix}
     I&O&I\\B+n_xE&-n_yE&B-n_xE\\D+n_yE&n_xE&D-n_yE
     \end{pmatrix},
 \end{equation*}
whose inverse is given by 
\begin{equation*}
    P^{-1} = \left(\frac{1}{2}\right)(I\otimes E^{-1})
    \begin{pmatrix}
        E-n_xB-n_yD&n_xI&n_yI\\
        2(n_yB-n_xD)&-2n_yI&2n_xI\\
        E+n_xB+n_yD&-n_xI&-n_yI\\
    \end{pmatrix}.
\end{equation*}
Recall that $\otimes$ represents the Kronecker product between matrices.

%\begin{equation*}
%    P^{-1} = \left(\frac{1}{2}\right)
%    \begin{pmatrix}
%        I-n_xE^{-1}B-n_yE^{-1}D&n_xE^{-1}&n_yE^{-1}\\
%        2(n_yE^{-1}B-n_xE^{-1}D)&-2n_yE^{-1}&2n_xE^{-1}\\
%        I+n_xE^{-1}B+n_yE^{-1}D&-n_xE^{-1}&-n_yE^{-1}\\
%    \end{pmatrix}.
%\end{equation*}
Using \eqref{eq:x-jacobian}-\eqref{eq:y-jacobian}, direct calculation leads to
\begin{equation}\label{eq:dfdu-similar}
    J\coloneqq\begin{aligned}
    P^{-1}\left(n_x\frac{\partial \widehat{F}}{\partial \hU}+n_y\frac{\partial \widehat{G}}{\partial \hU}\right)P
    =&\left(\frac{1}{2}\right)\begin{pmatrix}
        J_{11}&O&J_{13}\\
        O&J_{22}&O\\
        J_{31}&O&J_{33},
    \end{pmatrix}.
\end{aligned}
\end{equation}
where,
\begin{equation}\label{eq:submatrices}
    \begin{aligned}
        &J_{11} = 2E+n_x(B+E^{-1}AE^{-1})+n_y(D+E^{-1}CE^{-1}),\\
        &J_{13}= n_x(B-E^{-1}AE^{-1})+n_y(D-E^{-1}CE^{-1}) = J_{31},\\
        &J_{22} = 2n_xE^{-1}AE^{-1}+2n_yE^{-1}CE^{-1},\\
        &J_{33} = -2E+n_x(B+E^{-1}AE^{-1})+n_y(D+E^{-1}CE^{-1}).
    \end{aligned}
\end{equation}
From \eqref{eq:matrixnotations}, the matrices $A, B, C,D$ and $E$ are all symmetric. Thus, it can be verified that $J_{11}, J_{13}, J_{22}, J_{31},$ and $J_{33}$ are all symmetric. Since $J_{13} = J_{31}$, the matrix $J$ is symmetric and therefore the similar matrix $\left(n_x\frac{\partial \widehat{F}}{\partial \hU}+n_y\frac{\partial \widehat{G}}{\partial \hU}\right)$ is real diagonalizable.\qed
\end{pf}
\begin{rmk}
     When the system is deterministic, the matrix $P$ and the matrix $J$ in the proof above reduce to the eigenmatrix that symmetrizes the deterministic Jacobian matrix and the diagonal matrix after symmetrization, respectively.
\end{rmk}
From \Cref{thm:hyperbolicity}, the hyperbolicity of the SG system relies only on the positive definiteness of the matrix $\mathcal{P}(\hh)$, a stochastic variant of the hyperbolic condition for the deterministic 2D SWE. In \cite{doi:10.1137/20M1360736}, a computationally convenient sufficient condition to guarantee $\mathcal{P}(\hh) > 0$ is given by the following theorem.
\begin{thm}\label{thm:h-positivity}
 Given $\Lambda$, let nodes $\xi_m$ and weights $\tau_m$ satisfying $\{(\xi_m, \tau_m)\}_{m=1}^M \subset \R^d \times (0, \infty)$ represent any $M$-point positive quadrature rule that is exact on $P_{\Lambda}^3$, i.e., 
  \begin{align}\label{eq:P3-exactness}
    \int_{\R^d} p(\xi) \rho(\xi) d \xi &= \sum_{m=1}^M p(\xi_m) \tau_m, & p &\in P_\Lambda^3.
  \end{align}
  If 
  \begin{align}\label{eq:h-positivity}
    h_\Lambda(x,y,t,\xi_m) > 0 \;\; \forall\; m = 1, \ldots, M,
  \end{align}
  then the SG system of the 2D SWE \eqref{eq:swesg4} is hyperbolic.
\end{thm}
\begin{pf}
    See Theorem 3.4 in \cite{doi:10.1137/20M1360736}.\qed
\end{pf}
Thus, to preserve the hyperbolicity of \eqref{eq:swesg4}, we only need to preserve a \textit{finite} number of positivity conditions on the stochastic water height. In the worst-case scenario, the number $M$ of quadrature points used in the stochastic space is of order $O(\text{dim} P_\Lambda^3) \lesssim K^3$ (see \eqref{eq:PL} for the definition of $P_\Lambda^3$). For one-dimensional stochastic space, one can use the $M$-point Gaussian quadrature associated with the density $\rho$ of $\xi$, with $M$ any integer satisfying $M\ge\left\lceil \frac{3 K}{2} \right\rceil - 1$.

Instead of preserving the positive definiteness of $\mathcal{P}(\hh)$, one could choose to preserve the positivity of the expectation $\hh_1$ of $h_\Lambda$. The following corollary states that this condition is guaranteed by $\mathcal{P}(\hh)>0$.
\begin{corollary}
  %If the conditions of \cref{thm:h-positivity} are satisfied, then the expectation $\hh_1$ of $h_\Lambda$ is positive.
If the assumptions of Theorem \eqref{thm:h-positivity} hold, then $\hh_1>0$.
\end{corollary}
\begin{pf}
    See Corollary 3.6 in \cite{doi:10.1137/20M1360736}.\qed
\end{pf}
%%%%%%%%%%%%%%%%%%%%%%%%%%%%%%%%%%%%%%%%%%%%%%%%%%%%%%%%%%%%%%%%%%%%%%%%%%%%%
%%% 5. Numerical scheme
%%%%%%%%%%%%%%%%%%%%%%%%%%%%%%%%%%%%%%%%%%%%%%%%%%%%%%%%%%%%%%%%%%%%%%%%%%%%%
\section{Numerical Scheme for Stochastic Shallow Water Equations}\label{sec:scheme}
The family of central-upwind schemes (\cite{kurganov2000new, kurganov2001semidiscrete}, etc.) are Riemann-solver-free, second-order, high-resolution MUSCL schemes \cite{van1974towards}. The generalized minmod limiter is used for the piecewise linear approximation of the solution. In this section, we present a well-balanced central-upwind scheme that preserves the hyperbolicity of the system \eqref{eq:swesg4} at every time step. 
\subsection{Central-Upwind Scheme for the SG System} 
For simplicity, we consider only rectangular spatial domains. Let $\left\{\mathcal{C}_{i,j}\right\}_{i\in [N_x], j\in [N_y]}$ be a uniform rectangular partition over the rectangular computational domain, where
\begin{align*}
    \mathcal{C}_{i,j} \coloneqq \left[ x_{i-\frac{1}{2}}, x_{i+\frac{1}{2}}\right]\times\left[ y_{j-\frac{1}{2}}, y_{j+\frac{1}{2}}\right],
\end{align*}
$[N_x]$ is the index set $\{1,2,\ldots,N_x\}$, and $[N_y]$ is defined similarly. Define $\Delta x \coloneqq x_{i+\frac{1}{2}}-x_{i-\frac{1}{2}}, \Delta y = y_{j+\frac{1}{2}}-y_{j-\frac{1}{2}}$. Let $\vert \mathcal{C}_{i,j}\vert$ denote the size $\Delta x\Delta y$ of cell $(i,j)$. Define the approximation to the cell average of the vector $\hU$ over cell $\mathcal{C}_{i,j}$ to be 
\begin{align*}
    \bU_{i,j}(t) \coloneqq
    \begin{pmatrix}
    \bh_{i,j}(t)\\ \bqx_{i,j}(t)\\ \bqy_{i,j}(t)
    \end{pmatrix}
    \approx
    \frac{1}{\vert \mathcal{C}_{i,j}\vert}\bigintsss_{\mathcal{C}_{i,j}}
    \begin{pmatrix}
    \hh(x,y,t)\\
    \hqx(x,y,t)\\
    \hqy(x,y,t)
    \end{pmatrix}
    dx\;dy\in \mathbbm{R}^{3K}.
\end{align*}
Since $\hU$ with subscripts has been used for the polynomial moments, we use the bold letter $\boldsymbol{U}$ with superscripts and subscripts here and in the following subsections to represent the cell averages and the pointwise values of the reconstructions, respectively, for the vector function $\hU(x,y,t)$. For example, $\boldsymbol{U}^{W}_{i,j}$ is the reconstructed value of $\hU$ at the west-side $(x_{i-\frac{1}{2}}, y_{j})$ of the cell $\mathcal{C}_{i,j}$. A similar reasoning applies to the pairs $(\boldsymbol{h}, \hh)$, $(\boldsymbol{q}^x, \hqx)$, and $(\boldsymbol{q}^y, \hqy)$. To minimize clutter, we will notationally suppress $t$ dependence from here onward.

For a continuous bottom topography, the pointwise values of the reconstructions are given by
\begin{align*}
    &\boldsymbol{U}^{E}_{i,j} = \bU_{i,j}+\frac{\Delta x}{2}(\boldsymbol{U}_x)_{i,j},
    &&\boldsymbol{U}^{W}_{i,j} = \bU_{i,j}-\frac{\Delta x}{2}(\boldsymbol{U}_x)_{i,j},\\
    &\boldsymbol{U}^{N}_{i,j} = \bU_{i,j}+\frac{\Delta y}{2}(\boldsymbol{U}_y)_{i,j},
    &&\boldsymbol{U}^{S}_{i,j} = \bU_{i,j}-\frac{\Delta y}{2}(\boldsymbol{U}_y)_{i,j},\\    
\end{align*}
where $\boldsymbol{U}_x$ and $\boldsymbol{U}_y$ are computed using the generalized minmod limiter,
\begin{align*}
    &(\boldsymbol{U}_x)_{i,j} = \text{minmod}\left(\theta\frac{\bU_{i,j}-\bU_{i-1,j}}{\Delta x},\frac{\bU_{i+1,j}-\bU_{i-1,j}}{2\Delta x},\theta\frac{\bU_{i+1,j}-\bU_{i,j}}{\Delta x}\right),\\
    &(\boldsymbol{U}_y)_{i,j} = \text{minmod}\left(\theta\frac{\bU_{i,j}-\bU_{i,j-1}}{\Delta y},\frac{\bU_{i,j+1}-\bU_{i,j-1}}{2\Delta y},\theta\frac{\bU_{i,j+1}-\bU_{i,j}}{\Delta y}\right).  
\end{align*}
The minmod function for scalars is defined by
\begin{equation*}
  \text{minmod}(c_1,c_2,\ldots, c_n) \coloneqq \left\{\begin{aligned}
       &\max\nolimits_{j} \{c_j\},&& \text{if } c_j<0\;\;\forall j,\\
       &\min\nolimits_{j} \{c_j\},&& \text{if } c_j>0\;\;\forall j,\\
       &0,&&  \text{otherwise},\\
  \end{aligned}\right.
\end{equation*}
and the minmod function for vector is simply applying the scalar minmod function for each component. The parameter $\theta\in [1,2]$ controls how dissipative the limiter is. The limiter is less dissipative for larger value of $\theta$, but it remains \textit{non-oscillatory} for any choice of $\theta$ \cite{kurganov2001semidiscrete}.

The speeds of local propagation are then estimated by, 
\begin{equation}\label{eq:pspeed}
    \begin{aligned}
        &a^{-}_{i+\frac{1}{2},j} = \min\left\{\lambda_1\left(\frac{\partial \widehat{F}}{\partial \hU}(\boldsymbol{U}^{W}_{i+1,j})\right),\lambda_1\left(\frac{\partial \widehat{F}}{\partial \hU}(\boldsymbol{U}^{E}_{i,j})\right),0\right\},\\
        &a^{+}_{i+\frac{1}{2},j} =  \max\left\{\lambda_{3K}\left(\frac{\partial \widehat{F}}{\partial \hU}(\boldsymbol{U}^{W}_{i+1,j})\right),\lambda_{3K}\left(\frac{\partial \widehat{F}}{\partial \hU}(\boldsymbol{U}^{E}_{i,j})\right),0\right\},\\
        &b^{-}_{i,j+\frac{1}{2}} = \min\left\{\lambda_1\left(\frac{\partial \widehat{G}}{\partial \hU}(\boldsymbol{U}^{S}_{i,j+1})\right),\lambda_1\left(\frac{\partial \widehat{G}}{\partial \hU}(\boldsymbol{U}^{N}_{i,j})\right),0\right\},\\
        &b^{+}_{i,j+\frac{1}{2}} =  \max\left\{\lambda_{3K}\left(\frac{\partial \widehat{G}}{\partial \hU}(\boldsymbol{U}^{S}_{i,j+1})\right),\lambda_{3K}\left(\frac{\partial \widehat{G}}{\partial \hU}(\boldsymbol{U}^{N}_{i,j})\right),0\right\},  
    \end{aligned}
\end{equation}
where $\lambda_1$ and $\lambda_{3K}$ are the smallest and the largest eigenvalues of the matrix, respectively. The semi-discrete version of the central-upwind scheme is given by
\begin{align}\label{eq:semidiscretewsg}
  \dfrac{d}{dt}\bU_{i,j} = -\dfrac{\mathcal{F}_{i+\frac{1}{2},j}-\mathcal{F}_{i-\frac{1}{2},j}}{\Delta x}-\dfrac{\mathcal{G}_{i,j+\frac{1}{2}}-\mathcal{G}_{i,j-\frac{1}{2}}}{\Delta y}+\overline{\boldsymbol{S}}_{i,j},
\end{align}
where $\overline{\boldsymbol{S}}_{i,j}\approx\frac{1}{\vert \mathcal{C}_{i,j}\vert}\int_{\mathcal{C}_{i,j}}\widehat{S}(\boldsymbol{U},\boldsymbol{B})dxdy$
is the well-balanced discretization of the source term to be discussed in the following section. The numerical fluxes for $x$- and $y$-direction are given by 
\begin{equation}\label{eq:x-flux}
    \mathcal{F}_{i+\frac{1}{2},j} \coloneqq \dfrac{a^{+}_{i+\frac{1}{2},j}\widehat{F}(\boldsymbol{U}^{E}_{i,j})-a^{-}_{i+\frac{1}{2},j}\widehat{F}(\boldsymbol{U}^{W}_{i+1,j})}{a^{+}_{i+\frac{1}{2},j}-a^{-}_{i+\frac{1}{2},j}}+\dfrac{a^{+}_{i+\frac{1}{2},j}a^{-}_{i+\frac{1}{2},j}}{a^{+}_{i+\frac{1}{2},j}-a^{-}_{i+\frac{1}{2},j}}\left[\boldsymbol{U}^{W}_{i+1,j}-\boldsymbol{U}^{E}_{i,j}\right],
\end{equation}
and
\begin{equation}\label{eq:y-flux}
    \mathcal{G}_{i+\frac{1}{2},j} \coloneqq \dfrac{b^{+}_{i,j+\frac{1}{2}}\widehat{G}(\boldsymbol{U}^{N}_{i,j})-b^{-}_{i,j+\frac{1}{2}}\widehat{G}(\boldsymbol{U}^{S}_{i,j+1})}{b^{+}_{i,j+\frac{1}{2}}-b^{-}_{i,j+\frac{1}{2}}}+\dfrac{b^{+}_{i,j+\frac{1}{2}}b^{-}_{i,j+\frac{1}{2}}}{b^{+}_{i,j+\frac{1}{2}}-b^{-}_{i,j+\frac{1}{2}}}\left[\boldsymbol{U}^{S}_{i,j+1}-\boldsymbol{U}^{N}_{i,j}\right],
\end{equation}
respectively.
\subsection{Well-Balanced Flux Discretization}
For deterministic 2D SWE, the numerical scheme should accurately capture the ``lake-at-rest'' steady state solution, or small perturbations to this state on coarse meshes. The scheme is said to be \textit{well-balanced} if it possesses such a property. A stochastic variant to the ``lake-at-rest'' property is given by
\begin{equation}\label{eq:lake-at-rest}
  q^x_\Lambda(x,y,t,\xi) = q^y_\Lambda(x,y,t,\xi)\equiv 0,\quad h_\Lambda + \Pi_\Lambda[B](x,y,t,\xi) \equiv C(\xi),
\end{equation}
where $C(\xi)$ is a random variable depending only on $\xi$. The solution to \eqref{eq:lake-at-rest} represents a still water surface with a flat but stochastic water surface. The corresponding vector equations for the PCE coefficients are
\begin{equation}
\hqx = \hqy \equiv \boldsymbol{0},\quad\hh+\widehat{B} \equiv \widehat{C}.
\end{equation}

Within the framework of central-upwind scheme, the bottom function $\widehat{B}$ is first replaced by its continuous piecewise bilinear interpolant $\widetilde{\widehat{B}}$,
\begin{equation}\label{eq:reconstructedbot}
\begin{aligned}
    \widetilde{\widehat{B}}(x,y) &= \widehat{B}_{i-\frac{1}{2},j-\frac{1}{2}} + \left(\widehat{B}_{i+\frac{1}{2},j-\frac{1}{2}}-\widehat{B}_{i-\frac{1}{2},j-\frac{1}{2}}\right)\cdot\frac{x-x_{i-\frac{1}{2}}}{\Delta x}\\
    &\hspace{5.5em}+\left(\widehat{B}_{i-\frac{1}{2},j+\frac{1}{2}}-\widehat{B}_{i-\frac{1}{2},j-\frac{1}{2}}\right)\cdot\frac{y-y_{j-\frac{1}{2}}}{\Delta y}\\
    &+\left(\widehat{B}_{i+\frac{1}{2},j+\frac{1}{2}} - \widehat{B}_{i+\frac{1}{2},j-\frac{1}{2}} - \widehat{B}_{i-\frac{1}{2},j+\frac{1}{2}} +
    \widehat{B}_{i-\frac{1}{2},j-\frac{1}{2}}\right)\\
    &\hspace{5.5em}\cdot \frac{(x-x_{i-\frac{1}{2}})(y-y_{j-\frac{1}{2}})}{\Delta x\Delta y}, \quad(x,y)\in \mathcal{C}_{i,j}.
\end{aligned}
\end{equation}
For a general (possibly discontinuous) bottom topography function, the values at the corners of the cells are estimated by the average of the lower and the upper limits,
\begin{align}\label{eq:bcorners}
    \widehat{B}_{i\pm\frac{1}{2},j\pm\frac{1}{2}} = \frac{1}{2}\left(\liminf_{(x,y)\to (x_{i\pm\frac{1}{2}},y_{j\pm\frac{1}{2}})}\widehat{B}(x,y) + \limsup_{(x,y)\to (x_{i\pm\frac{1}{2}},y_{j\pm\frac{1}{2}})}\widehat{B}(x,y)\right).
\end{align}
For continuous $\widehat{B}(x,y)$, \eqref{eq:bcorners} reduces to
\begin{align*}
    \widehat{B}_{i\pm\frac{1}{2},j\pm\frac{1}{2}} = \widehat{B}(x_{i\pm\frac{1}{2}},y_{j\pm\frac{1}{2}}).
\end{align*}
The cell average of $\widetilde{\widehat{B}}$ in the cell $\mathcal{C}_{i,j}$ is equal to the averages of the values of $\widetilde{\widehat{B}}$ at centers of the edges of the cell, i.e.,
\begin{equation*}
\begin{aligned}
\overline{\boldsymbol{B}}_{i,j}&\coloneqq \frac{1}{\vert\mathcal{C}_{i,j}\vert}\iint_{\mathcal{C}_{i,j}}\widetilde{\widehat{B}}(x,y)dxdy,\\
&=\frac{1}{4}(\boldsymbol{B}_{i+\frac{1}{2},j}+\boldsymbol{B}_{i-\frac{1}{2},j}+\boldsymbol{B}_{i,j+\frac{1}{2}}+\boldsymbol{B}_{i,j-\frac{1}{2}}),
\end{aligned}
\end{equation*}
where 
\begin{align}\label{eq:piecewisebottom}
    &\boldsymbol{B}_{i+\frac{1}{2},j}\coloneqq \widetilde{\widehat{B}}(x_{i+\frac{1}{2}}, y_j),&&\boldsymbol{B}_{i,j+\frac{1}{2}}\coloneqq \widetilde{\widehat{B}}(x_{i}, y_{j+\frac{1}{2}}).
\end{align}
Note that, since $\widetilde{\widehat{B}}$ is a piecewise continuous function, the cell average $\overline{\boldsymbol{B}}_{i,j}$ is equal to the value of the interpolant at the center of the cell
$\widetilde{\widehat{B}}(x_i, y_j)$.

At every time step, we compute the cell averages of the PCE vector of the water surface $\overline{\boldsymbol{\eta}}_{i,j}\coloneqq \overline{\boldsymbol{h}}_{i,j} + \overline{\boldsymbol{B}}_{i,j}$ and their second-order-accurate pointwise values $\boldsymbol{\eta}^{W,E,S,N}_{i,j}$ of the reconstructions. The pointwise values of the reconstructions for the water heights are then computed by
\begin{equation}\label{eq:hreconstruction}
\begin{aligned}
    &\boldsymbol{h}^{W}_{i,j}\coloneqq \boldsymbol{\eta}^{W}_{i,j}-\boldsymbol{B}_{i+\frac{1}{2},j}, 
    &&\boldsymbol{h}^{E}_{i,j}\coloneqq \boldsymbol{\eta}^{E}_{i,j}-\boldsymbol{B}_{i-\frac{1}{2},j},\\
    &\boldsymbol{h}^{S}_{i,j}\coloneqq \boldsymbol{\eta}^{S}_{i,j}-\boldsymbol{B}_{i,j-\frac{1}{2}},
    &&\boldsymbol{h}^{N}_{i,j}\coloneqq \boldsymbol{\eta}^{N}_{i,j}-\boldsymbol{B}_{i,j+\frac{1}{2}}.
\end{aligned}    
\end{equation}
The numerical fluxes $\left\{\left(\mathcal{F}_{i+\frac{1}{2},j},\mathcal{G}_{i,j+\frac{1}{2}}\right)\right\}_{i\in[N_x],j\in [N_y]}$ are then computed using the heights defined in \eqref{eq:hreconstruction}. With the reconstructions for the heights and the fluxes, we are now ready to propose the well-balanced discretization for the source term.
\begin{proposition}
The central-upwind scheme \eqref{eq:semidiscretewsg} is well-balanced if we choose $\overline{\boldsymbol{S}}_{ij} = \left((\overline{\boldsymbol{S}}^{(1)}_{i,j})^\top,(\overline{\boldsymbol{S}}^{(2)}_{i,j})^\top,(\overline{\boldsymbol{S}}^{(3)}_{i,j})^\top\right)^\top$, where
  \begin{equation}\label{eq:wbproperty}
  \left\{
    \begin{aligned}
      \overline{\boldsymbol{S}}^{(1)}_{i,j} &= \boldsymbol{0},\\
      \overline{\boldsymbol{S}}^{(2)}_{i,j} &=  -g\mathcal{P}(\overline{\boldsymbol{h}}_{i,j})\left(\dfrac{\boldsymbol{B}_{i+\frac{1}{2},j}-\boldsymbol{B}_{i-\frac{1}{2},j}}{\Delta x}\right),\\
      \overline{\boldsymbol{S}}^{(3)}_{i,j} &=  -g\mathcal{P}(\overline{\boldsymbol{h}}_{i,j})\left(\dfrac{\boldsymbol{B}_{i,j+\frac{1}{2}}-\boldsymbol{B}_{i,j-\frac{1}{2}}}{\Delta y}\right).      
    \end{aligned}\right.
\end{equation}
\end{proposition}
\begin{pf}
    Let $\boldsymbol{\eta}^{W,E,S,N}_{i,j}$ be the pointwise values of the reconstructions for the water surface at cell $\mathcal{C}_{i,j}$. Assume at time $t$, the stochastic water surface is flat and the water is still, that is, 
    \begin{align*}
     &\overline{\boldsymbol{\eta}}_{i,j} = \boldsymbol{\eta}^{*} = \text{const vector},   &&\overline{\boldsymbol{q}^x}_{i,j} =\overline{\boldsymbol{q}^y}_{i,j} = \boldsymbol{0},
    \end{align*}
     for all cell $\mathcal{C}_{i,j}$. The piecewise linear reconstructions yield $\boldsymbol{\eta}^{W,E,S,N}_{i,j} = \boldsymbol{\eta}^{*}$ and $(\boldsymbol{q}^x)^{W,E,S,N}_{i,j} = \boldsymbol{0}$. Hence, the numerical fluxes $\left\{\left(\mathcal{F}_{i+\frac{1}{2},j},\mathcal{G}_{i,j+\frac{1}{2}}\right)\right\}_{i\in[N_x],j\in [N_y]}$ are 
    \begin{equation}\label{eq:sgxflux}
        \mathcal{F}_{i+\frac{1}{2},j} = \begin{pmatrix}\boldsymbol{0}\\\frac{g}{2}\mathcal{P}(\boldsymbol{\eta}^*-\boldsymbol{B}_{i+\frac{1}{2},j})(\boldsymbol{\eta}^*-\boldsymbol{B}_{i+\frac{1}{2},j})\\\boldsymbol{0}\end{pmatrix},
    \end{equation}
    and
    \begin{equation}\label{eq:sgyflux}
        \mathcal{G}_{i,j+\frac{1}{2}} = \begin{pmatrix}\boldsymbol{0}\\\boldsymbol{0}\\\frac{g}{2}\mathcal{P}(\boldsymbol{\eta}^*-\boldsymbol{B}_{i,j+\frac{1}{2}})(\boldsymbol{\eta}^*-\boldsymbol{B}_{i,j+\frac{1}{2}})\end{pmatrix}.
    \end{equation}
    The semidiscrete form for the stochastic ``lake-at-rest'' problem is thus given by
    {\small 
    \begin{equation}\label{eq:wbsemi}
        \begin{aligned}
            \dfrac{d }{dt}\overline{\boldsymbol{h}}_{i,j} &= \overline{\boldsymbol{S}}^{(1)}_{i,j}  \stackrel{\eqref{eq:wbproperty}}{=}\boldsymbol{0},\\ 
            \dfrac{d }{dt}\overline{\boldsymbol{q}^x}_{i,j} &= 
            -\dfrac{1}{\Delta x}\frac{g}{2}\left[\mathcal{P}(\boldsymbol{\eta}^*-\boldsymbol{B}_{i+\frac{1}{2},j})(\boldsymbol{\eta}^*-\boldsymbol{B}_{i+\frac{1}{2},j})-\mathcal{P}(\boldsymbol{\eta}^*-\boldsymbol{B}_{i-\frac{1}{2},j})(\boldsymbol{\eta}^*-\boldsymbol{B}_{i-\frac{1}{2},j})\right]+\overline{\boldsymbol{S}}^{(2)}_{i,j},\\
            \dfrac{d }{dt}\overline{\boldsymbol{q}^y}_{i,j} &= 
            -\dfrac{1}{\Delta y}\frac{g}{2}\left[\mathcal{P}(\boldsymbol{\eta}^*-\boldsymbol{B}_{i,j+\frac{1}{2}})(\boldsymbol{\eta}^*-\boldsymbol{B}_{i,j+\frac{1}{2}})-\mathcal{P}(\boldsymbol{\eta}^*-\boldsymbol{B}_{i,j-\frac{1}{2}})(\boldsymbol{\eta}^*-\boldsymbol{B}_{i,j-\frac{1}{2}})\right]+\overline{\boldsymbol{S}}^{(3)}_{i,j},     
        \end{aligned}
    \end{equation}
    }  
    The well-balanced property requires that the numerical fluxes to be balanced by the source term, i.e., the right-hand-side of \eqref{eq:wbsemi} needs to vanish. Clearly, $\overline{\boldsymbol{S}}^{(1)}_{i,j} = \boldsymbol{0}$. Using the linearity of the operator $\mathcal{P}$ and the property \eqref{eq:pmatrixproperty}, one can show that if we choose $\overline{\boldsymbol{S}}^{(2)}_{i,j}$ and $\overline{\boldsymbol{S}}^{(3)}_{i,j}$ by \eqref{eq:wbproperty}, the right-hand-side of \eqref{eq:wbsemi} will vanish and the well-balanced property of the scheme can be satisfied (see also \cite{doi:10.1137/20M1360736}, Lemma 4.1).\qed
\end{pf}
    The numerical approximations in \eqref{eq:wbproperty} can be treated as the approximation of the cell average of the source term using midpoint quadrature rule with the partial derivatives $\boldsymbol{B}_x(x_i,y_j), \boldsymbol{B}_y(x_i,y_j)$ approximated by the central difference.
\subsection{Hyperbolicity-Preserving Time Discretization}
The hyperbolicity of the system \eqref{eq:swesg4} depends only on the vector $\hh$ and so we focus on the first $K$ equations in \eqref{eq:semidiscretewsg} in this subsection,
\begin{align}\label{eq:semidiscreteh}
  \dfrac{d}{dt}\bh_{i,j} = -\dfrac{\mathcal{F}^{\hh}_{i+\frac{1}{2},j}-\mathcal{F}^{\hh}_{i-\frac{1}{2},j}}{\Delta x}-\dfrac{\mathcal{G}^{\hh}_{i,j+\frac{1}{2}}-\mathcal{G}^{\hh}_{i,j-\frac{1}{2}}}{\Delta y},
\end{align}
where
\begin{equation}\label{eq:xfluxcu}
    \mathcal{F}^{\hh}_{i+\frac{1}{2},j} = \dfrac{a^{+}_{i+\frac{1}{2},j}(\boldsymbol{q}^x)_{i,j}^{E}-a^{-}_{i+\frac{1}{2},j}(\boldsymbol{q}^x)_{i+1,j}^{W}}{a^{+}_{i+\frac{1}{2},j}-a^{-}_{i+\frac{1}{2},j}}+\dfrac{a^{+}_{i+\frac{1}{2},j}a^{-}_{i+\frac{1}{2},j}}{a^{+}_{i+\frac{1}{2},j}-a^{-}_{i+\frac{1}{2},j}}\left[\boldsymbol{h}_{i+1,j}^{W}-\boldsymbol{h}_{i,j}^{E}\right],
\end{equation}
and
\begin{equation}\label{eq:yfluxcu}
    \mathcal{G}^{\hh}_{i,j+\frac{1}{2}} = \dfrac{b^{+}_{i,j+\frac{1}{2}}(\boldsymbol{q}^y)_{i,j}^{N}-b^{-}_{i,j+\frac{1}{2}}(\boldsymbol{q}^y)_{i,j+1}^{S}}{b^{+}_{i,j+\frac{1}{2}}-b^{-}_{i,j+\frac{1}{2}}}+\dfrac{b^{+}_{i,j+\frac{1}{2}}b^{-}_{i,j+\frac{1}{2}}}{b^{+}_{i,j+\frac{1}{2}}-b^{-}_{i,j+\frac{1}{2}}}\left[\boldsymbol{h}_{i,j+1}^{S}-\boldsymbol{h}_{i,j}^{N}\right].
\end{equation}
The following lemma presents a computable restriction on a hyperbolicity-preserving time step for a forward Euler time discretization.
\begin{lem}\label{lemma:CFL}
  Let $\left\{\xi_m\right\}_{m\in [M]}$ be the nodes of a quadrature rule satisfying the conditions of \Cref{thm:h-positivity} and $\boldsymbol{\Phi}(\xi) = (\phi_1(\xi),\ldots,\phi_K(\xi))^\top$. Denote the numerical approximation to $\bh_{i,j}(t^n)$ by $\bh_{i,j}^n$ and $\Delta t^n \coloneqq t^{n+1} - t^n$.
  Assume that $\overline{\boldsymbol{h}}_{i,j}^n(\xi_m)\coloneqq (\overline{\boldsymbol{h}}_{i,j}^n)^\top\Phi(\xi_m)\;  >  \;0$ for $m\in [M]$. If $\Delta t^n$ satisfies
  {\small
  \begin{align}\label{eq:CFL}
    \Delta t^n &< \Delta t^n_{h} \coloneqq 
    \min_{\substack{m\in[M]\\i\in[N_x],j\in [N_y]}}\left\{ \left|\dfrac{(\overline{\boldsymbol{h}}^{n}_{i,j})^{\top}\boldsymbol{\Phi}(\xi_m)}{\left[\frac{\mathcal{F}^{\hh}_{i+\frac{1}{2},j}(t_n)-\mathcal{F}^{\hh}_{i-\frac{1}{2},j}(t_n)}{\Delta x}+\frac{\mathcal{G}^{\hh}_{i,j+\frac{1}{2}}(t_n)-\mathcal{G}^{\hh}_{i,j-\frac{1}{2}}(t_n)}{\Delta y}\right]^{\top}\boldsymbol{\Phi}(\xi_m)}\right|\right\},
  \end{align}
  }
  then the flux Jacobian \eqref{eq:dfdu-similar}, $J\left(\overline{\boldsymbol{U}}_{i,j}^{n+1}\right)$ is diagonalizable with real eigenvalues when \eqref{eq:swesg4} is discretized using forward Euler solver. 
\end{lem}
\begin{pf}
    From \Cref{thm:h-positivity}, we only need to guarantee that $(\overline{\boldsymbol{h}}_{i,j}^{n+1})^\top\boldsymbol{\Phi}(\xi_m)>0$ for all $m\in [M]$, which, for a forward Euler discretization, gives the inequality
    {\small
    \begin{align}\label{eq:dtpositivity}
         0 <& (\overline{\boldsymbol{h}}^{n+1}_{i,j})^{\top}\boldsymbol{\Phi}(\xi_m)\nonumber\\ =& (\overline{\boldsymbol{h}}^{n}_{i,j})^{\top}\boldsymbol{\Phi}(\xi_m)-
        \frac{\Delta t^n}{\Delta x} \left[\mathcal{F}^{\hh}_{i+\frac{1}{2},j}(t_n)-\mathcal{F}^{\hh}_{i-\frac{1}{2},j}(t_n)\right]^{\top}\boldsymbol{\Phi}(\xi_m)\\
        &-\frac{\Delta t^n}{\Delta y}\left[\mathcal{G}^{\hh}_{i,j+\frac{1}{2}}(t_n)-\mathcal{G}^{\hh}_{i,j-\frac{1}{2}}(t_n)\right]^{\top}\boldsymbol{\Phi}(\xi_m).
      \end{align}    
      }
      Direct calculation and minimization over all the cells as well as all the quadrature points yields \eqref{eq:CFL}. \qed
\end{pf}
Combining the hyperbolicity-preserving condition \eqref{eq:CFL} with the wave-speed-based CFL stability conditions, we choose
\begin{align}\label{eq:tstep}
    \Delta t^n = 0.9\min\left\{\Delta t^n_{h}, \frac{\min\{\Delta x, \Delta y\}}{2a}\right\}.
\end{align}
Here, $a = \max_{i,j}\{a_{i+\frac{1}{2},j},b_{i,j+\frac{1}{2}}\}$ is the maximum possible magnitude of the propagation speeds at time level $t^n$, where
\begin{align*}
    &a_{i+\frac{1}{2},j} = \max\{a^{+}_{i+\frac{1}{2},j},-a^{-}_{i+\frac{1}{2},j}\},&&
    b_{i,j+\frac{1}{2}} = \max\{b^{+}_{i,j+\frac{1}{2}},-b^{-}_{i,j+\frac{1}{2}}\}.
\end{align*}
For now, we assume the local propagation speeds $a^{\pm}_{i+\frac{1}{2},j}, b^{\pm}_{i+\frac{1}{2},j}$ are real, i.e., the SG system is hyperbolic at the pointwise values of the reconstructions. The timestep condition analysis for \eqref{eq:tstep} is limited to the cell averages. For second-order or higher-order central upwind schemes, we need an additional treatment for the pointwise values $\boldsymbol{U}^{W,E,S,N}_{i,j}$ of the reconstructions. Further discussion on the procedure to preserve the hyperbolicity at $\boldsymbol{U}^{W,E,S,N}_{i,j}$ will be provided in \Cref{subsubsect:filtering}.

Our numerical scheme uses high-order strong stability-preserving Runge-Kutta schemes \cite{gottlieb2001strong} for time discretization, which corresponds to a convex combination of several forward Euler steps and thus the analysis above for deriving a hyperbolicity-preserving time step restriction still holds. An adaptive time control is introduced for the time discretization \cite{chertock2015well,kurganov2018finite}. Note that our analysis above can be applied to any other finite volume scheme.
\begin{rmk}
     If the signs of the fluxes at quadrature points are taken into account, the condition \eqref{eq:CFL} can be relaxed and the scheme will be more efficient.
\end{rmk}

\subsubsection{Hyperbolicity-Preserving Filtering for the Pointwise Values of the Reconstructions}\label{subsubsect:filtering}
Assuming $(\overline{\boldsymbol{h}}_{i,j}^n)^\top\boldsymbol{\Phi}(\xi_m) > 0$, we are able to guarantee that $(\overline{\boldsymbol{h}}_{i,j}^{n+1})^\top\boldsymbol{\Phi}(\xi_m) > 0$ for $m\in [M]$ by the condition \eqref{eq:tstep}, provided the local propagation speeds are real. However, the local propagation speeds are evaluated by the eigenvalues of the matrices $\frac{\partial \widehat{F}}{\partial \hU}$ and $\frac{\partial \widehat{G}}{\partial \hU}$ (see \eqref{eq:x-jacobian}-\eqref{eq:y-jacobian}) at $\boldsymbol{U}^{W,E,S,N}_{i,j}$ and, unfortunately, the hyperbolicity at the pointwise values of the reconstructions is not guaranteed by positivity of $(\overline{\boldsymbol{h}}_{i,j}^n)^\top\boldsymbol{\Phi}(\xi_m)$ due to the nonlinearity of the generalized minmod limiter. The problem can be resolved by filtering the higher-order stochastic moments \cite{schlachter2018hyperbolicity,doi:10.1137/20M1360736}.

Let $p_{\widehat{z}}(\xi) = \sum_{k=1}^{K}\widehat{z}_k\phi_k(\xi)$ be a polynomial in $P_\Lambda$ with positive moment $\widehat{z}_1$. Our goal is to determine the smallest possible weight $\mu'$ such that the weighted average of the polynomial $p_{\widehat{z}}(\xi)$ and the moment $\widehat{z}_1$ are nonnegative at given quadrature points $\{\xi_m\}_{m \in [M]}$, i.e.,
\begin{equation}\label{eq:filtering}
    \mu' \widehat{z}_1+(1-\mu')p_{\widehat{z}}(\xi_m)\ge 0\Leftrightarrow \widehat{z}_1+\sum_{k=2}^{K}(1-\mu')\widehat{z}_k\phi_k(\xi_m)\ge 0, m \in [M],
\end{equation}
and we filter the coefficients in the following way,
\begin{align}\label{eq:filterPCE}
    &\widehat{\mathsf{z}}_1 = \widehat{z}_1, &\widehat{\mathsf{z}}_k = (1-\mu)\widehat{z}_k, \;\;k = 2,\ldots, K,
\end{align}
where $\mu = \min\{\mu'+\delta, 1\}$, and we select $\delta = 10^{-10}$ in our scheme. The filtered polynomial $p_{\widehat{\mathsf{z}}}(\xi) = \sum_{k=1}^{K}\widehat{\mathsf{z}}_k\phi(\xi)$ is positive at quadrature points $\{\xi_m\}_{m\in[M]}$.

The hyperbolicity at the pointwise values of the reconstructions requires the positivity of the corresponding reconstructed water height to satisfy $(\boldsymbol{h}_{i,j}^{W,E,N,S})^\top\boldsymbol{\Phi}(\xi_m)>0$ for all $m\in [M]$. We will use the filter defined in \eqref{eq:filtering}-\eqref{eq:filterPCE} to filter the moments of these polynomials. Our strategy is to take the smallest filtering parameter $\mu^n_{i,j}$, among the filtering parameters of all the pointwise values $\boldsymbol{h}_{i,j}^{W,E,N,S}$ of the reconstructions in the same cell $\mathcal{C}_{i,j}$, and use it as the common filtering parameter for $\boldsymbol{h}_{i,j}^{W,E,N,S}$,
\begin{equation}\label{eq:filterh}
\begin{aligned}
      &\left(\mathsf{h}^{W,E,N,S}_{i,j}\right)_1 = \left(\boldsymbol{h}^{W,E,N,S}_{i,j}\right)_1,\\
  &\left(\mathsf{h}^{W,E,N,S}_{i,j}\right)_k = (1-\mu^n_{i,j})\left(\boldsymbol{h}^{W,E,N,S}_{i,j}\right)_k,k=2,\ldots, K.
\end{aligned}    
\end{equation}
The corresponding cell averages are recomputed accordingly,
\begin{equation}\label{eq:filterhavg}
    \overline{\mathsf{h}}^n_{i,j} = \frac{1}{4}\left(\mathsf{h}^{W}_{i,j}+\mathsf{h}^{E}_{i,j}+\mathsf{h}^{N}_{i,j}+\mathsf{h}^{S}_{i,j}\right).
\end{equation}
With \eqref{eq:filterhavg}, we are also able to preserve the first moment of the cell averages before filtering. There are other choices of filters to further reduce the oscillations (e.g., see \cite{kusch2020filtered,schlachter2020weighted,kusch2020oscillation}), but we adopt this filter because of its simplicity. Note that, in \eqref{eq:filterh}-\eqref{eq:filterhavg}, we only filter the polynomial moments for $\boldsymbol{h}^{W,E,N,S}$ and do not filter the discharges. 
\begin{rmk}\label{rmk:filter}
  We use a filtering strategy because we cannot theoretically guarantee that our second-order reconstructions $\bh^{W,E,N,S}_{i,j}$ are positive, even if the corresponding cell averages $\bh_{i,j}$ and positive and the stochastic pointwise assumption \eqref{eq:h-positivity} from \Cref{thm:h-positivity} is satisfied.
     If a first-order numerical scheme is used, i.e., the approximated solutions are piecewise constant functions, the filtering procedure is not necessary as long as $\bh^n_{i,j}$ satisfies \Cref{thm:h-positivity} for all $i,j$. A positivity-preserving condition for such an approximation can be derived in a similar manner to \Cref{lemma:CFL}.
\end{rmk}

The above strategy requires that the first moments of $\boldsymbol{h}_{i,j}^{W,E,N,S}$ are positive during computation, but this is not necessarily the case. Nonpositive first moments may lead to the failure of the filtering procedure \eqref{eq:filterh}. This is an issue that happens even for the deterministic system of 2D SWE \cite{kurganov2007second}. Therefore, we propose the following correction for the PCE vector of the reconstructed water height with nonpositive first moment,
\begin{align*}
  \textrm{if } \left(\boldsymbol{h}^{W}_{i,j}\right)_1\le0 \hskip 5pt 
  \textrm{then take } \boldsymbol{h}^{W}_{i,j}=\boldsymbol{0},\;\;
                      \boldsymbol{h}^{E}_{i,j}=2\overline{\boldsymbol{h}}^n_{i,j},\\
  \textrm{if } \left(\boldsymbol{h}^{E}_{i,j}\right)_1\le0 \hskip 5pt 
  \textrm{then take } \boldsymbol{h}^{E}_{i,j}=\boldsymbol{0},\;\;
                      \boldsymbol{h}^{W}_{i,j}=2\overline{\boldsymbol{h}}^n_{i,j},\\
  \textrm{if } \left(\boldsymbol{h}^{N}_{i,j}\right)_1\le0 \hskip 5pt 
  \textrm{then take } \boldsymbol{h}^{N}_{i,j}=\boldsymbol{0},\;\;
                      \boldsymbol{h}^{S}_{i,j}=2\overline{\boldsymbol{h}}^n_{i,j},\\
  \textrm{if } \left(\boldsymbol{h}^{S}_{i,j}\right)_1\le0 \hskip 5pt 
  \textrm{then take } \boldsymbol{h}^{N}_{i,j}=\boldsymbol{0},\;\;
                      \boldsymbol{h}^{N}_{i,j}=2\overline{\boldsymbol{h}}^n_{i,j}.
\end{align*}
In this way, the cell averages remain consistent, i.e., the cell average is equal to the average of the pointwise values for the corresponding cell. Note that, in the case when the first moment is corrected, the filter will not be applied. This correction reduces to a similar procedure in the central-upwind scheme for the deterministic shallow water equations \cite{kurganov2007second}.
\subsubsection{Velocity Desingularization}
In computation, the matrix $\mathcal{P}(\boldsymbol{h}^{W,E,N,S}_{i,j})$ may be near singular at the pointwise values of the reconstructions and may result in large round-off errors when determining the PCE vector for the velocity \eqref{eq:uPCE}. An approach for stochastic velocity desingularization is proposed in \cite{doi:10.1137/20M1360736}. To recall the idea, we consider $\mathcal{P}(\boldsymbol{h}^{W}_{i,j})$ as an example. But our process can be extended naturally to the other point values of the reconstructions. Let 
\begin{align}
    \mathcal{P}(\boldsymbol{h}^{W}_{i,j}) = Q^{\top}\Sigma Q,
\end{align}
be the eigenvalues decomposition of $\mathcal{P}(\boldsymbol{h}^{W}_{i,j})$, where $\Sigma = \text{diag}(\sigma_1,\ldots,\sigma_K)$. For $k\in [K]$ and a prescribed $\epsilon$, defined the corrected inverse of $\Sigma$ to be, 
\begin{align}\label{eq:invcorrected}
  \Sigma^{\co} &= \text{diag}(\sigma^{\co}_1,\ldots,\sigma^{\co}_K),&
    \sigma^{\co}_k = \frac{\sqrt{2}\sigma_{k}}{\sqrt{\sigma^4_{k}+\max\{\sigma^4_{k},\epsilon^4\}}}.%k = 1,\ldots,K.
\end{align}
In our simulation, $\epsilon=\min\{\Delta x, \Delta y\}$. The PCE coefficients for the velocity vector are then computed by,
\begin{align}\label{eq:udesin}
    &\boldsymbol{u}^{W}_{i,j} = Q^{\top}\Sigma^{\co} Q(\boldsymbol{q}^x)^{W}_{i,j},&&\boldsymbol{v}^{W}_{i,j} = Q^{\top}\Sigma^{\co} Q(\boldsymbol{q}^y)^{W}_{i,j}.
\end{align}
For well-conditioned $\mathcal{P}(\boldsymbol{h}^{W}_{i,j})$, the correction \eqref{eq:udesin} reduces to \eqref{eq:uPCE} and no further reconstruction is needed. In the case when $\mathcal{P}(\boldsymbol{h}^{W}_{i,j})$ is ill-conditioned, to remain the consistency of the PCE vector of the discharges, we recompute the discharges using the corrected velocity
\begin{align}
    &(\boldsymbol{q}^x)^{W}_{i,j} = \mathcal{P}(\boldsymbol{h}^{W}_{i,j})\boldsymbol{u}^{W}_{i,j},
    &&(\boldsymbol{q}^y)^{W}_{i,j} = \mathcal{P}(\boldsymbol{h}^{W}_{i,j})\boldsymbol{v}^{W}_{i,j}.
\end{align}
If there is no uncertainty in the system (all moments beyond the first vanish), the desingularization procedure reduces to the deterministic velocity desingularization procedure proposed in \cite{kurganov2007second,kurganov2018finite}. Readers can also find other desingularization options there.
\section{Numerical Results}\label{sec:results}
%
%\begin{itemize}
%    \item Outline of numerical example
%    \begin{itemize}
%        \item example 1: order of convergence,
%        \item well-balanced example 1
%        \begin{itemize}
%            \item collocation comparison I (mean and std contours of the water)
%        \end{itemize}
%        \item well-balanced example 1 (tensor version 2): (disk glyphs, different setups of the uncertainties)        
%        \item well-balanced example 2: discontinuous bottom
%        \item example 1 (tensor version)
%    \end{itemize}
%\end{itemize}
In this section, we will demonstrate the capability of our scheme to capture the solutions to small perturbations of the stochastic ``lake-at-rest'' solution for different choices of the random variables $\xi$. The distributions we consider are from the families of distributions with Beta density. For a single random variable $\xi^{(1)}$, its density is given by,
{\small
\begin{align}\label{eq:betafamily}
  \rho(\xi^{(1)}) \coloneqq \rho^{(\alpha,\beta)}(\xi^{(1)}) &= C(\alpha,\beta) (1-\xi^{(1)})^\alpha (1+\xi^{(1)})^\beta, & 
  C(\alpha,\beta)^{-1} &= 2^{\alpha+\beta+1} \mathcal{B}(\beta+1, \alpha+1),
\end{align}
}
where $\mathcal{B}(\cdot,\cdot)$ is the Beta function, and the parameters $\alpha, \beta > -1$ can be chosen freely and control how mass concentrates at $\xi^{(1)} = 1$ and $\xi^{(1)} = - 1$, respectively. The associated orthonormal polynomial basis is the Jacobi polynomials with parameter $\alpha$ and $\beta$. In particular, $\alpha = \beta = 0$ corresponds to the uniform distribution $\mathcal{U}(-1,1)$ and the Legendre polynomial basis. 

In our numerical tests, we consider the cases with both the one-dimensional random variable $\xi = \xi^{(1)}$ and the two-dimensional random variable $\xi = (\xi^{(1)},\xi^{(2)})^\top$, where $\xi^{(1)}$ and $\xi^{(2)}$ are independent one-dimensional random variables with different distributions, and thus the density function for $\xi$ is defined by $\rho(\xi)\coloneqq\rho(\xi^{(1)})\rho(\xi^{(2)})$. The multi-index set $\Lambda$ associated with the polynomial space $P_\Lambda$ is chosen to be the tensor-product set,
\begin{align}\label{eq:index1}
    \Lambda = \{(\nu^{(1)},\nu^{(2)})\in \mathbbm{N}^2\;\;\vert\;\;0\le\nu^{(1)},\nu^{(2)}\le 3\},
\end{align}
for all the two-dimensional random field numerical examples.

Throughout the numerical section, the gravitational constant is $g = 1$, the parameter for the generalized minmod limiter is set to be $\theta = 1.3$ (except for \Cref{sect:sub-island}, where $\theta = 1.0$). The filtering parameter is set to be $\delta = 10^{-10}$. 
We further observe that, in all of our numerical examples, only \Cref{sect:filter} requires filters.

\subsection{Example 1: Filtering}\label{sect:filter}
As we mentioned in Remark \ref{rmk:filter}, filtering is employed to guarantee that the pointwise reconstructions $\boldsymbol{h}_{i,j}^{{W,E,N,S}}$ satisfy the sufficient condition \Cref{thm:h-positivity}. This section demonstrates the necessity of filtering by providing a comparison of the filtered and nonfiltered numerical results. 

\begin{figure}[htbp]
        \centering
        \subfloat{\includegraphics[width=0.45\textwidth]{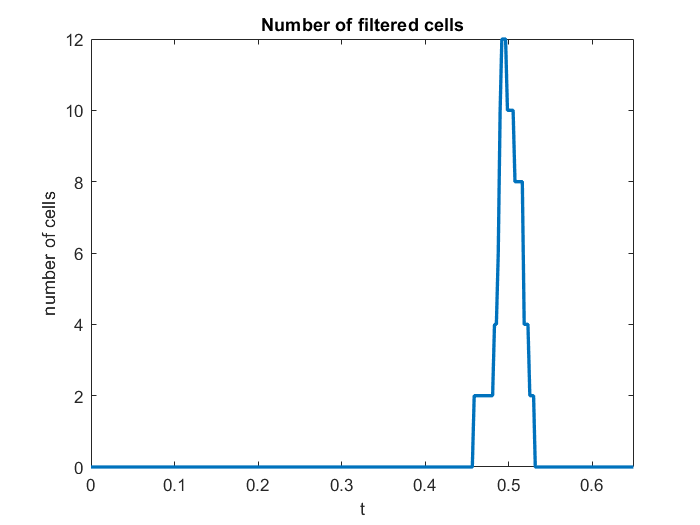}}        
        \subfloat{\includegraphics[width=0.45\textwidth]{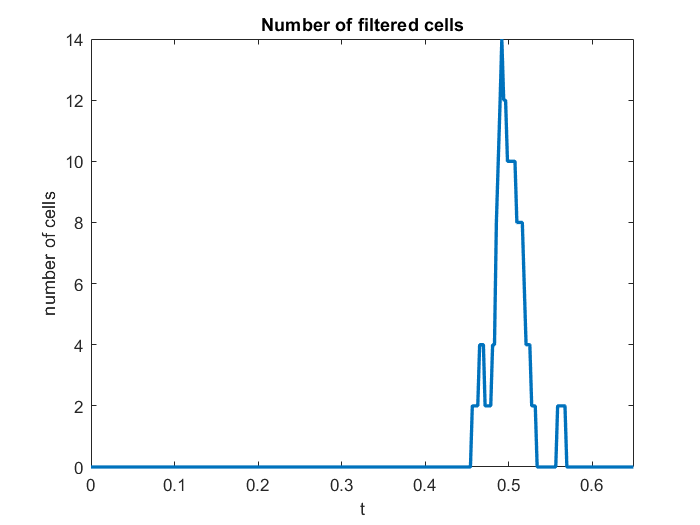}}\\        
        \subfloat{\includegraphics[width=0.45\textwidth]{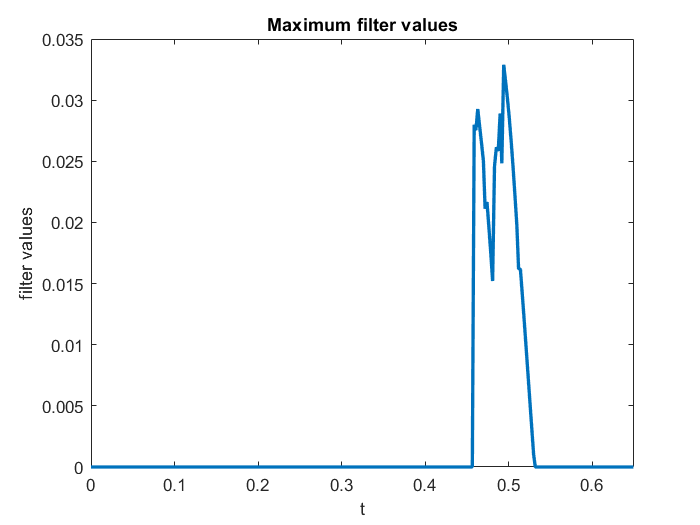}}
       \subfloat{\includegraphics[width=0.45\textwidth]{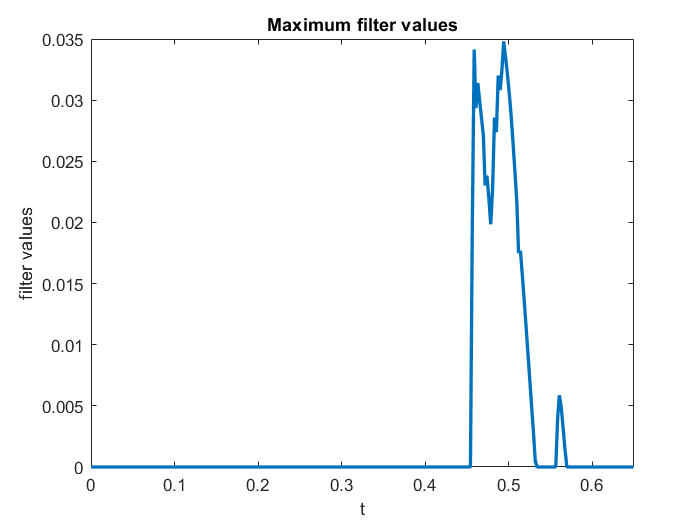}}
     \caption{Study of the filtering procedure for \eqref{eq:initwb3-1-uq}-\eqref{eq:bottomwb3-1-uq}, grid size $200\times 200$. Top: number of filtered cells; Bottom: max magnitude of the filter $\mu$. Left: $K=4$; Right: $K=8$.}
        \label{fig:filter-compare}   
    \end{figure}

We consider a deterministic water surface with  perturbation
\begin{equation}\label{eq:initwb3-1-uq}
   \eta(x,y,0) = \left\{\begin{aligned}
    &1+0.0001(\xi+1),&&\text{if } -0.4<x<-0.3,\\
    &1,&&\text{otherwise},
    \end{aligned}\right.,\quad u(x,y,0) = v(x,y,0) = 0,
\end{equation}
and with  bottom topography,
\begin{equation}\label{eq:bottomwb3-1-uq}
    B(x,y,\xi) = \left\{\begin{aligned}
    &0.9998,&& r\le 0.1,\\
    &9.998(0.2-r),&& 0.1<r\le0.2,\\
    &0.0001,&&\text{otherwise},
    \end{aligned}\right. 
\end{equation}
The one-dimensional random variable is chosen to be $\xi\sim \mathcal{U}(-1,1)$. We run the numerical test for the following $4$ sets of parameters up to time $t = 0.65$ for both $K=4$ and $K=8$,
\begin{itemize}
    \item $(nx,ny) = (100,100)$, filtered = True,
    \item $(nx,ny) = (200,200)$, filtered = True,
    \item $(nx,ny) = (100,100)$, filtered = False,
    \item $(nx,ny) = (200,200)$, filtered = False.
\end{itemize}
We first present a table of comparison of running time (in seconds) \Cref{tab:timecomp}. The text ``Failed'' indicates that \Cref{eq:swesg4} loses hyperbolicity at some pointwise reconstructions before $t = 0.65$. We find that, for a $100\times 100$ grid, the filtering procedure was never activated ($\mu\equiv 0$ in \Cref{eq:filterPCE}). The slight difference in time for the $100\times 100$ grid is from determining whether a reconstruction satisfies the sufficient condition \Cref{thm:h-positivity}. The filtering scheme is more robust as it consistently avoids possible loss of hyperbolicity at pointwise reconstructed values. We further keep track of the number of filtered cells and the maximum magnitude of the filter values $\mu$ at each time step (\Cref{fig:filter-compare}). The number of cells that need filtering is relatively small ($\sim 0.03\%$ of the total cells) and the filter values $\mu$ are relative small. 
%Qualitatively, the numerical results for filtered and nonfiltered schemes are similar \Cref{fig:filter-compare}.

\begin{table}[htbp]
\centering
  \begin{tabular}{ccccc}
    \toprule
    \multirow{2}{*}{Grid size} &
      \multicolumn{2}{c}{$K=4$} &
      \multicolumn{2}{c}{$K=8$} \\
      & {Filtered} & {Nonfiltered} & {Filtered} & {Nonfiltered} \\
      \midrule
    $100\times 100$ & 450.74s
     & 441.27s & 1597.25s & 1556.59s \\
    $200\times 200$ & 3868.11s & Failed &
    13459.43s & Failed \\
    \bottomrule
  \end{tabular}
\caption{Comparison of time for \eqref{eq:initwb3-1-uq}-\eqref{eq:bottomwb3-1-uq} with different choices of parameters.}
    \label{tab:timecomp}
\end{table}

%    \begin{figure}[htbp]
%        \centering
%        \subfloat{\includegraphics[width=0.45\textwidth]{wb3-1-dglyph-mean-std-N4-nx100ny100T65.png}}        
%        \subfloat{\includegraphics[width=0.45\textwidth]{wb3-1-nf-dglyph-mean-std-N4-nx100ny100T65.png}}\\        
%        \subfloat{\includegraphics[width=0.45\textwidth]{wb3-1-dglyph-mean-std-N8-nx100ny100T65.png}}
%        \subfloat{\includegraphics[width=0.45\textwidth]{wb3-1-nf-dglyph-mean-std-N8-nx100ny100T65.png}}
%        \caption{Numerical solutions to \eqref{eq:initwb3-1-uq}-\eqref{eq:bottomwb3-1-uq}, water surface at $t = 0.65$, disk-glyph over mean contours, where the radii of the disks indicate the magnitude of the standard deviation, grid size $100\times 100$. Top: $K=4$; Bottom: $K=8$. Left: filtered; Right: nonfiltered. The largest disks are corresponding to the standard deviation values 2.25e-5, 2.25e-5, 2.25e-5, and 2.25e-5, respectively.}
%        \label{fig:filter-compare}    
%    \end{figure}    

\subsection{Example 2: Accuracy Test}\label{sect:accuracy}
\begin{figure}[htbp]
    \centering
    \subfloat{\includegraphics[width=0.45\textwidth]{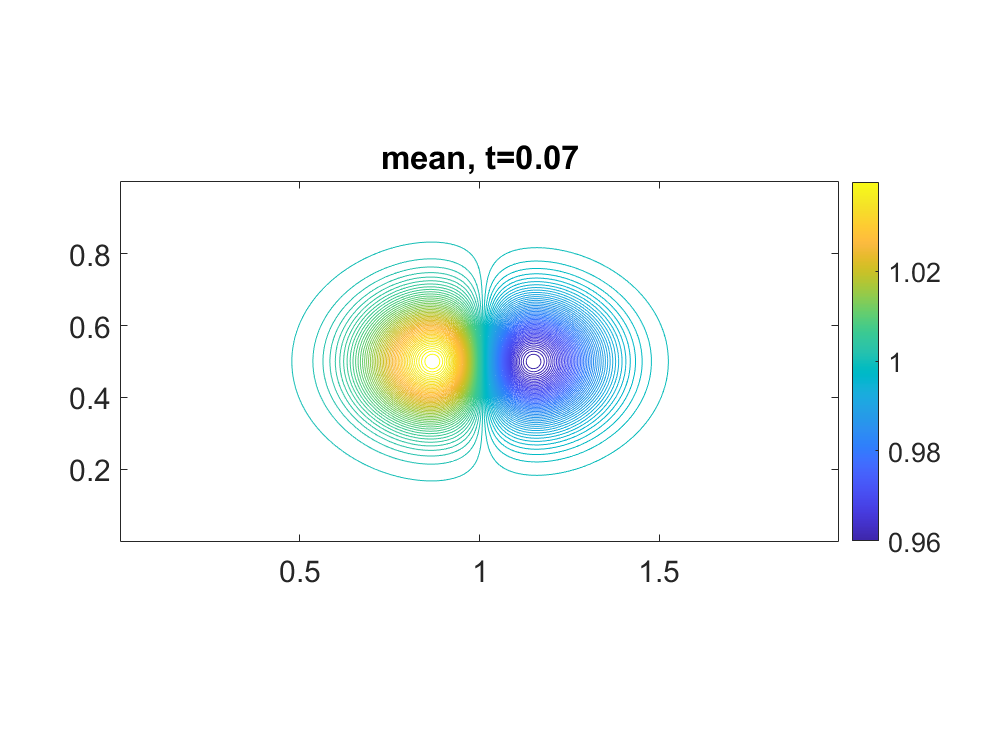}}
    \subfloat{\includegraphics[width=0.45\textwidth]{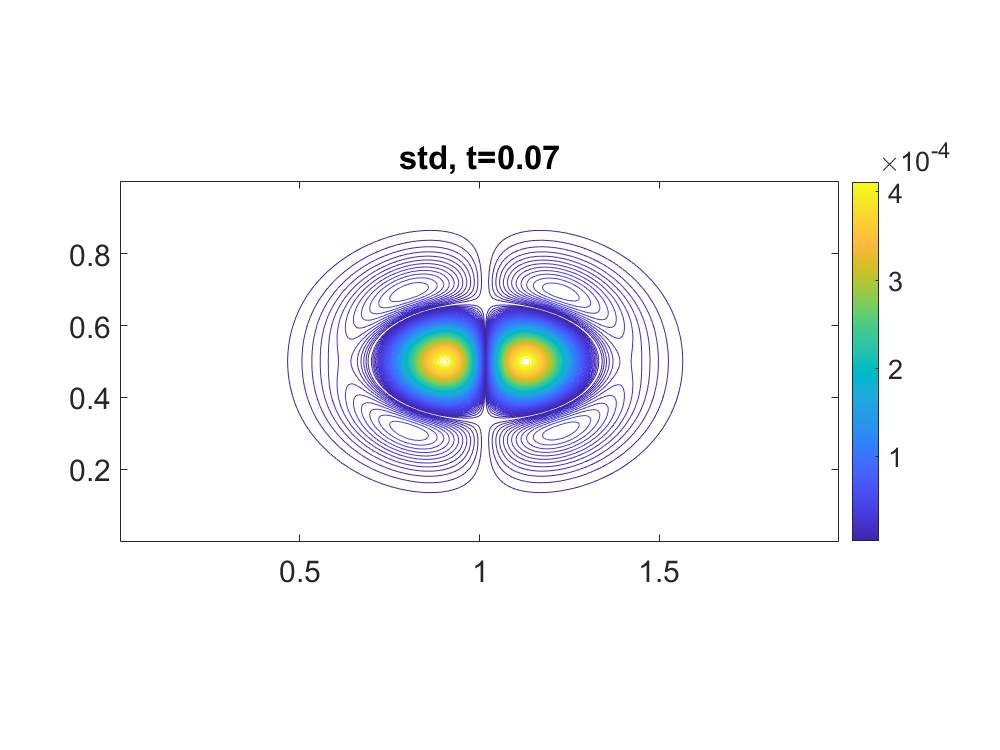}} 
    \caption{Contours for the water surface of the reference solution. Left: mean. Right: standard deviation at $t = 0.07$}
    \label{fig:ex1-1-uq-time}    
\end{figure}

This subsection is devoted to checking the order of accuracy of the proposed scheme for our SG system of 2D SWE. We consider a stochastic version of the accuracy test in \cite{M2AN_2011__45_3_423_0,SHIRKHANI201625,ghazizadeh2020adaptive} with deterministic initial water surface,
\begin{equation}\label{eq:initex1-1uq-conv}
    \eta(x,y,0,\xi) = 1,\quad u(x,y,0,\xi) = 0.3,\quad v(x,y,0,\xi) = 0,
\end{equation}
and stochastic elliptic-shaped hump bottom,
\begin{equation}\label{eq:bottomex1-1uq-conv}
    B(x,y) = 0.5e^{-25(x-1)^2-50(y-0.5)^2} + 0.1(\xi+1).
\end{equation}
The one-dimensional random variable is chosen to be $\xi\sim \mathcal{U}(-1,1)$. The computational domain is $[0,2]\times [0,1]$. All the boundary conditions are obtained by zero-order extrapolation. The reference solution is computed on a $800\times 800$ uniform rectangular grid. 
%By $t = 0.07$, the solution converges to the steady state, which, for this example, is non-constant but smooth. 
The contour plots of the reference solution for the mean water surface and the standard deviation of the water surface at $t = 0.07$ are shown in \Cref{fig:ex1-1-uq-time}. 

We computed the error between the reference solution and the numerical solution using,
\begin{align}\label{eq:err}
    \text{err}(U_d) = \Vert U_d(x,y,\xi) - U_{\text{ref}}(x,y,\xi)\Vert_{L^1(\Omega_{xy};L^2_{\rho}(\Omega_\xi;\mathbbm{R}^3))},
\end{align}
where $\Omega_{xy}$ is the physical domain, $\Omega_\xi$ is the stochastic domain, $U_d$ is the numerical $K$-term PCE solution to the $\xi$-parameterized equation \eqref{eq:swesg1}, and $U_{\text{ref}}$ is the $K$-term PCE reference solution to the $\xi$-parameterized equation \eqref{eq:swesg1}. Here, for the stochastic state variable in $K$-term PCE, $U = (h_\Lambda,q^x_\Lambda, q^y_\Lambda)^\top$, the norm $\Vert\cdot\Vert_{L^1(\Omega_{xy};L^2_{\rho}(\Omega_\xi;\mathbbm{R}^3))}$ is defined to be, 
\begin{align*}
    &\Vert U(x,y,t,\xi)\Vert_{L^1(\Omega_{xy};L^2_{\rho}(\Omega_\xi;\mathbbm{R}^3))}\coloneqq\\ &\int_{\Omega_{xy}}\left\vert \Vert h_\Lambda(x,y,t,\xi)\Vert_{L^2_\rho}+\Vert q^x_\Lambda(x,y,t,\xi)\Vert_{L^2_\rho}+\Vert q^y_\Lambda(x,y,t,\xi)\Vert_{L^2_\rho} \right\vert dx\;dy.
\end{align*}
From the orthogonality of the basis functions in the PCE, we have, 
$$
\Vert h_\Lambda(x,y;\xi)\Vert_{L^2_\rho} = \Vert\hh\Vert_{2},
$$ 
where $\hh$ is the vector of the PCE coefficients of $h_\Lambda$, and $\Vert\cdot\Vert_2$ is the Euclidean norm. The computations for the norm of $q^x_\Lambda$ and $q^y_\Lambda$ are similar. \rone{We compute orders of convergence for the solution with $K=4$, $K=8$ at various times $t$ in  \Cref{tab:conv}. The rate of convergence is comparable to the rates of the deterministic benchmark reported in \cite{M2AN_2011__45_3_423_0,SHIRKHANI201625,ghazizadeh2020adaptive} and we expect that using a finer grid for the reference solution will yield second-order convergence. We also note that, for this example, increasing the polynomial order will not have a significant effect on the order of the convergence, especially when $K$ is relatively large ($K>8$).}

\begin{table}
\centering
\rone{
  \begin{tabular}{cccccc}
    \toprule
    \multirow{2}{*}{Time}&
    \multirow{2}{*}{Grid size} &
    \multicolumn{2}{c}{$K=4$} &
    \multicolumn{2}{c}{$K=8$} \\
    && {Error} & {Order} & {Error} & {Order} \\
    \midrule
    \multirow{3}{*}{$t=0.067$}&$100\times 100$ & 1.475875e-05
     & $-$ & 1.456866e-05 & $-$ \\
    &$200\times 200$ & 4.287569e-06 & 1.764646 &
    4.287573e-06 & 1.764635 \\
    &$400\times 400$ & 1.279305e-06 & 1.744799 & 
    1.279306e-06 &  1.744799 \\
    \hline
    \multirow{3}{*}{$t=0.07$}&$100\times 100$ & 1.475875e-05
     & $-$ & 1.475865e-05 & $-$ \\
    &$200\times 200$ & 4.343711e-06 & 1.764571 &
    4.343714e-06 & 1.764559 \\
    &$400\times 400$ & 1.296122e-06 & 1.744727 & 
    1.296122e-06 &  1.744727 \\
    \hline    
    \multirow{3}{*}{$t=0.073$}&$100\times 100$ & 1.494536e-05
     & $-$ & 1.494524e-05 & $-$ \\
    &$200\times 200$ & 4.398829e-06 & 1.764506 &
    4.398832e-06 & 1.764494 \\
    &$400\times 400$ & 1.312592e-06 & 1.744700 & 
    1.312593e-06 &  1.744701 \\
    \hline
    \multirow{3}{*}{$t=0.14$}&$100\times 100$ & 1.808672e-05
     & $-$ & 1.808651e-05 & $-$ \\
    &$200\times 200$ & 5.329017e-06 & 1.762990 &
    5.329023e-06 & 1.762971 \\
    &$400\times 400$ & 1.590317e-06 & 1.744555 & 
    1.590320e-06 & 1.744554 \\    
    \bottomrule
  \end{tabular}}
    \caption{Order of convergence for \eqref{eq:initex1-1uq-conv}-\eqref{eq:bottomex1-1uq-conv} with different grid sizes and different $K$ values \rone{at different times.}}
    \label{tab:conv}
\end{table}
\subsection{Example 3: Comparison with Collocation Method}\label{sect:sc-comparison}
\begin{figure}[htbp]
        \centering
        \subfloat{\includegraphics[width=0.45\textwidth]{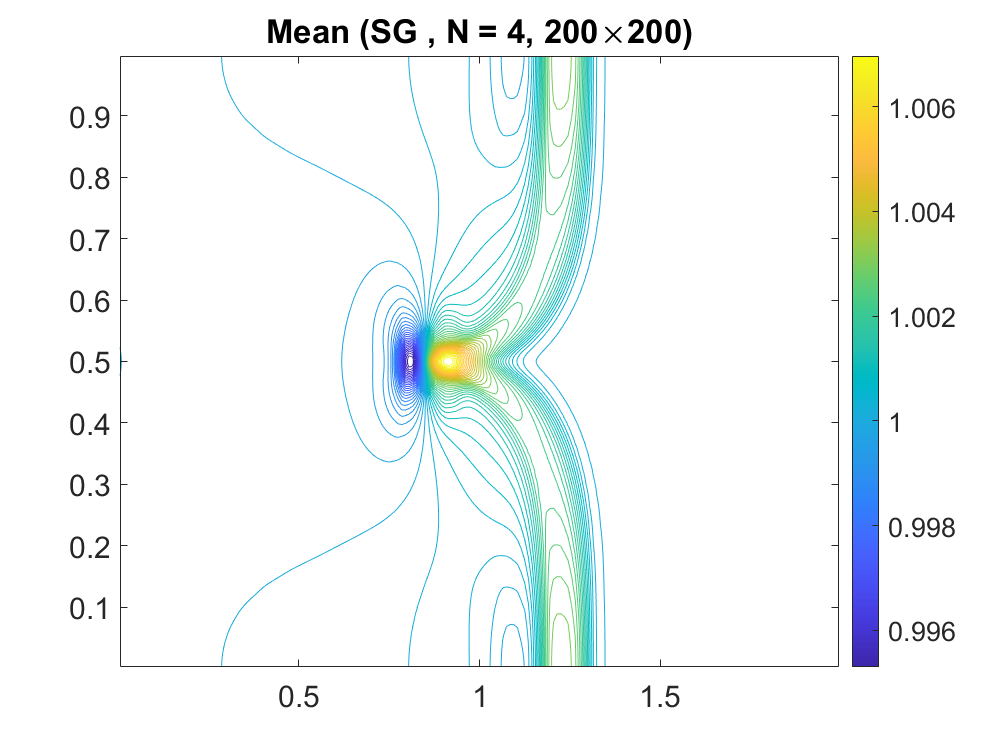}}        
        \subfloat{\includegraphics[width=0.45\textwidth]{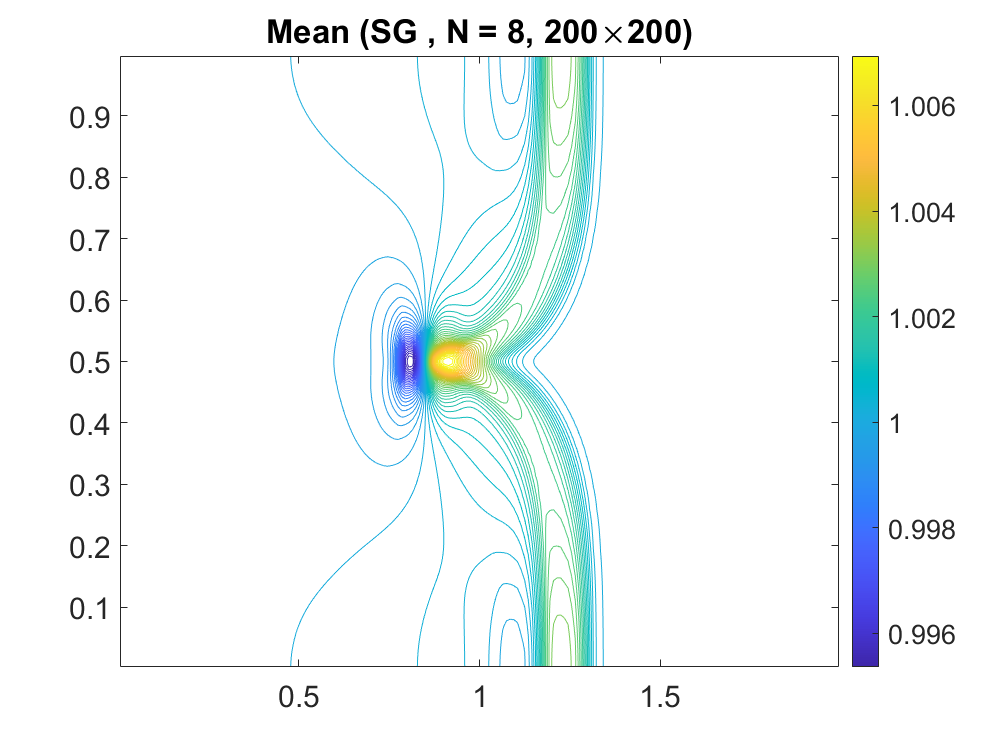}}\\        
        \subfloat{\includegraphics[width=0.45\textwidth]{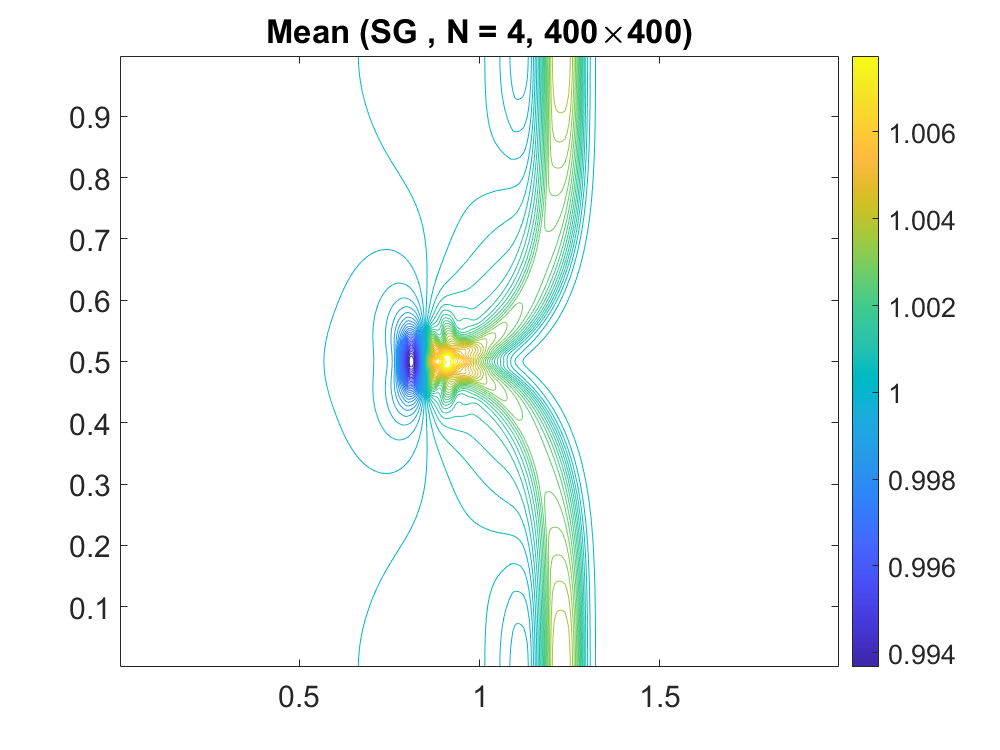}}
        \subfloat{\includegraphics[width=0.45\textwidth]{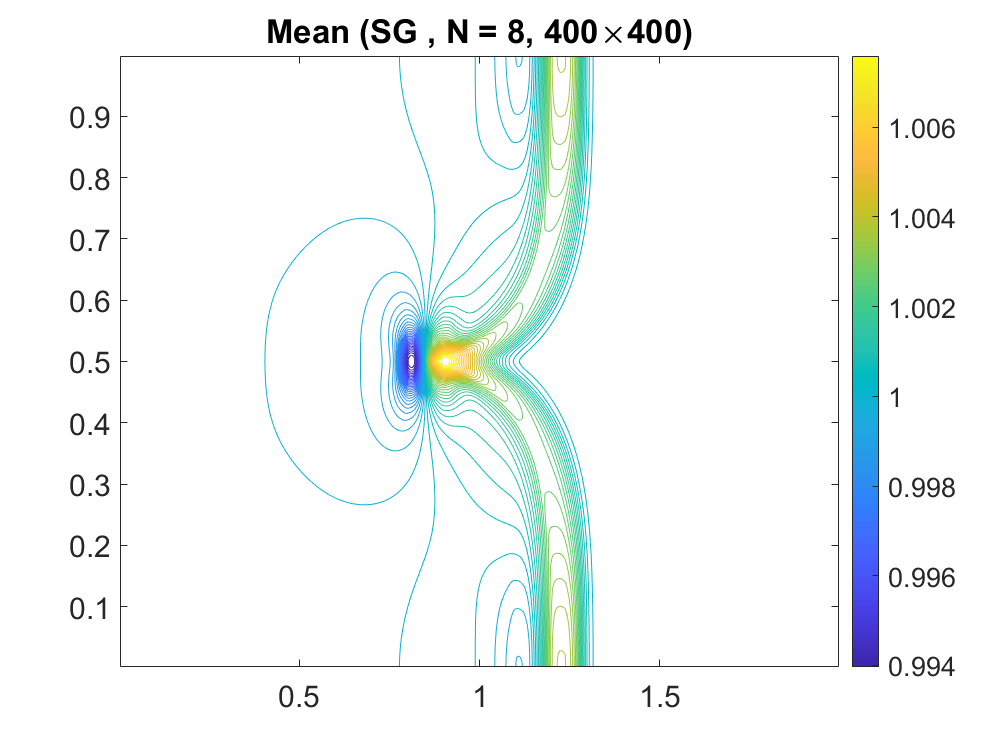}}\\
        \subfloat{\includegraphics[width=0.45\textwidth]{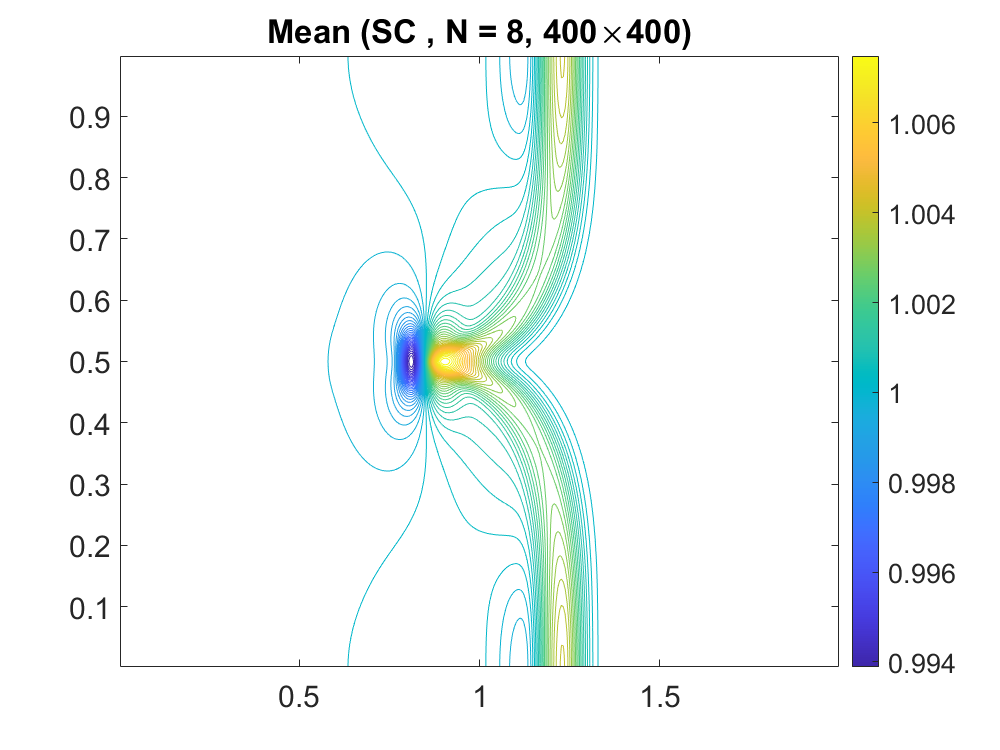}}
        \caption{Numerical solutions to \eqref{eq:initwb1-1uq-conv}-\eqref{eq:bottomwb1-1uq-conv}, mean water surface at $t = 1.2$, grid size $400\times 400$. Top left: SG solution, $K=4$, $200\times 200$ grid; Top right: SG solution, $K = 8$, $200\times 200$ grid. Middle left: SG solution, $K=4$, $400\times 400$ grid; Middle right: SG solution, $K = 8$, $400\times 400$ grid. Bottom: SC solution, $K = 8$, $400\times 400$ grid. The number of collocation points is $100$.}
        \label{fig:wb1-1-uq-mean-conv}    
    \end{figure}

\begin{figure}[htbp]
        \centering
        \subfloat{\includegraphics[width=0.45\textwidth]{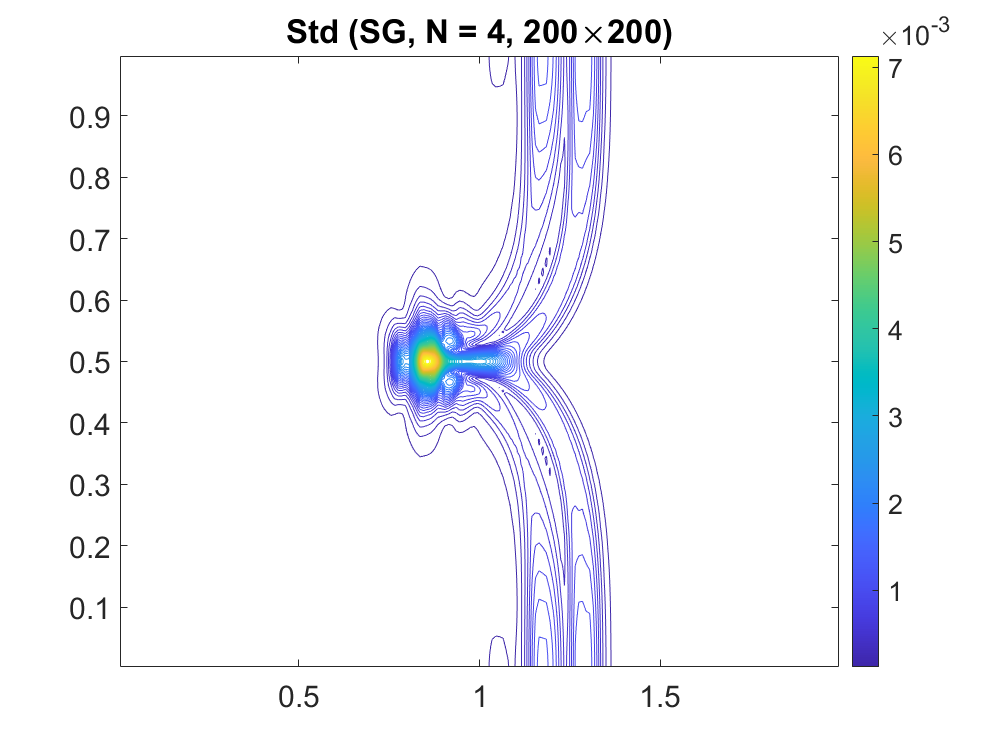}}
        \subfloat{\includegraphics[width=0.45\textwidth]{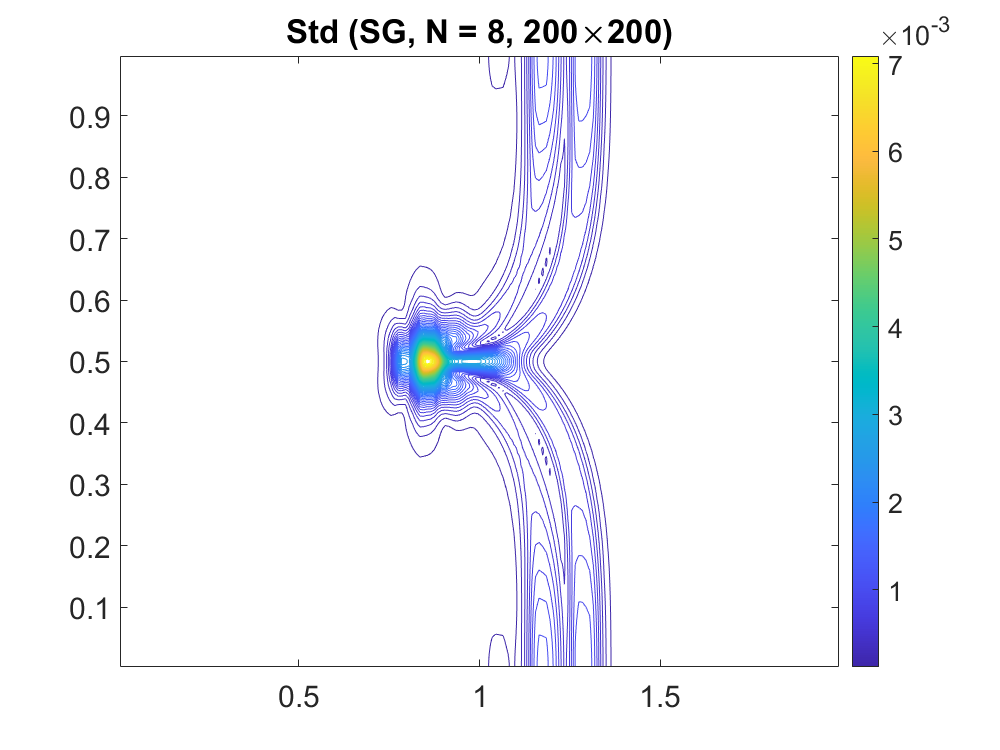}}\\        
        \subfloat{\includegraphics[width=0.45\textwidth]{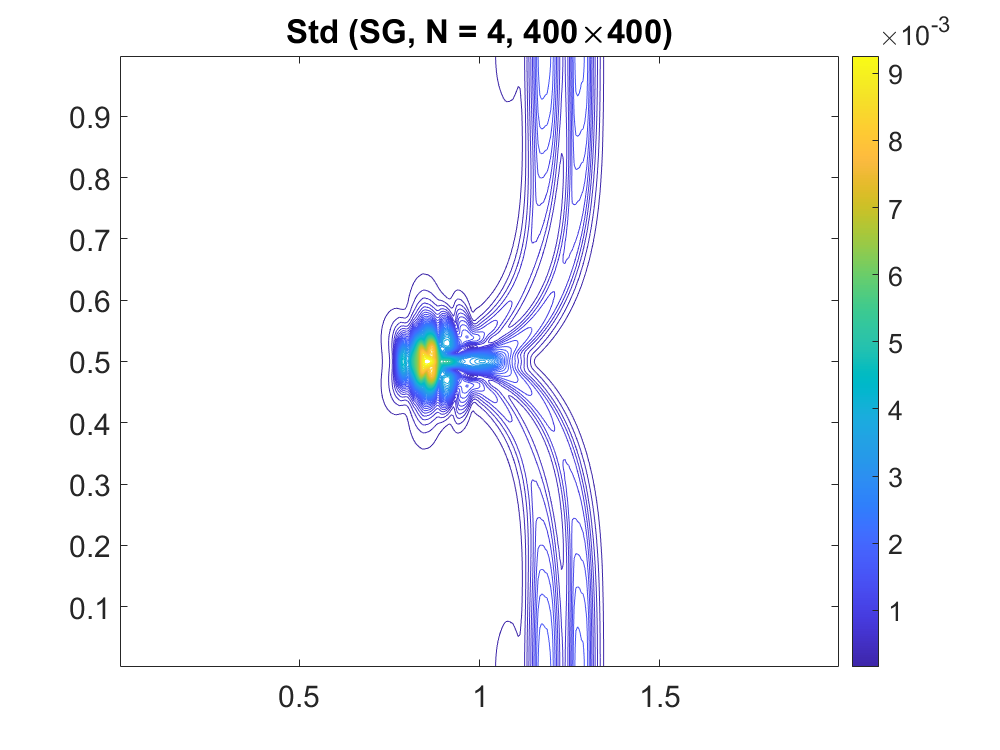}}
        \subfloat{\includegraphics[width=0.45\textwidth]{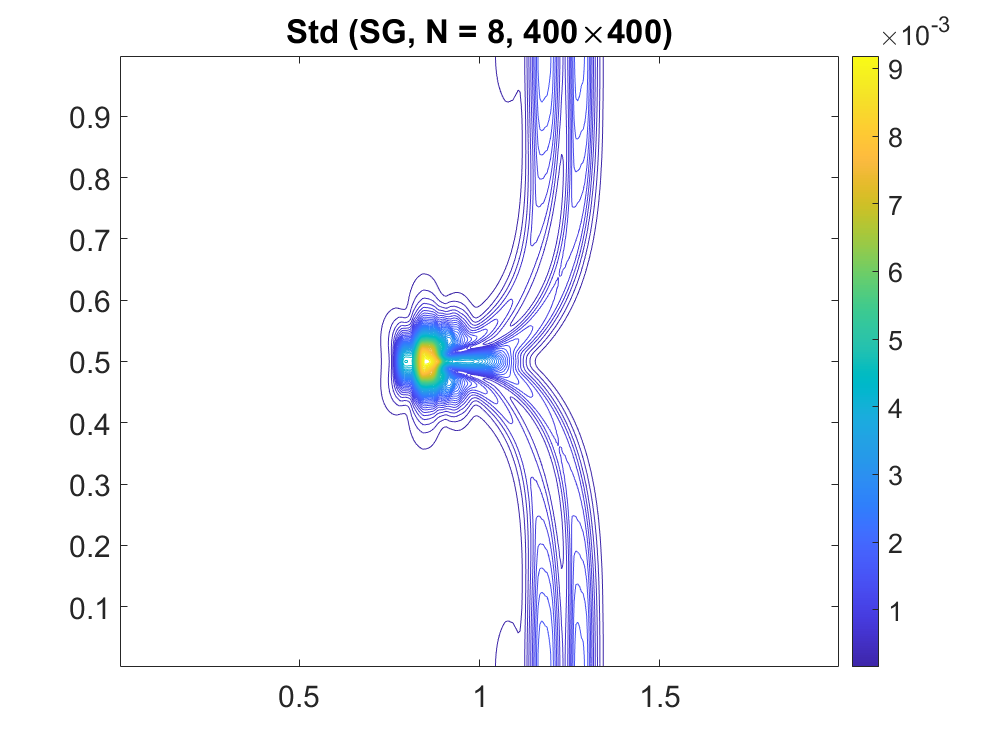}}\\
        \subfloat{\includegraphics[width=0.45\textwidth]{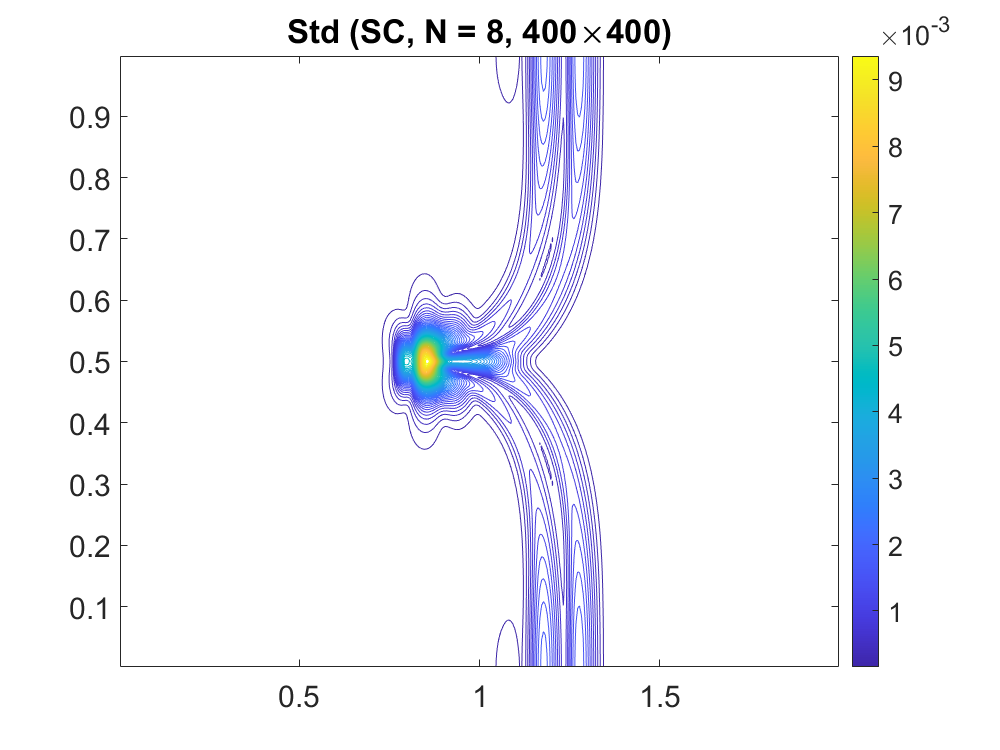}}
        \caption{Numerical solutions to \eqref{eq:initwb1-1uq-conv}-\eqref{eq:bottomwb1-1uq-conv}, standard deviation of the water surface at $t = 1.2$.Top left: SG solution, $K=4$, $200\times 200$ grid; Top right: SG solution, $K = 8$, $200\times 200$ grid. Middle left: SG solution, $K=4$, $400\times 400$ grid; Middle right: SG solution, $K = 8$, $400\times 400$ grid. Bottom: SC solution, $K = 8$, $400\times 400$ grid. The number of collocation points is $100$.}
        \label{fig:wb1-1-uq-std-conv}    
    \end{figure}

The second numerical example serves to compare our SG formulations \eqref{eq:swesg4} as well as our scheme with a collocation procedure. Consider a deterministic initial water surface
    \begin{equation}\label{eq:initwb1-1uq-conv}
        \eta(x,y,0,\xi) = \left\{\begin{aligned}
        &1.01,&&\text{if } 0.05<x<0.15,\\
        &1,&&\text{otherwise},
        \end{aligned}\right.\quad u(x,y,0,\xi) = v(x,y,0,\xi) = 0,
    \end{equation}
    with an uncertain bottom function
    \begin{equation}\label{eq:bottomwb1-1uq-conv}
        B(x,y,\xi) = 0.8e^{-5(x-0.9)^2-50(y-0.5)^2} + 0.1(\xi+1),
    \end{equation}    
where $\xi$ has Beta density with parameters $(\alpha,\beta) = (1,3)$. At $t=0$, there is a small but deterministic perturbation to the ``lake-at-rest'' state $\eta\equiv 1$ at $0.05<x<0.15$. The uncertainty of the bottom skews toward the positive direction, i.e., the probability for $\xi$ to be a positive number is higher. The computational domain is $[0,2]\times [0,1]$. The boundary conditions for the upper and the lower boundaries are periodic boundary conditions. The boundary conditions for the left and right boundaries are outflow boundary conditions. We compute the SG solution for $K=4, 8$ on $200\times 200$ and $400\times 400$ grid up to time $T = 1.2$. As a comparison, we compute the stochastic collocation (SC) solution $U_{SC}$ for $K=8$ with $M = 100$ collocation points on $400\times 400$ grid for the same terminal time, where the SC solution is obtained by 
{\small
\begin{align*}
    &U_{SC}(x,y,t,\zeta) = \sum_{k=1}^{K}\widehat{U}_{SC,k}(x,y,t)\phi_k(\zeta),&&\widehat{U}_{SC,k}(x,y,t) \coloneqq  \sum_{m=1}^{M}U_{SC}(x,y,t,\zeta_m)\phi_k(\zeta_m)z_m.
\end{align*}
}
Here, $\{\zeta_m, z_m\}_{m\in [M]}$ is an $M$-point $\rho$-Gaussian quadrature rule. The comparison of the contours of the statistics of the water surface is shown in \Cref{fig:wb1-1-uq-mean-conv} (mean) and \Cref{fig:wb1-1-uq-std-conv} (standard deviation). Since the SG solutions and the SC solutions are obtained from numerically solving different partial differential equations, we do not expect that one solution will converge to the other solution by, for example, increasing the number $M$ of quadrature points and/or refining the mesh. The contour lines look different near the intermediate region from peak to valley of the mean water surface, but overall they are similar, especially between the SG solutions on $400\times 400$ grid with $K=8$ and the SC solution.
\subsection{Example 4: A Two-Dimensional Random Variable Example}\label{sect:wb1-2d}
In this example, we will consider a more difficult variant of \eqref{eq:initwb1-1uq-conv}-\eqref{eq:bottomwb1-1uq-conv}. Consider a deterministic initial water surface
\begin{equation}\label{eq:initwb1-2-ten-uq}
    \eta(x,y,0,\xi) = \left\{\begin{aligned}
    &1.01,&&\text{if } 0.05<x<0.15,\\
    &1,&&\text{otherwise},
    \end{aligned}\right.\quad u(x,y,0,\xi) = v(x,y,0,\xi) = 0,
\end{equation}
and a stochastic bottom
\begin{equation}\label{eq:bottomwb1-2-ten-uq}
    B(x,y, \xi) = 0.8e^{-5(x-0.9+0.1\xi^{(1)})^2-50(y-0.5 + 0.1\xi^{(2)})^2},
\end{equation}
where $\xi = (\xi^{(1)},\xi^{(2)})^\top$ is a two-dimensional random variable, with $\xi^{(1)}$ having Beta density with parameters $(\alpha,\beta) = (1,3)$, and $\xi^{(2)}\sim\mathcal{U}(-1,1)$. The randomness is put on the position of the center of the Gaussian-shape ``hump''. 

\begin{figure}[htbp]
    \centering
    \subfloat{\includegraphics[width=0.45\textwidth]{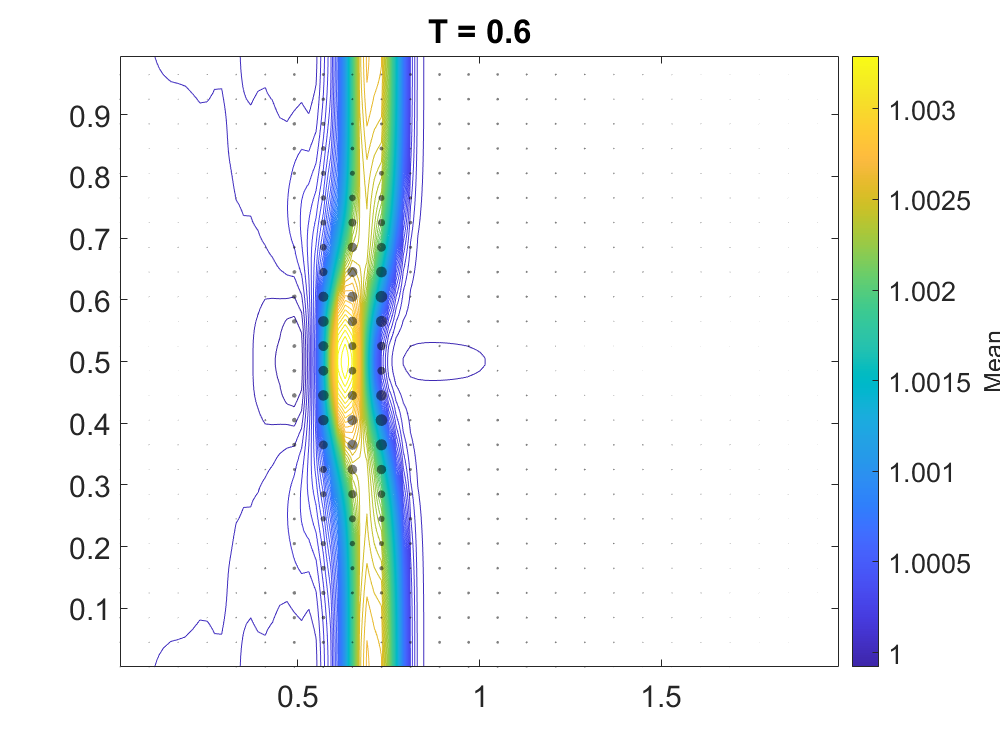}}
    \subfloat{\includegraphics[width=0.45\textwidth]{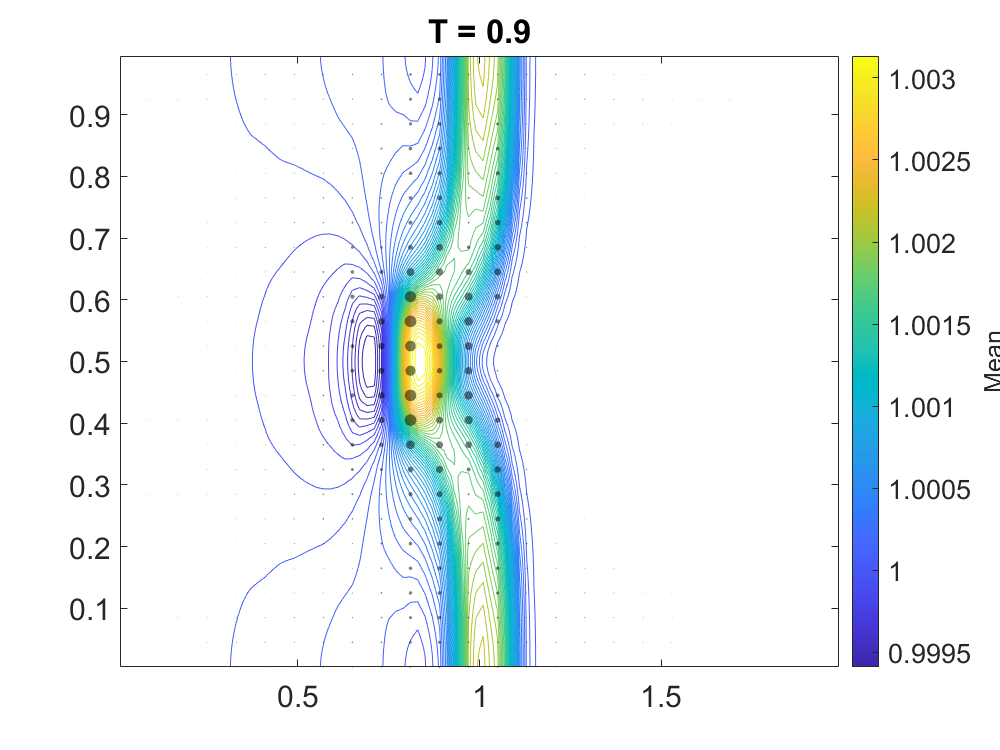}}
    \\
    \subfloat{\includegraphics[width=0.45\textwidth]{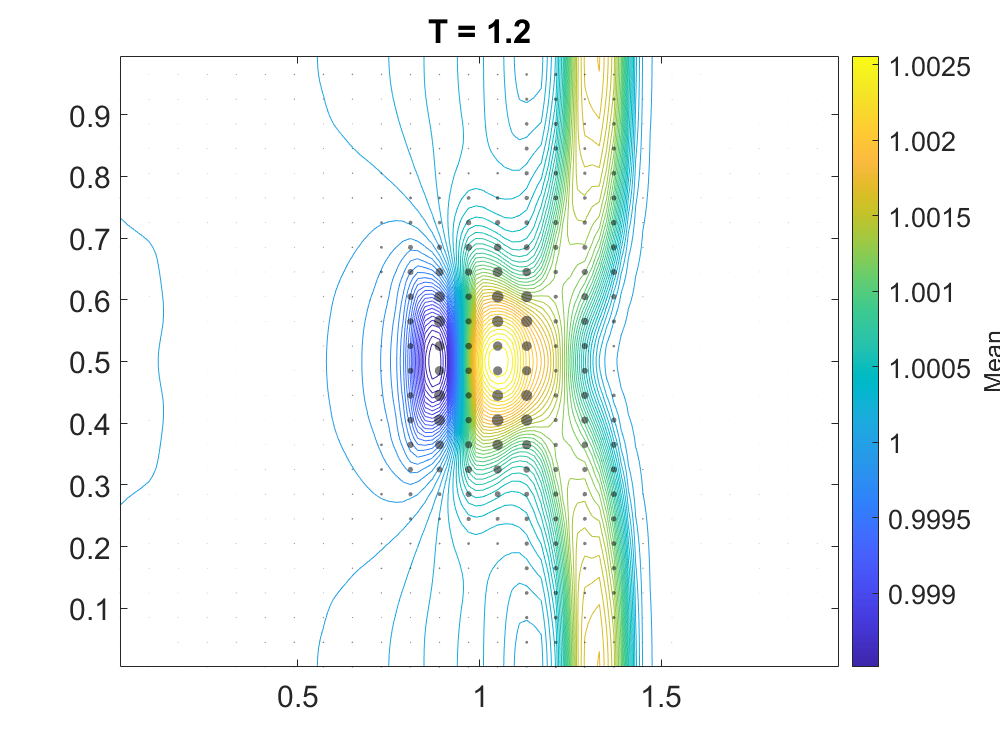}}
    \subfloat{\includegraphics[width=0.45\textwidth]{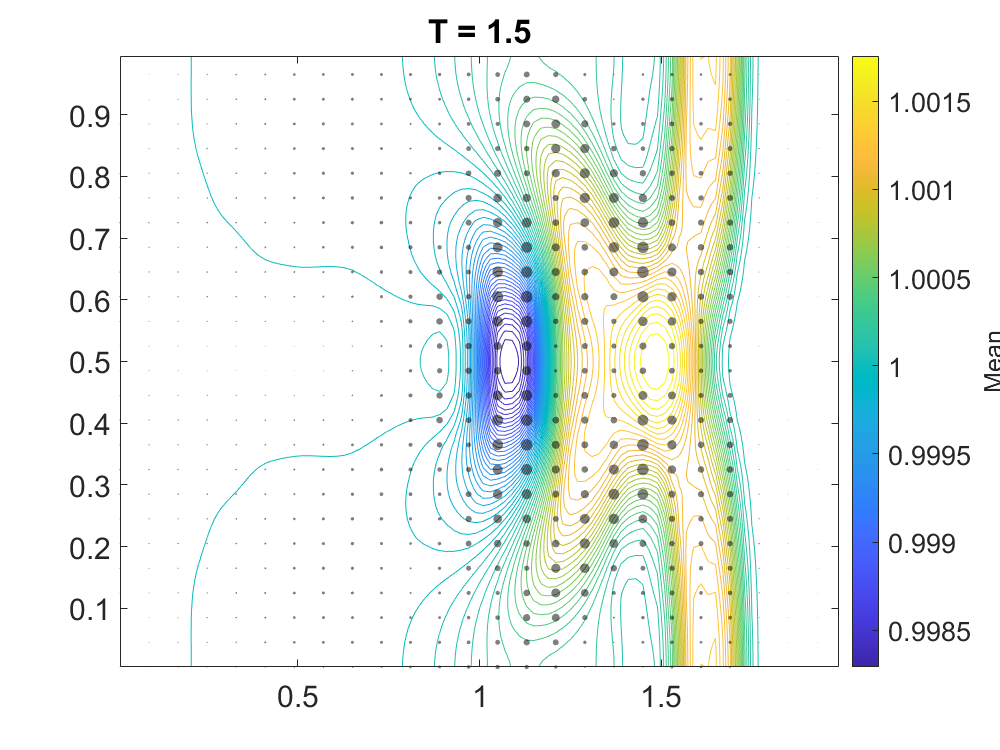}}\\
    \subfloat{\includegraphics[width=0.45\textwidth]{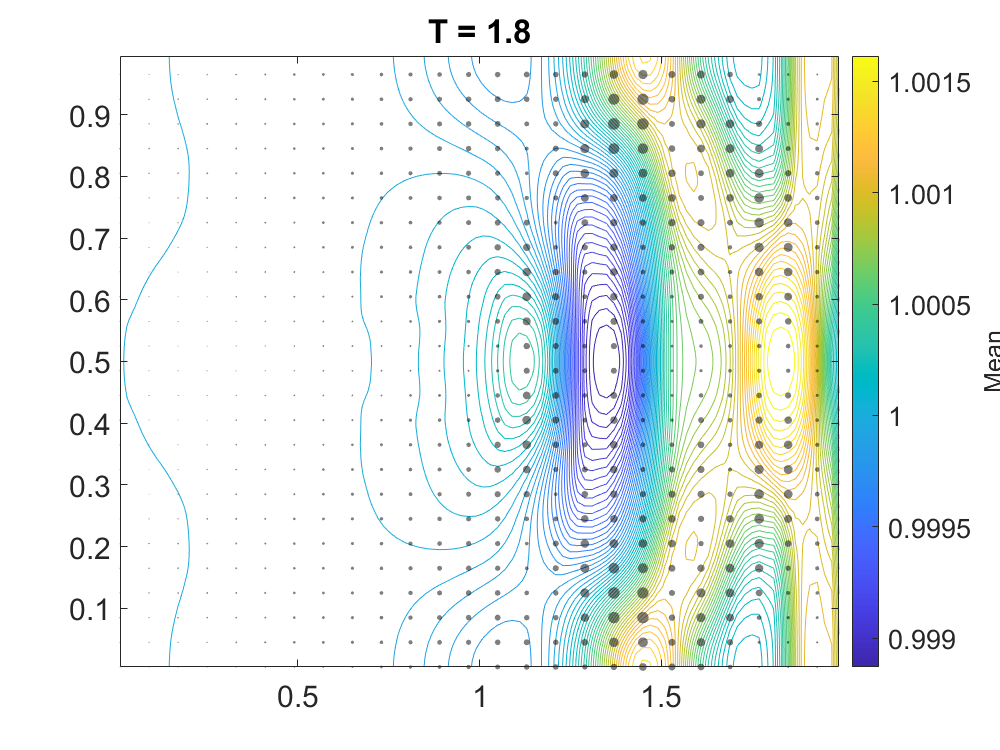}}    
    \caption{Numerical solution to \eqref{eq:initwb1-2-ten-uq}-\eqref{eq:bottomwb1-2-ten-uq}, water surface, $100\times 100$ grid, disk-glyph over mean contours, where the radii of the disks indicate the magnitude of the standard deviation, $t = 0.6, 0.9, 1.2, 1.5, 1.8$ (Top left, top right, middle left, middle right, bottom). The largest disks are corresponding to the standard deviation values 6.97e-4, 1.50e-3, 1.03e-3, 5.16e-4, and 4.81e-4, respectively.}
    \label{fig:wb1-2-ten-100-uqmean}    
\end{figure}

\begin{figure}[htbp]
    \centering
    \subfloat{\includegraphics[width=0.45\textwidth]{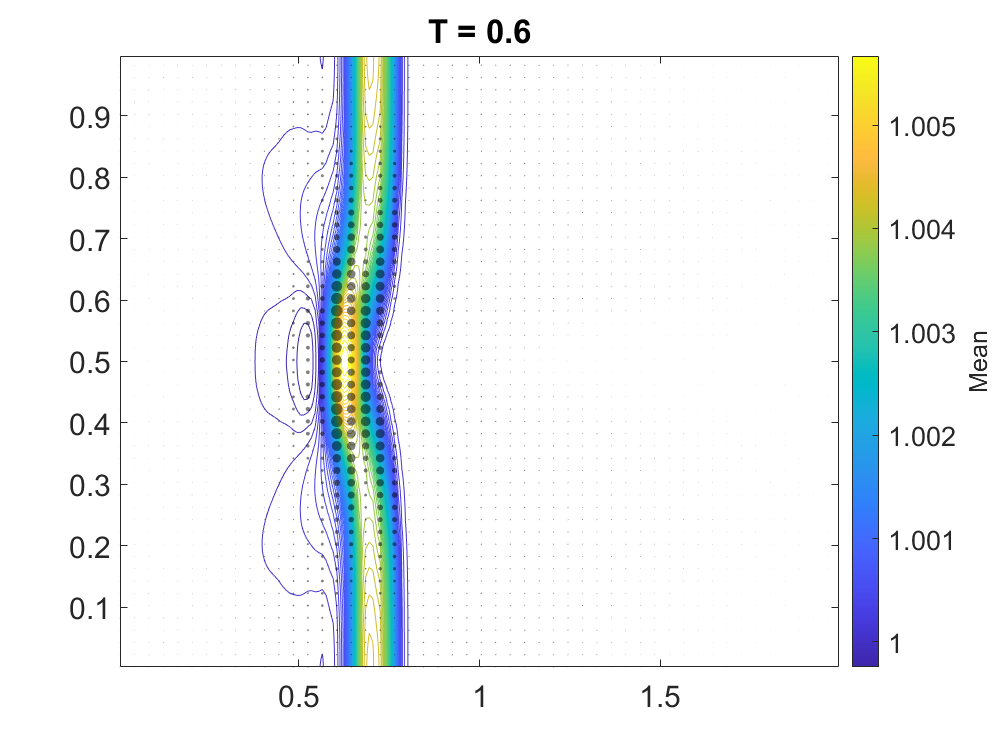}}
    \subfloat{\includegraphics[width=0.45\textwidth]{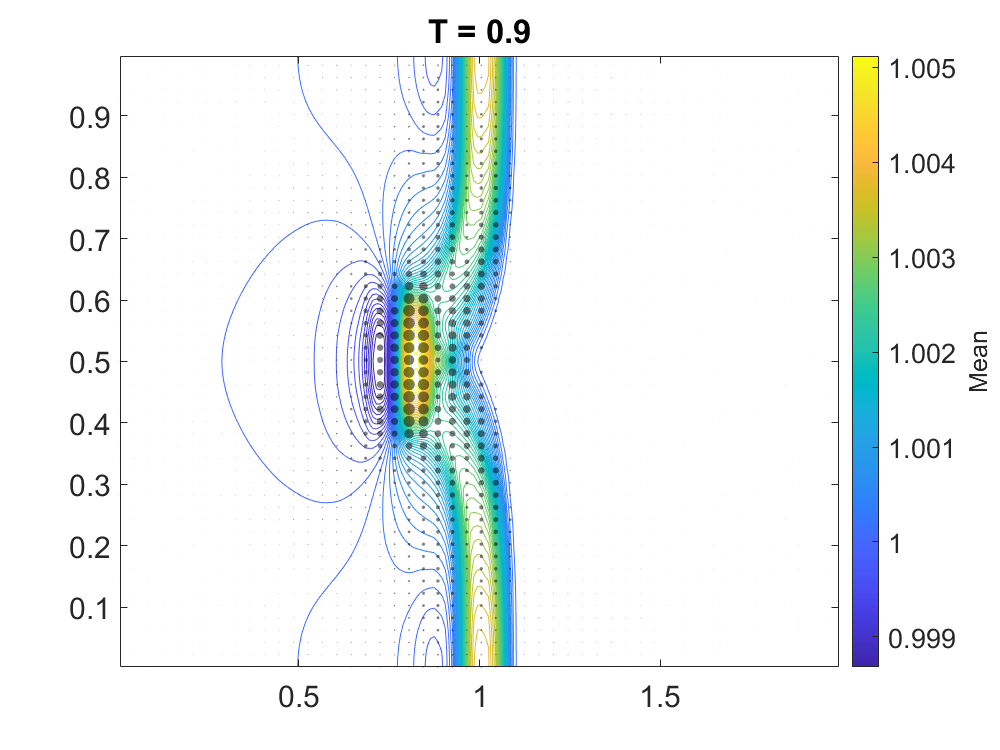}}
    \\
    \subfloat{\includegraphics[width=0.45\textwidth]{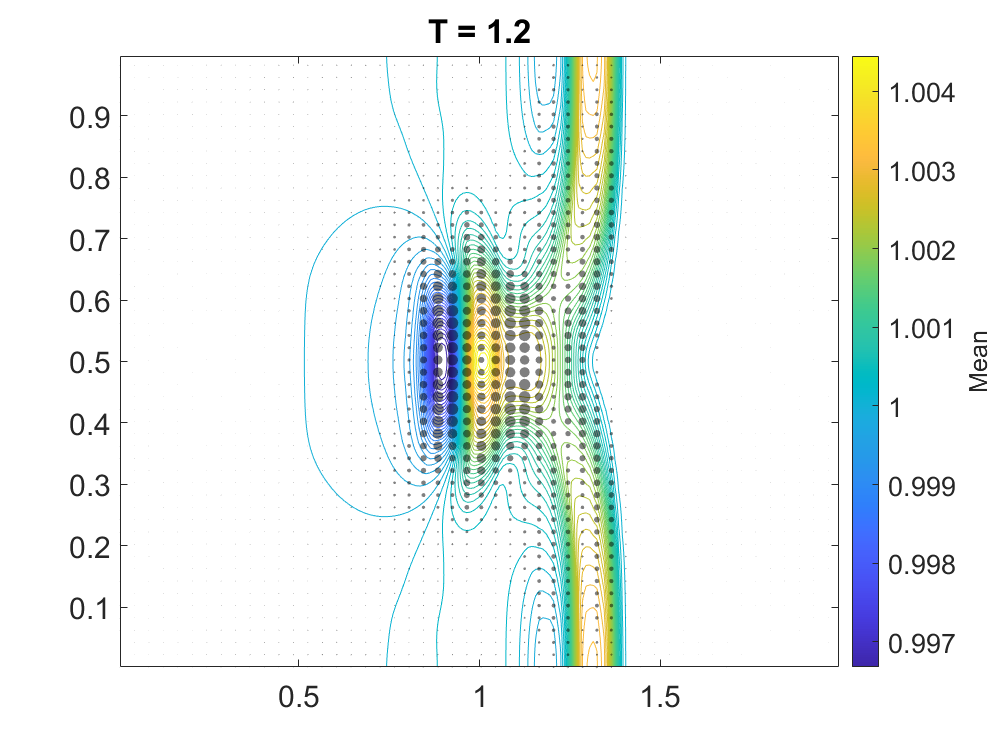}}
    \subfloat{\includegraphics[width=0.45\textwidth]{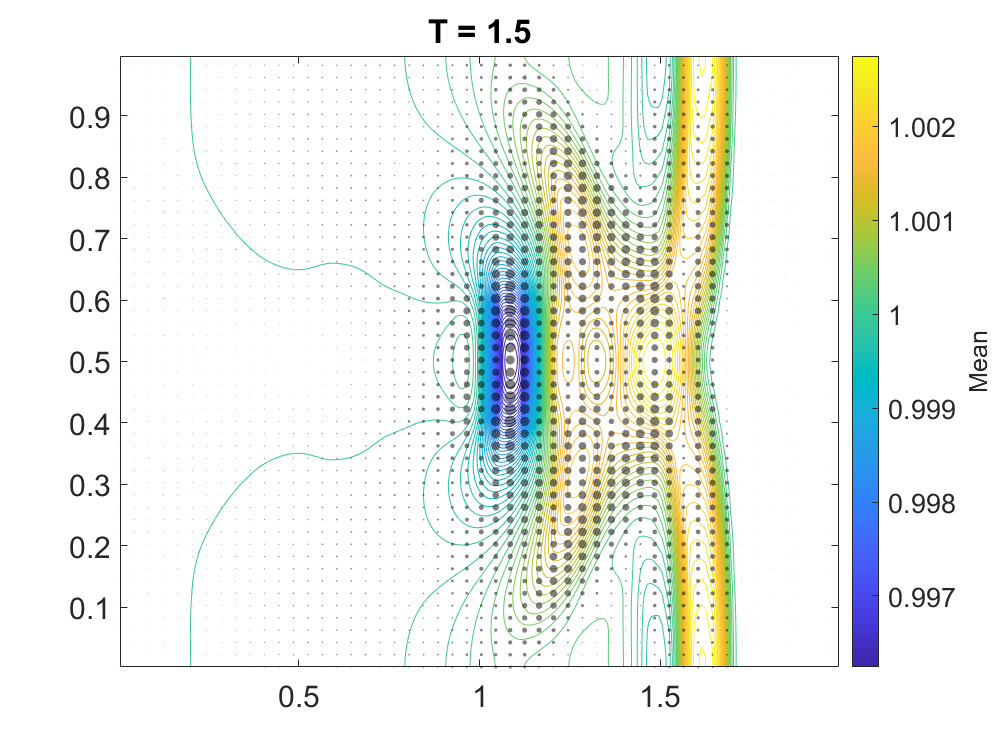}}\\
    \subfloat{\includegraphics[width=0.45\textwidth]{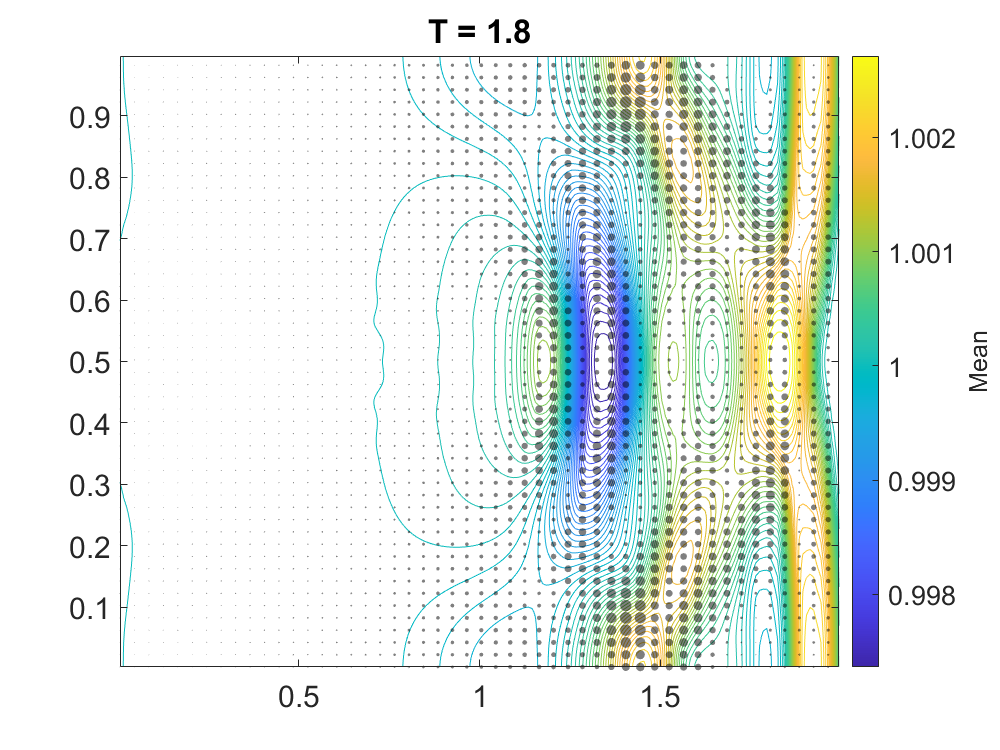}}    
    \caption{Numerical solution to \eqref{eq:initwb1-2-ten-uq}-\eqref{eq:bottomwb1-2-ten-uq}, $200\times 200$ grid, water surface, disk-glyph over mean contours, where the radii of the disks indicate the magnitude of the standard deviation, $t = 0.6, 0.9, 1.2, 1.5, 1.8$ (Top left, top right, middle left, middle right, bottom). The largest disks are corresponding to the standard deviation values 2.20e-3, 3.35e-3, 2.00e-3, 1.58e-3, and 1.20e-3, respectively.}
    \label{fig:wb1-2-ten-200-uqmean}    
\end{figure}

The computational domain is $[0,2]\times [0,1]$. We compute the numerical solution on $100\times 100$ and $200\times 200$ grid at various time. We plot the contours of mean overlaid with disk glyphs (\Cref{fig:wb1-2-ten-100-uqmean} and \Cref{fig:wb1-2-ten-200-uqmean}) whose radius are proportional to the uncertainties at the corresponding cells. We uniformly downsample the cells for the disk glyphs to better present the uncertainties in the figures. The water propagates towards the right in the computational domain and the uncertainties remain small for most of the sampled cells at the beginning ($t=0.6$). As the water interacts with the hump, the water splits and propagates along all directions. The contours of the mean water surface generate more complicated structures and the uncertainties propagate along the directions of the mean waves. The largest uncertainties were concentrated near the peak of the water surface ($t=0.6,0.9$) but then gradually ``spread out'' following the mean waves ($t=1.2,1.5,1.8$). One can also compare the middle left contours (corresponding to the solution at $t = 1.2$) in \Cref{fig:wb1-2-ten-200-uqmean} or  \Cref{fig:wb1-2-ten-100-uqmean} with the contours of the mean water surfaces in \Cref{fig:wb1-1-uq-mean-conv}. The uncertainties imposed in the position shifts the peak of the mean towards the right. It is reasonable because $\xi^{(1)}$ follows a distribution that skews negatively so the mean center of the hump $(0.9-0.1\xi^{(1)})$ is skews towards the right. On the other hand, the uncertainties on the center increase the width of the peak and valley structures.

\subsection{Example 5: A Submerged Flat Plateau}\label{sect:sub-island}
\begin{figure}[htbp]
    \centering
     \subfloat{\includegraphics[width=0.45\textwidth]{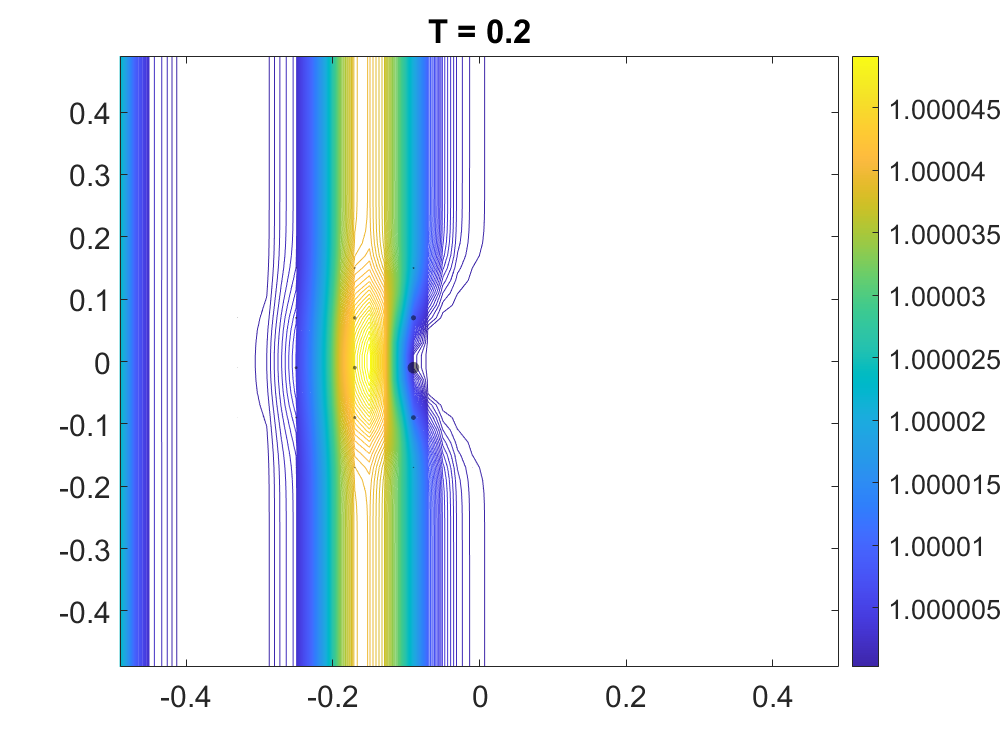}}
    \subfloat{\includegraphics[width=0.45\textwidth]{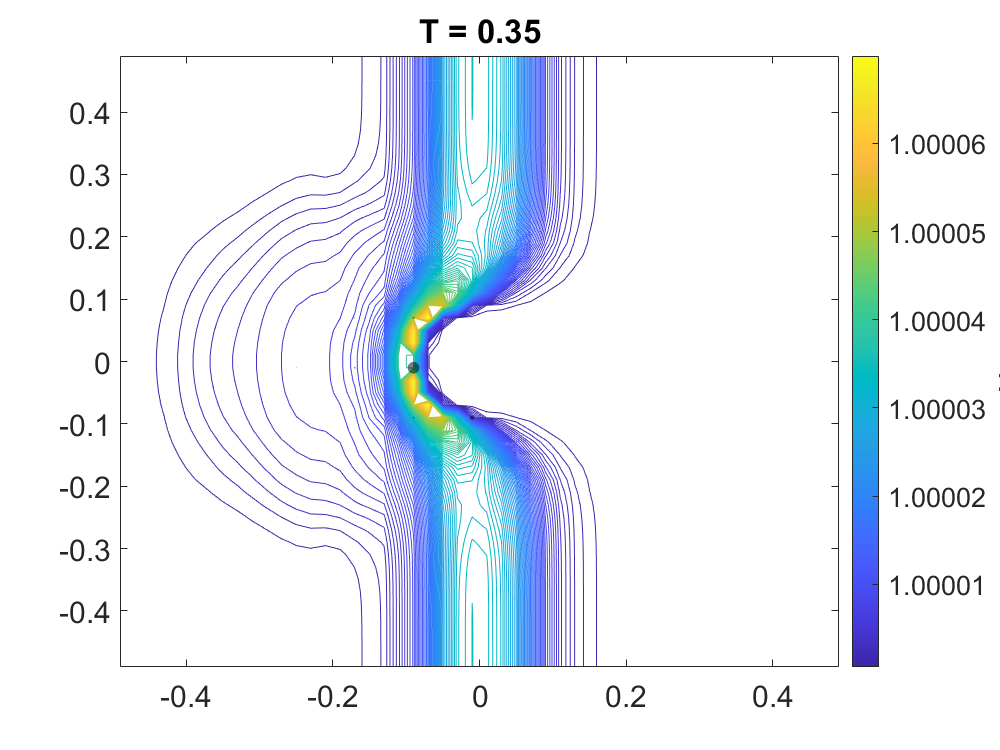}}
    \\
    \subfloat{\includegraphics[width=0.45\textwidth]{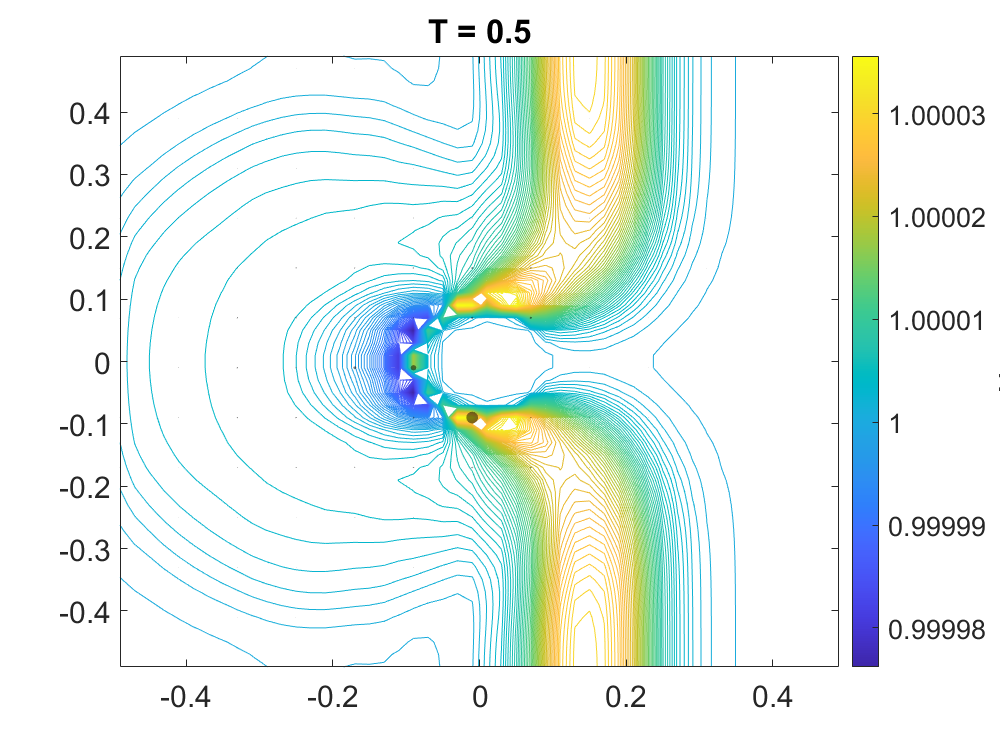}}
    \subfloat{\includegraphics[width=0.45\textwidth]{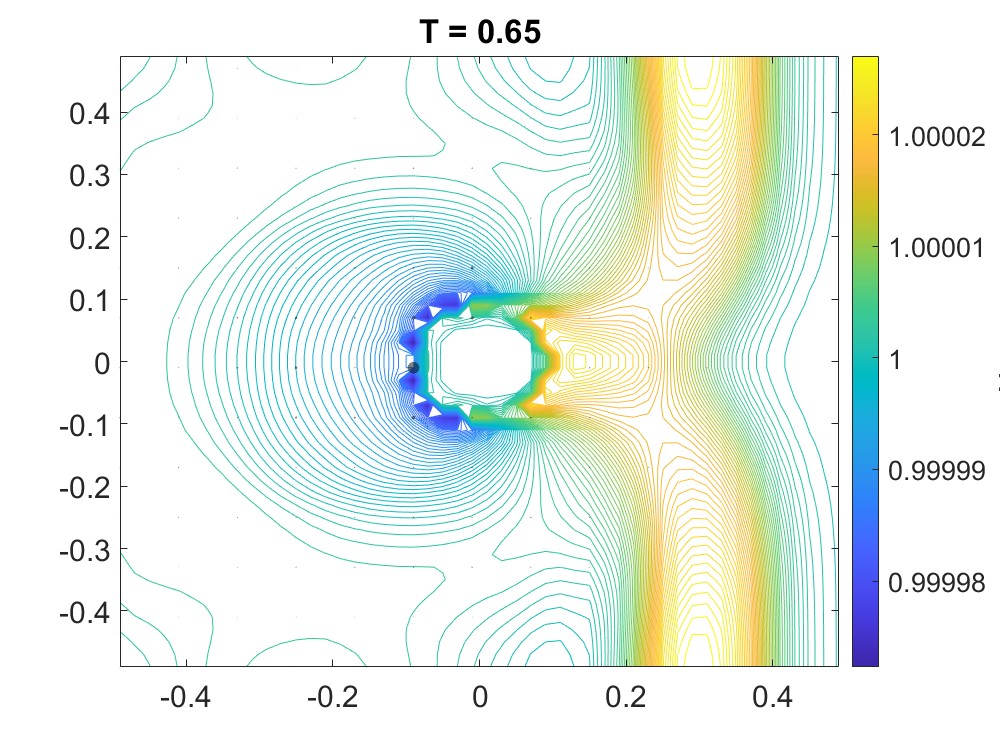}}\\
    \caption{Numerical solution to \eqref{eq:initwb3-4-tenuq}-\eqref{eq:bottomwb3-4-tenuq}, $50\times 50$, water surface, disk-glyph over mean contours, where the radii of the disks indicate the magnitude of the standard deviation, $t = 0.2, 0.35, 0.5, 0.65$ (top left, top right, middle left, middle right). The largest disks are corresponding to the standard deviation values 4.87e-7, 6.28e-7, 8.45e-7, and 5.37e-8, respectively.}
    \label{fig:wb3-4-ten-50-uqmean}    
\end{figure}

\begin{figure}[htbp]
    \centering
     \subfloat{\includegraphics[width=0.45\textwidth]{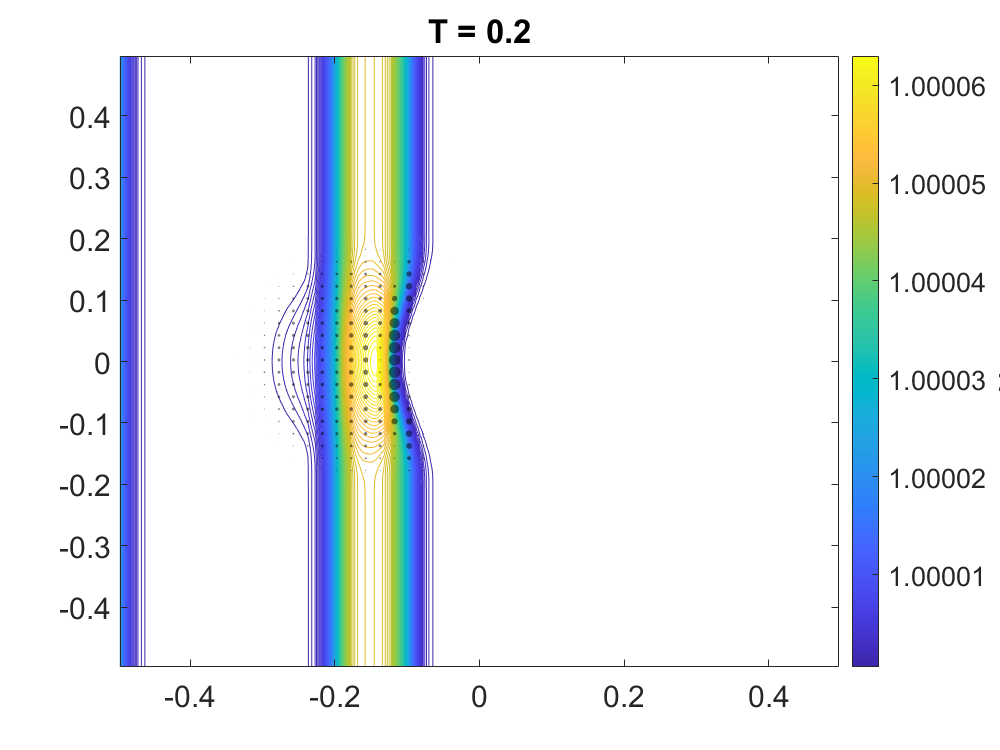}}
    \subfloat{\includegraphics[width=0.45\textwidth]{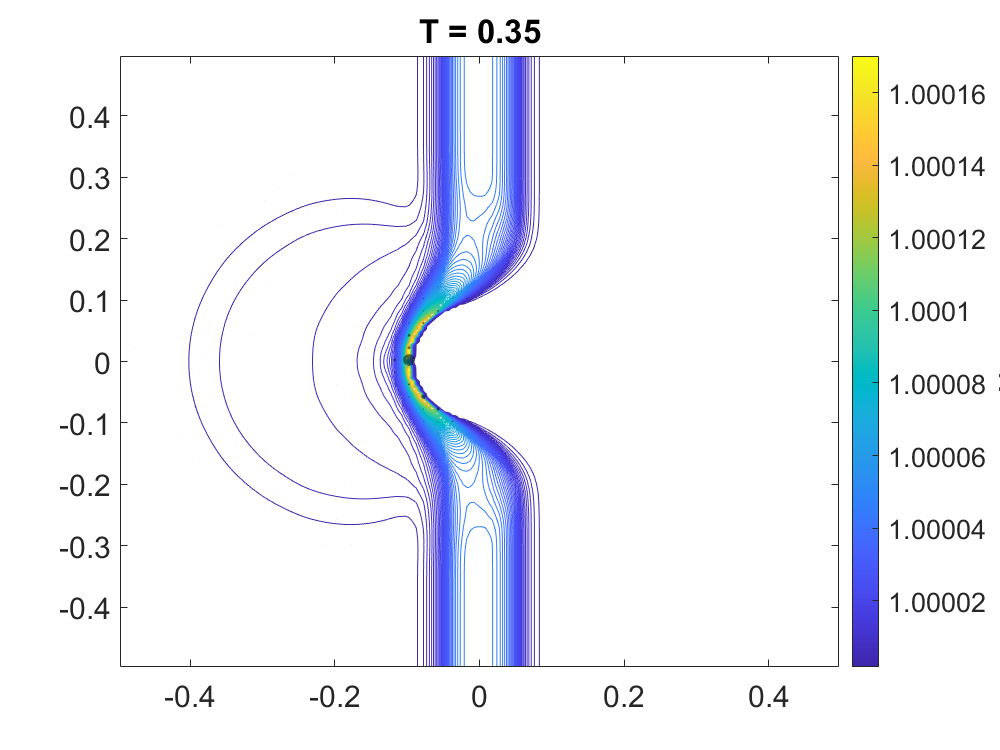}}
    \\
    \subfloat{\includegraphics[width=0.45\textwidth]{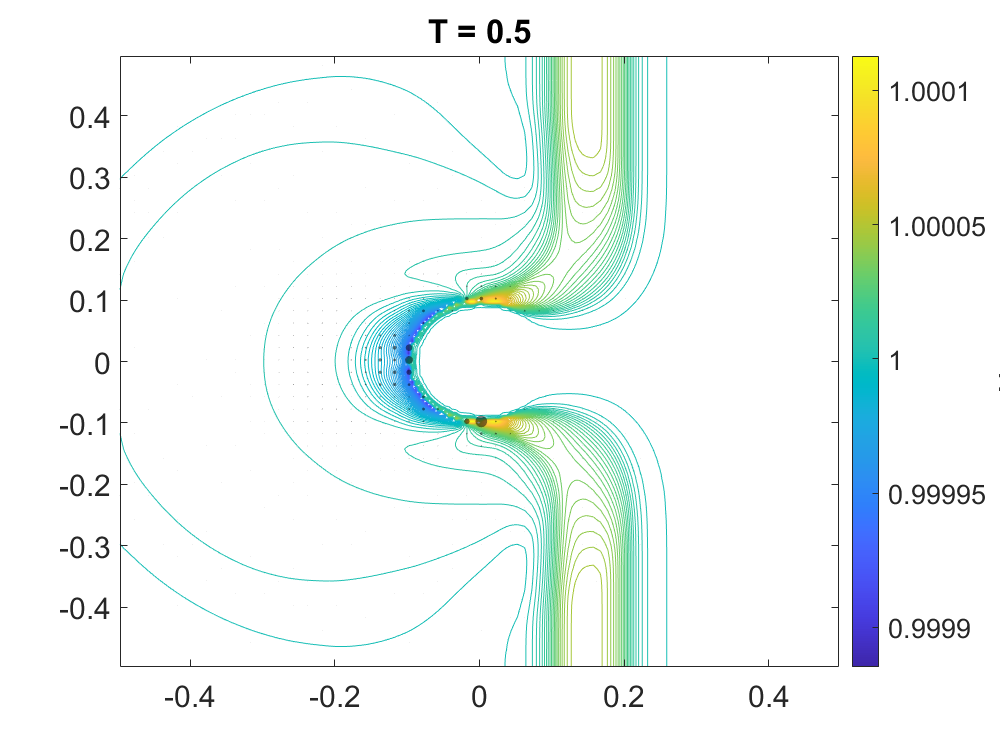}}
    \subfloat{\includegraphics[width=0.45\textwidth]{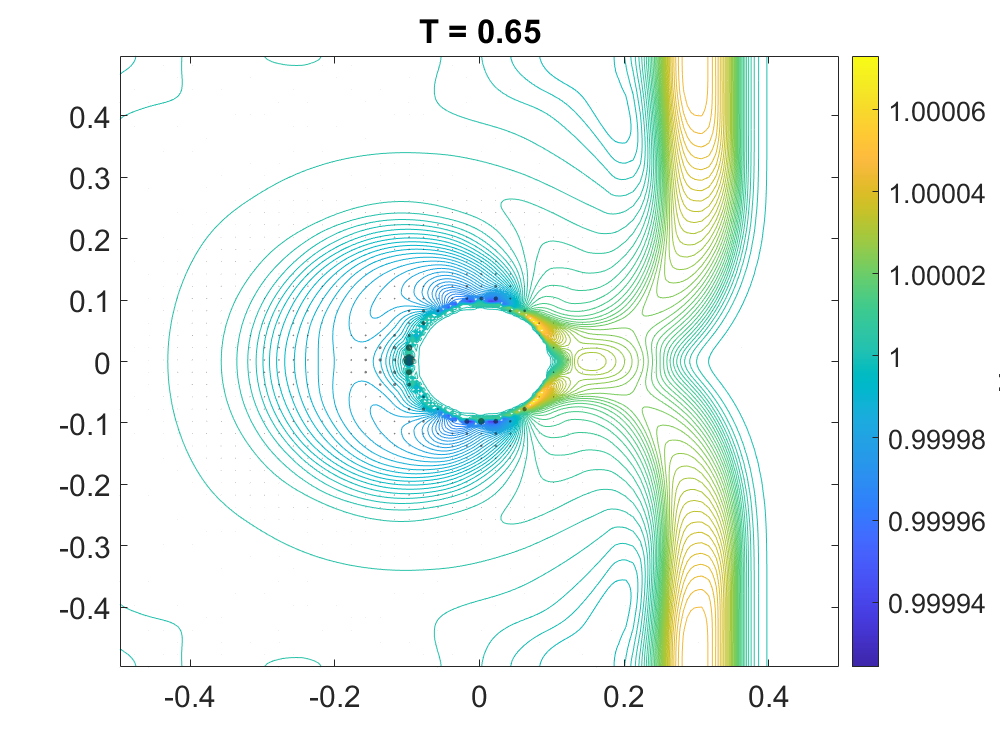}}\\
    \caption{Numerical solution to \eqref{eq:initwb3-4-tenuq}-\eqref{eq:bottomwb3-4-tenuq}, $200\times 200$, water surface, disk-glyph over mean contours, where the radii of the disks indicate the magnitude of the standard deviation, $t = 0.2, 0.35, 0.5, 0.65$ (top left, top right, middle left, middle right). The largest disks are corresponding to the standard deviation values 4.83e-8, 6.57e-6, 4.20e-6, and 3.50e-6, respectively.}
    \label{fig:wb3-4-ten-200-uqmean}    
\end{figure}

We consider a deterministic initial water surface with deterministic perturbation
\begin{equation}\label{eq:initwb3-4-tenuq}
    \eta(x,y,0) = \left\{\begin{aligned}
    &1.0001,&&\text{if } -0.4<x<-0.3,\\
    &1,&&\text{otherwise},
    \end{aligned}\right.,\quad u(x,y,0) = v(x,y,0) = 0,
\end{equation}
and with stochastic bottom topography,
\begin{equation}\label{eq:bottomwb3-4-tenuq}
    B(x,y,\xi) = \left\{\begin{aligned}
    &0.9998,&& r\le 0.1,\\
    &9.997(0.2-r)+0.0001(\xi^{(1)}+1),&& 0.1<r\le0.2,\\
    &0.0001,&&\text{otherwise},
    \end{aligned}\right.
\end{equation}
where $r \coloneqq\sqrt{x^2+y^2}+0.0001(\xi^{(2)}+1)$. The bottom is a flat plateau beneath, but very close to, the water surface. It is continuous but not smooth in the physical domain. The uncertainties $\xi^{(1)} \sim \mathcal{U}(-1,1)$ are on the heights of the bottom in the transition region between the plateau and the lowland in the water. The uncertainties $\xi^{(2)}$ (with Beta density $(\alpha,\beta) = (1,3)$) are on the radius of the disk-shaped plateau. We compute the solutions at various time on the computational domain $[-0.5,0.5]\times[-0.5,0.5]$ with a $50\times 50$ and $200\times 200$ grid. The mean contour plots with disk glyphs are presented in \Cref{fig:wb3-4-ten-50-uqmean} and \Cref{fig:wb3-4-ten-200-uqmean}. The initial deterministic perturbation is split into two waves. The left-going wave propagates out of the computational region while the right-going wave interacts with the plateau and generates complicated waves as the interaction starts. The initial deterministic perturbation on the water surface is $10^{-4}$ and largest level of the mean water displacement from the steady state stay in the same order during the computation. The uncertainties (standard deviation) remain small and are close to the island, which shows the capability of our scheme to capture the stochastic ``lake-at-rest'' solution \eqref{eq:lake-at-rest}.

\begin{figure}[htbp]
    \centering
     \subfloat{\includegraphics[width=0.45\textwidth]{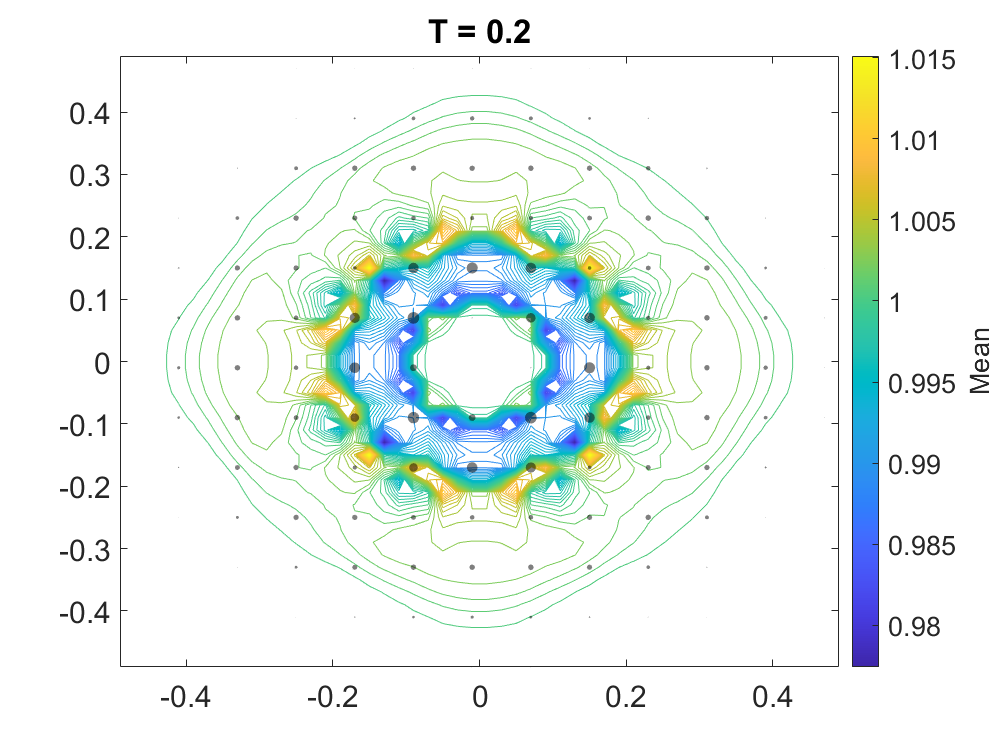}}
    \subfloat{\includegraphics[width=0.45\textwidth]{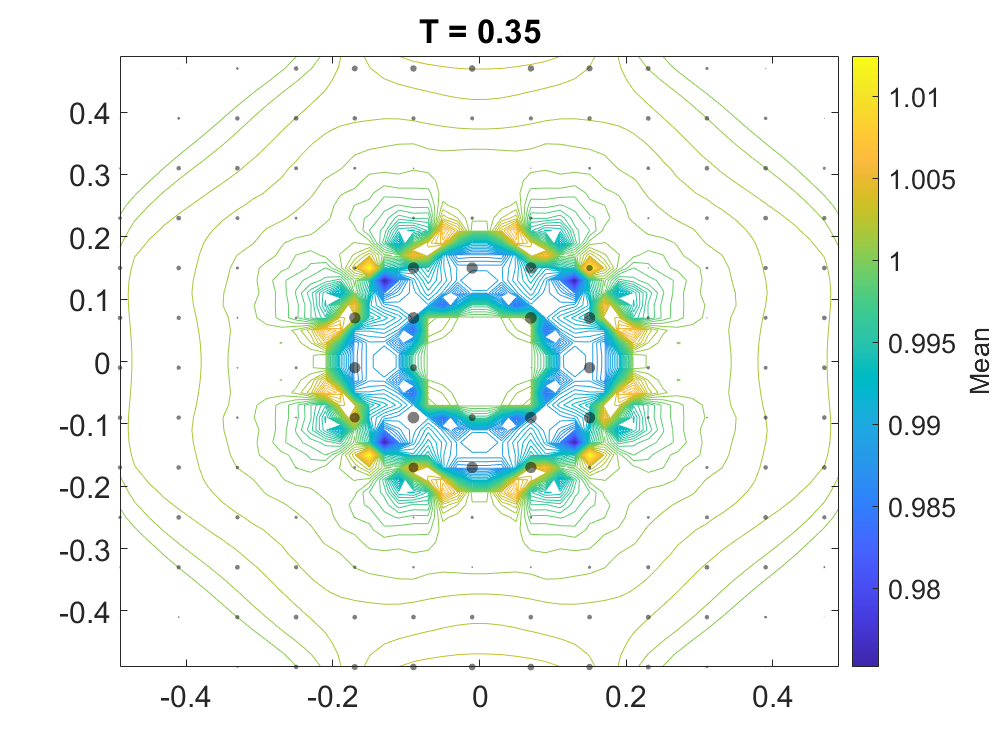}}
    \\
    \subfloat{\includegraphics[width=0.45\textwidth]{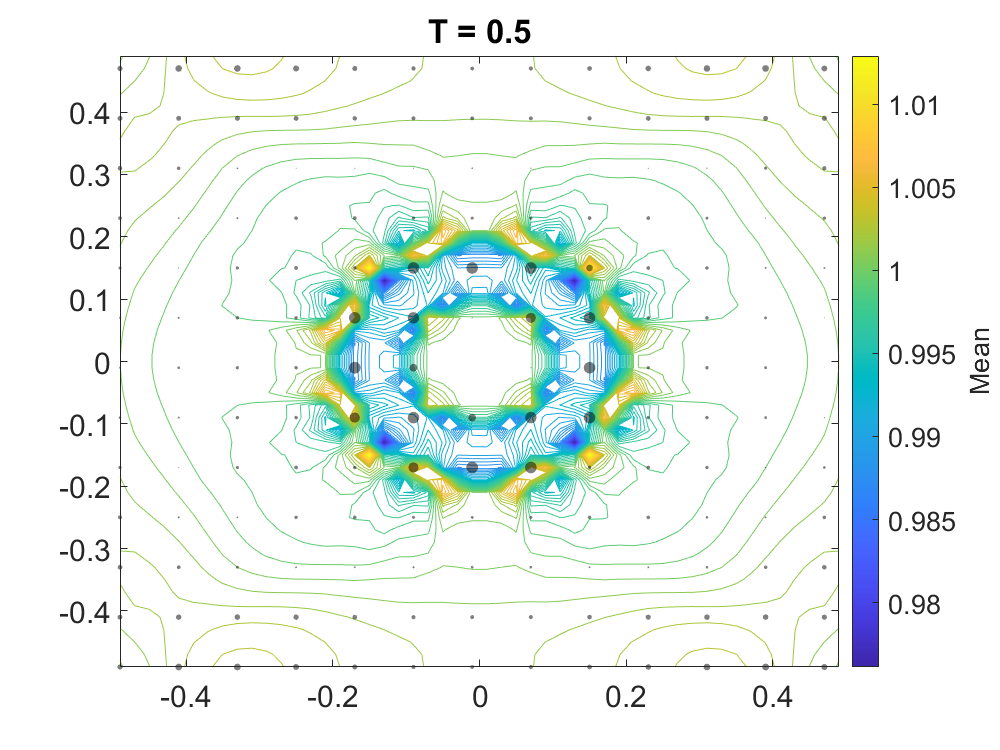}}
    \subfloat{\includegraphics[width=0.45\textwidth]{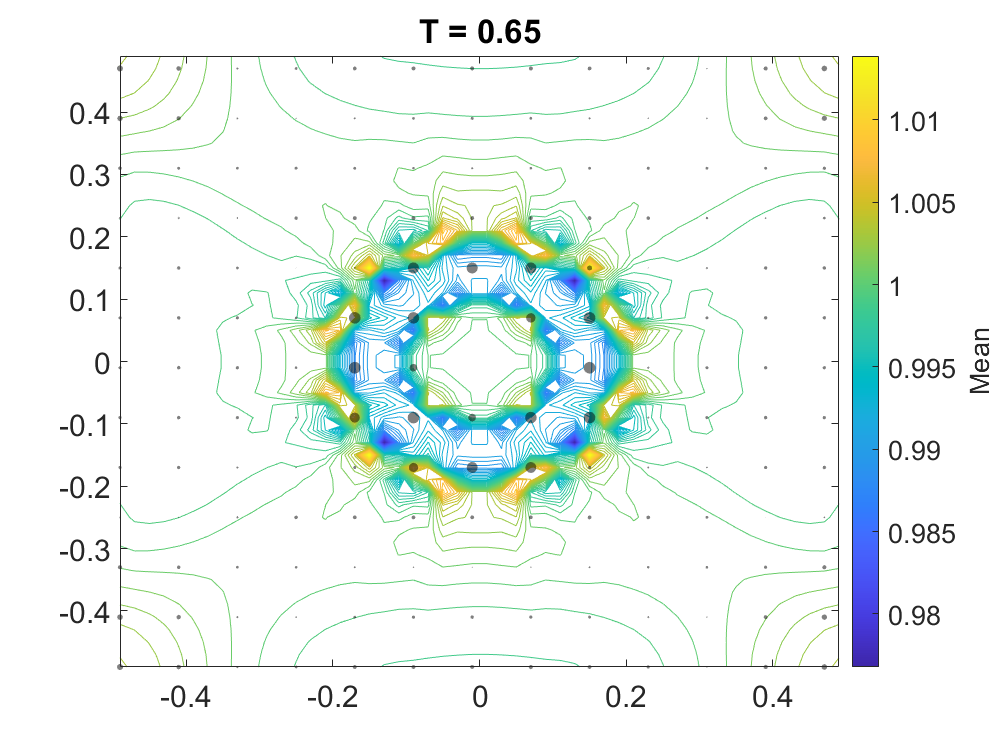}}\\
    \caption{Numerical solution to \eqref{eq:initwb3-4-tenuq}-\eqref{eq:bottomwb3-4-tenuq} for non-well-balanced scheme, $50\times 50$, water surface, disk-glyph over mean contours, where the radii of the disks indicate the magnitude of the standard deviation, $t = 0.2, 0.35, 0.5, 0.65$ (top left, top right, middle left, middle right). The largest disks are corresponding to the standard deviation values 4.38e-4, 4.45e-4, 3.92e-4, and 3.73e-6, respectively. The well-balanced counterpart is \Cref{fig:wb3-4-ten-50-uqmean}.}
    \label{fig:nwb3-4-ten-50-uqmean}   
\end{figure}

\begin{figure}[htbp]
    \centering
     \subfloat{\includegraphics[width=0.45\textwidth]{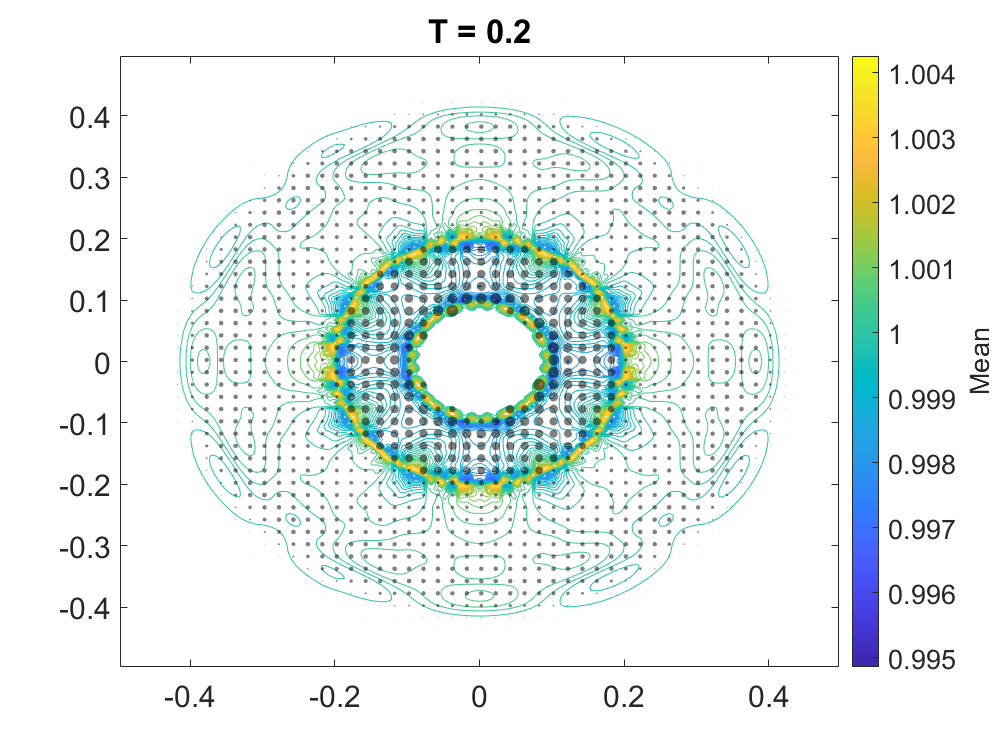}}
    \subfloat{\includegraphics[width=0.45\textwidth]{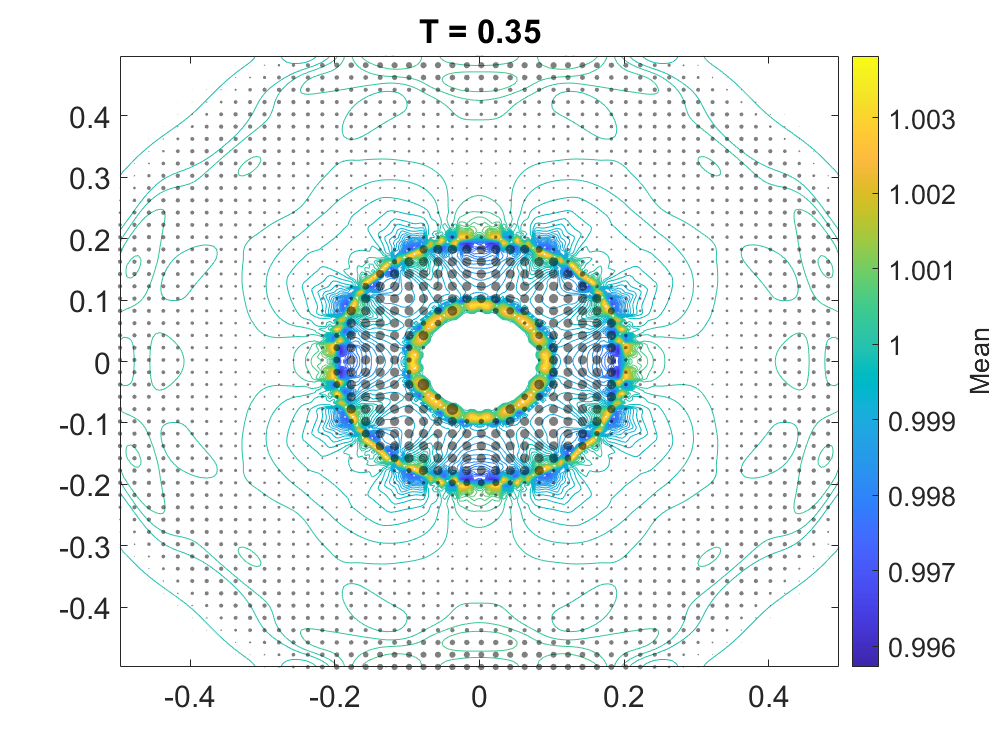}}
    \\
    \subfloat{\includegraphics[width=0.45\textwidth]{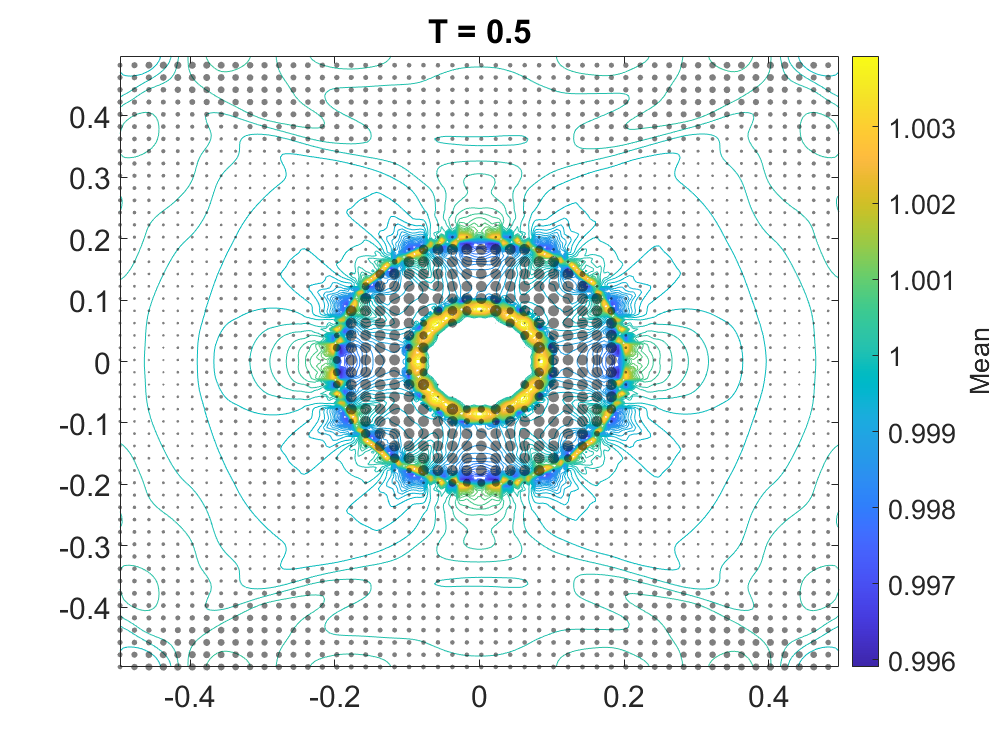}}
    \subfloat{\includegraphics[width=0.45\textwidth]{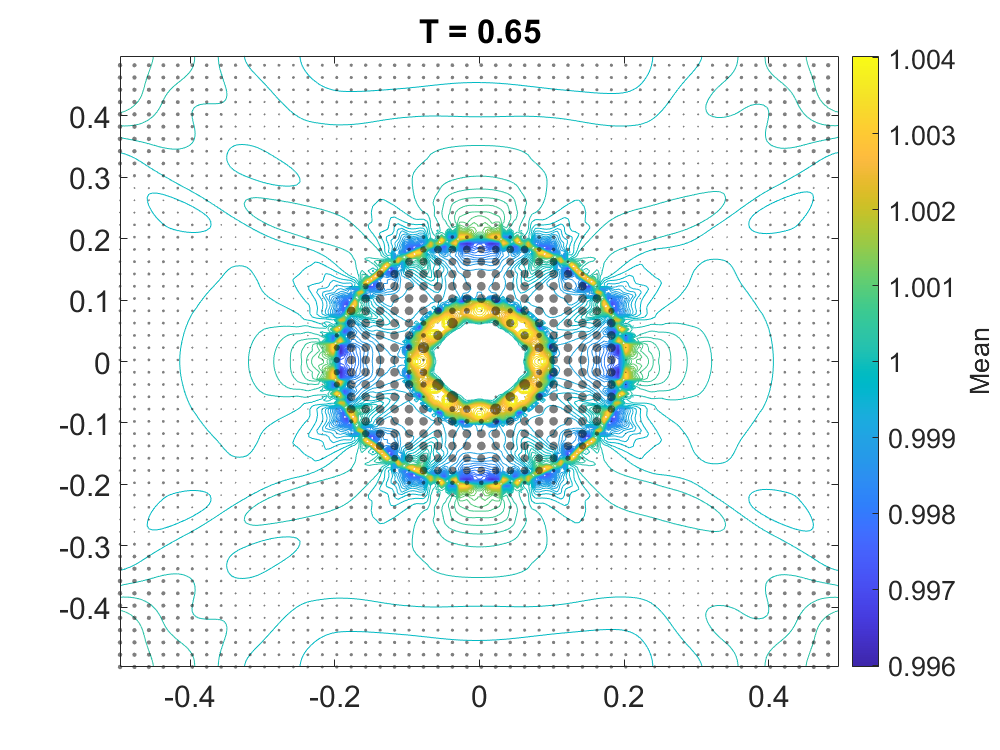}}\\
   \caption{Numerical solution to \eqref{eq:initwb3-4-tenuq}-\eqref{eq:bottomwb3-4-tenuq} for non-well-balanced scheme, $200\times 200$, water surface, disk-glyph over mean contours, where the radii of the disks indicate the magnitude of the standard deviation, $t = 0.2, 0.35, 0.5, 0.65$ (top left, top right, middle left, middle right). The largest disks are corresponding to the standard deviation values 7.07e-4, 6.56e-4, 4.22e-4, and 6.60e-4, respectively. The well-balanced counterpart is \Cref{fig:wb3-4-ten-200-uqmean}.}
    \label{fig:nwb3-4-ten-200-uqmean}   
\end{figure}

As a comparison, we also include the results of a non-well-balanced scheme, where the numerical source term is given by a straightforward midpoint quadrature rule:
  \begin{equation}\label{eq:nwbsource}
  \left\{
    \begin{aligned}
      \overline{\boldsymbol{S}}^{(1)}_{i,j} &= \boldsymbol{0},\\
      \overline{\boldsymbol{S}}^{(2)}_{i,j} &=  -g\mathcal{P}(\overline{\boldsymbol{h}}_{i,j})(\boldsymbol{B}_x)_{i,j},\\
      \overline{\boldsymbol{S}}^{(3)}_{i,j} &=  -g\mathcal{P}(\overline{\boldsymbol{h}}_{i,j})(\boldsymbol{B}_y)_{i,j},
    \end{aligned}\right.
\end{equation}
where $(\boldsymbol{B}_x)_{i,j}$ and $(\boldsymbol{B}_y)_{i,j}$ are the partial derivatives of the bottom function $\widehat{B}$ at the center of cell $\mathcal{C}_{i,j}$. The numerical results of this non-well-balanced scheme are presented in \Cref{fig:nwb3-4-ten-50-uqmean} and \Cref{fig:nwb3-4-ten-200-uqmean}. As we can observe, the perturbations in the non-well-balanced scheme are far larger than the corresponding well-balanced results. Refining the mesh can reduce the magnitude of the oscillation. However, we cannot observe the correct wave propagation from \Cref{fig:nwb3-4-ten-50-uqmean} and \Cref{fig:nwb3-4-ten-200-uqmean}. The uncertainty exercised on the radius of the island seems to be the major cause of the unphysical results since it affects the ``width'' of the transition region.

\subsection{Example 6: Another Two-Dimensional Example}\label{eq:ex1-2d-1}
\begin{figure}[htbp]
    \centering
    \subfloat{\includegraphics[width=0.45\textwidth]{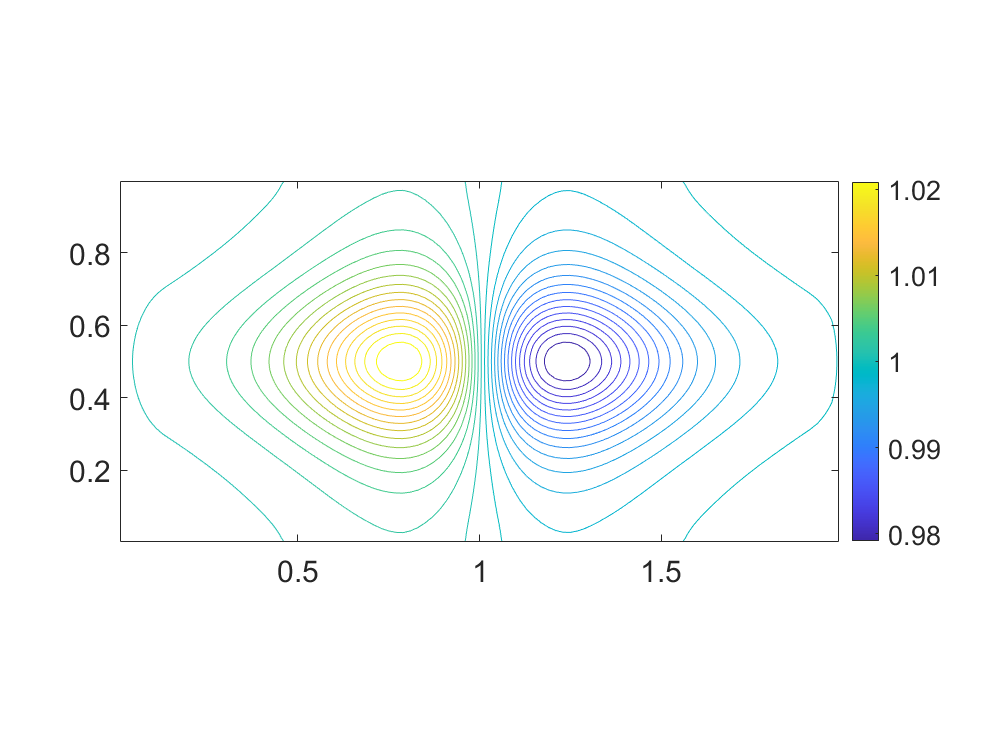}}
    \subfloat{\includegraphics[width=0.45\textwidth]{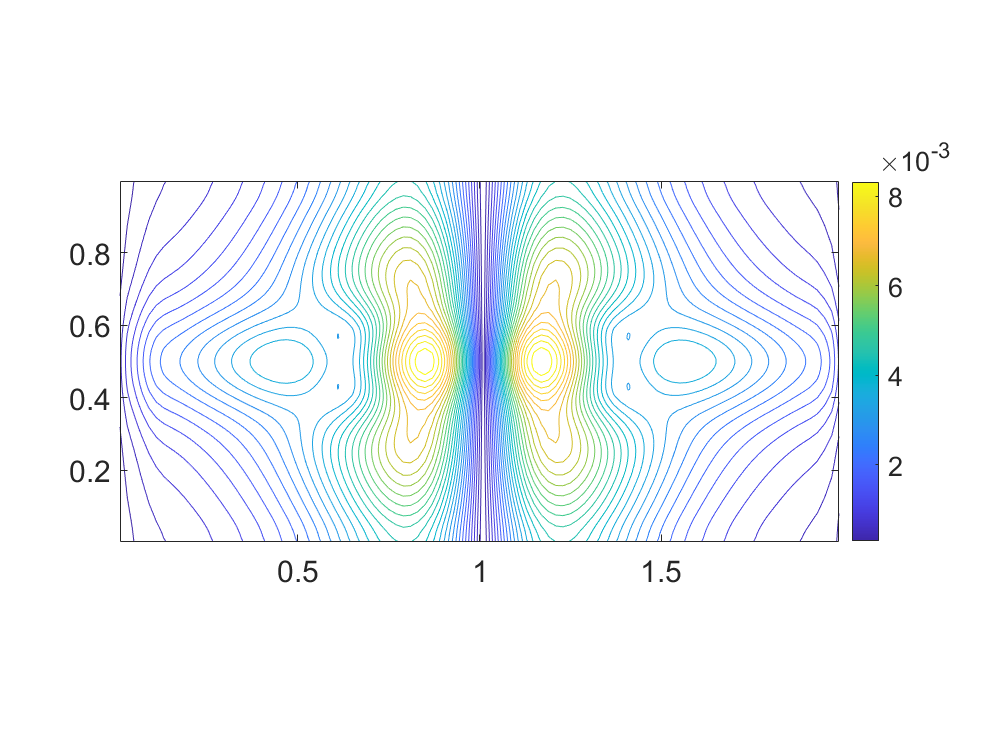}}
    \\
    \subfloat{\includegraphics[width=0.7\textwidth]{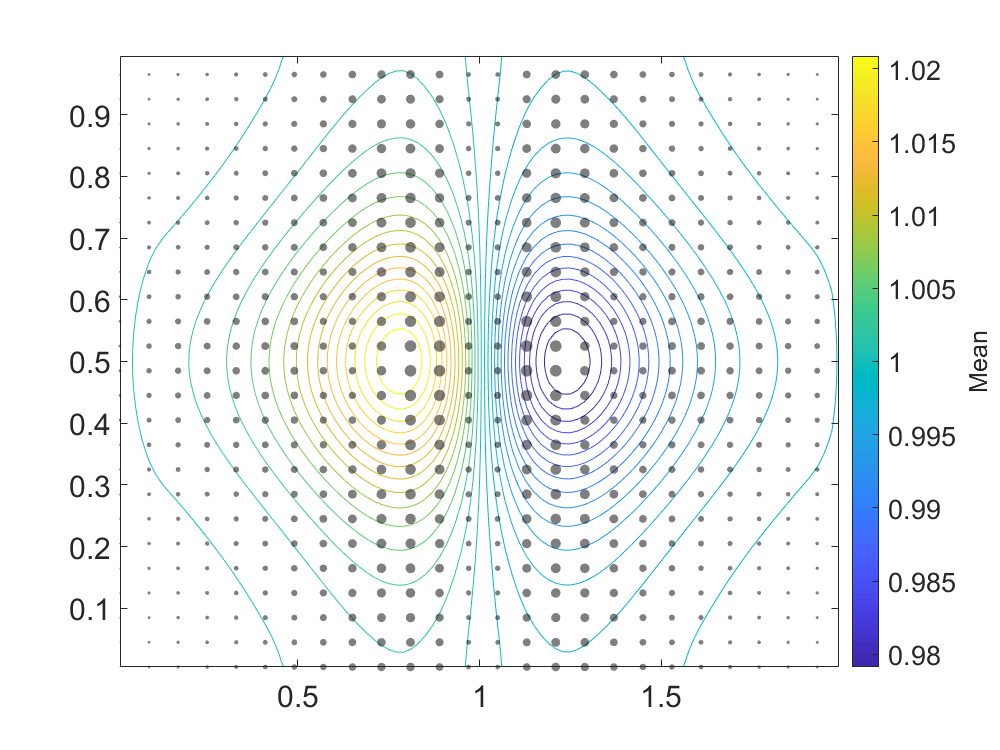}}
    \caption{Numerical solution to \eqref{eq:initex1-3-ten-uq}-\eqref{eq:bottomex1-3-ten-uq}, water surface $t = 0.07$, $100\times 100$ grid. Top left: mean water surface; Top right: the standard deviation of water surface; Bottom: disk-glyph over mean contours, where the radii of the disks indicate the magnitude of the standard deviation. The largest disk is corresponding to the standard deviation 8.10e-3.}
    \label{fig:ex1-3-ten-100-uq}    
\end{figure}

\begin{figure}[htbp]
    \centering
    \subfloat{\includegraphics[width=0.45\textwidth]{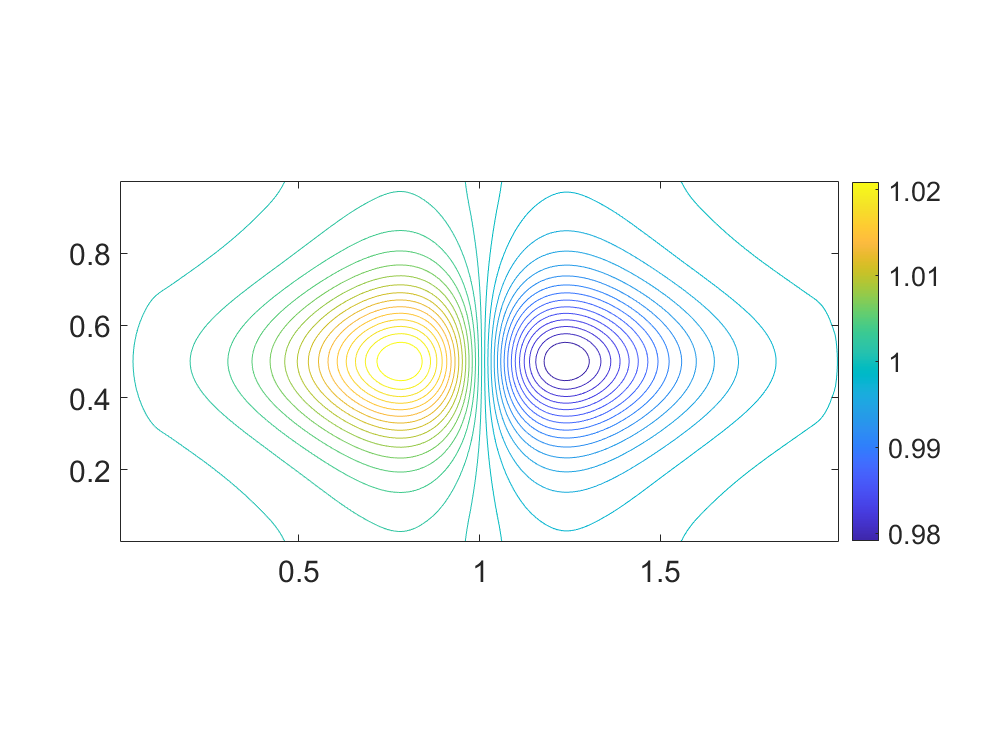}}
    \subfloat{\includegraphics[width=0.45\textwidth]{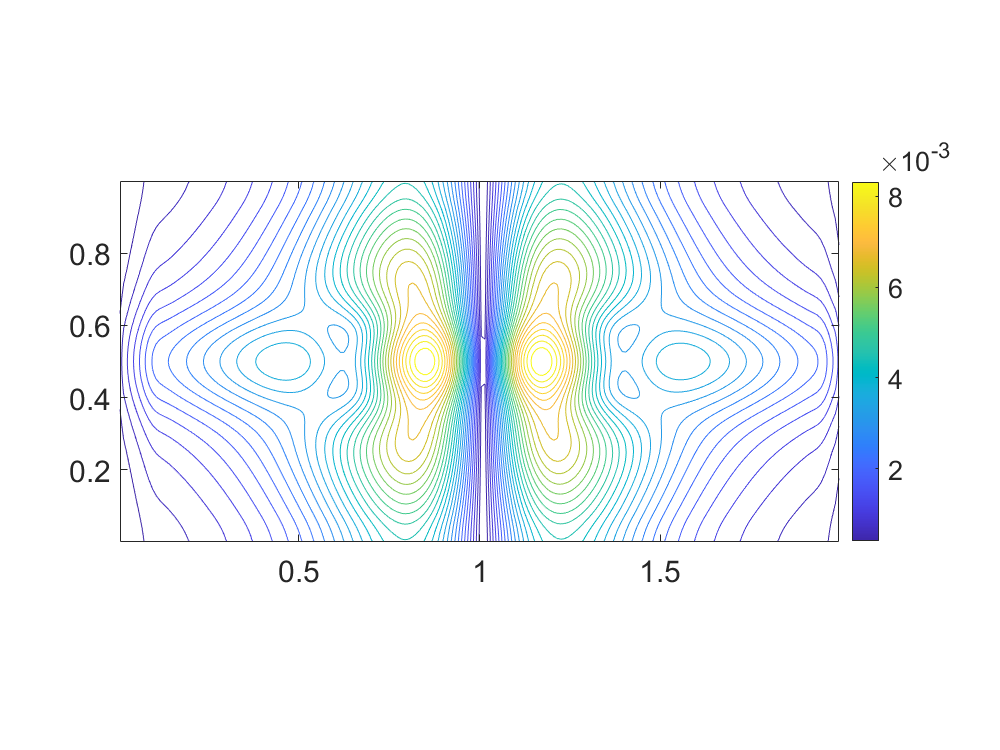}}
    \\
    \subfloat{\includegraphics[width=0.7\textwidth]{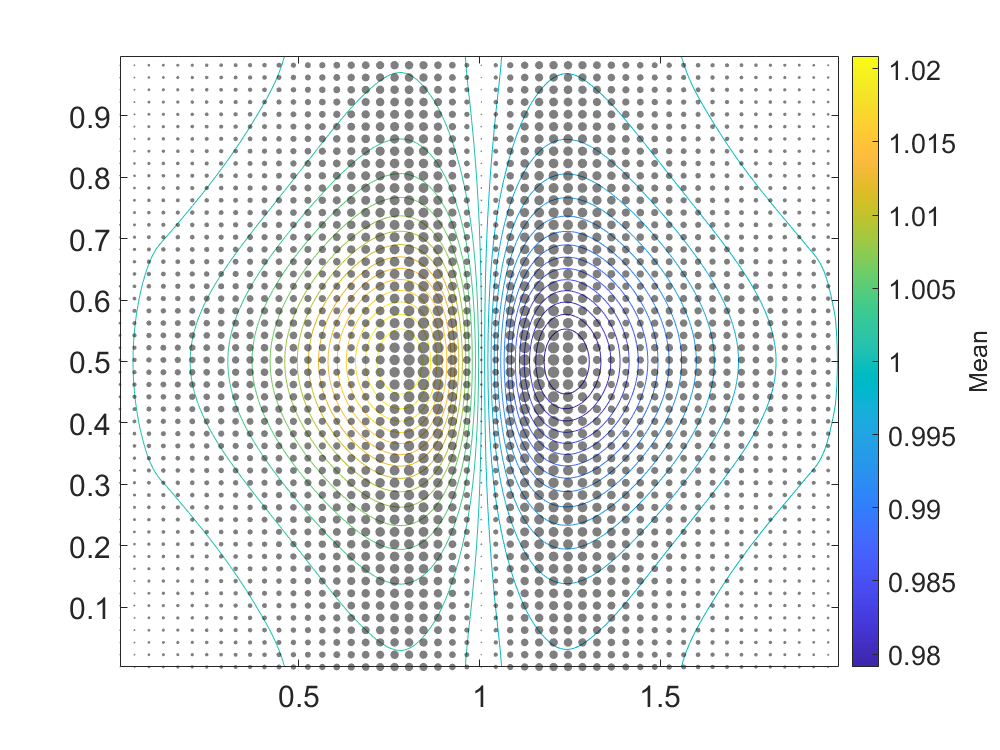}}
    \caption{Numerical solution to \eqref{eq:initex1-3-ten-uq}-\eqref{eq:bottomex1-3-ten-uq}, water surface $t = 0.07$, $200\times 200$ grid. Top left: mean water surface; Top right: the standard deviation of water surface; Bottom: disk-glyph over mean contours, where the radii of the disks indicate the magnitude of the standard deviation. The largest disk is corresponding to the standard deviation 8.58e-3.}
    \label{fig:ex1-3-ten-200-uq}    
\end{figure}

We consider the deterministic initial water surface 
\begin{equation}\label{eq:initex1-3-ten-uq}
    \eta(x,y,0,\xi) = 1,\quad u(x,y,0,\xi) = 0.3,\quad v(x,y,0,\xi) = 0,
\end{equation}
with stochastic bottom
\begin{equation}\label{eq:bottomex1-3-ten-uq}
    B(x,y,\xi) = 0.5e^{-12.5(\xi^{(1)}+1)(x-1)^2-25(\xi^{(2)}+1)(y-0.5)^2}.
\end{equation}
where $\xi = (\xi^{(1)},\xi^{(2)})$ is a random vector, $\xi^{(1)}$ is the distribution with Beta density $(\alpha,\beta) = (3,1)$, $\xi^{(2)}\sim\mathcal{U}(-1,1)$. This example is a multidimensional variant of \Cref{sect:accuracy}. Here, we put the uncertainties on the width of the Gaussian-shape hump. The computational domain is $[0,2]\times [0,1]$. We compute the solution numerically on $100\times 100$ (\Cref{fig:ex1-3-ten-100-uq}) and $200\times 200$ (\Cref{fig:ex1-3-ten-200-uq}) grids up to $t = 0.07$. Compared with the solution in \Cref{fig:ex1-1-uq-time}, where only a simple additive uncertainty is imposed, the contours in \Cref{fig:ex1-3-ten-100-uq} (or \Cref{fig:ex1-3-ten-200-uq}) are more complicated with the introduction of the two dimensional uncertainties on the width of the hump. The largest uncertainties occur near the peak and the lowest part of the mean waves. The region between them, however, has very small uncertainties.
\section{Conclusion}\label{sec:conclusion}
We have introduced a hyperbolicity-preserving SG formulations for the two-dimensional shallow water equations. A sufficient condition for preserving the hyperbolicity of our SG system has been introduced and can be interpreted as a ``stochastic'' positivity condition. Our formulation remains in conservative form. Based on our theoretical results, we develop a well-balanced hyperbolicity-preserving central-upwind solver and illustrate the effectiveness of our solver on some challenging numerical examples. However, our theoretical framework is not limited to central-upwind type schemes, and can be extended to other finite-volume-type schemes (e.g., \cite{shu2003high}) or discontinuous Galerkin schemes (\cite{xing2010positivity}) with modifications on the discretization of the well-balanced source term, 
as well as the CFL-type time condition. 

Currently, our scheme is not designed to handle the dry/wet interfaces, for example, a wave arriving or leaving the shore. The presence of dry or partially dry cells may cause computationally impractical time step even for deterministic problems. In deterministic cases, the ``draining'' time step strategy has been introduced to capture the wet/dry interfaces \cite{bollermann2011finite,bollermann2013well,liu2018well}. Future design of stochastic version of ``draining'' time step may adapt our scheme to more challenging problems.
Another possible extension is to preserve the strict positivity of the water height rather than the positivity at a finite number of quadrature points, which requires an efficient algorithm to preserve the positivity of  high-order polynomials (e.g., \cite{zala2020structure}).
\section*{Acknowledgement}
D.~Dai and A.~Narayan were partially supported by NSF DMS-1848508.
\bibliography{bibfile}

\begin{thebibliography}{54}
\expandafter\ifx\csname natexlab\endcsname\relax\def\natexlab#1{#1}\fi
\providecommand{\url}[1]{\texttt{#1}}
\providecommand{\href}[2]{#2}
\providecommand{\path}[1]{#1}
\providecommand{\DOIprefix}{doi:}
\providecommand{\ArXivprefix}{arXiv:}
\providecommand{\URLprefix}{URL: }
\providecommand{\Pubmedprefix}{pmid:}
\providecommand{\doi}[1]{\href{http://dx.doi.org/#1}{\path{#1}}}
\providecommand{\Pubmed}[1]{\href{pmid:#1}{\path{#1}}}
\providecommand{\bibinfo}[2]{#2}
\ifx\xfnm\relax \def\xfnm[#1]{\unskip,\space#1}\fi
%Type = Article
\bibitem[{Babuska et~al.(2004)Babuska, Tempone and
  Zouraris}]{babuska2004galerkin}
\bibinfo{author}{Babuska, I.}, \bibinfo{author}{Tempone, R.},
  \bibinfo{author}{Zouraris, G.E.}, \bibinfo{year}{2004}.
\newblock \bibinfo{title}{{Galerkin Finite Element Approximations of Stochastic
  Elliptic Partial Differential Equations}}.
\newblock \bibinfo{journal}{SIAM Journal on Numerical Analysis}
  \bibinfo{volume}{42}, \bibinfo{pages}{800--825}.
\newblock \DOIprefix\doi{https://doi.org/10.1137/S0036142902418680}.
%Type = Article
\bibitem[{Bollermann et~al.(2013)Bollermann, Chen, Kurganov and
  Noelle}]{bollermann2013well}
\bibinfo{author}{Bollermann, A.}, \bibinfo{author}{Chen, G.},
  \bibinfo{author}{Kurganov, A.}, \bibinfo{author}{Noelle, S.},
  \bibinfo{year}{2013}.
\newblock \bibinfo{title}{{A Well-Balanced Reconstruction of Wet/Dry Fronts for
  the Shallow Water Equations}}.
\newblock \bibinfo{journal}{Journal of Scientific Computing}
  \bibinfo{volume}{56}, \bibinfo{pages}{267--290}.
\newblock \DOIprefix\doi{https://doi.org/10.1007/s10915-012-9677-5}.
%Type = Article
\bibitem[{Bollermann et~al.(2011)Bollermann, Noelle and
  Luk{\'a}{\v{c}}ov{\'a}-Medvid’ov{\'a}}]{bollermann2011finite}
\bibinfo{author}{Bollermann, A.}, \bibinfo{author}{Noelle, S.},
  \bibinfo{author}{Luk{\'a}{\v{c}}ov{\'a}-Medvid’ov{\'a}, M.},
  \bibinfo{year}{2011}.
\newblock \bibinfo{title}{{Finite Volume Evolution Galerkin Methods for the
  Shallow Water Equations with Dry Beds}}.
\newblock \bibinfo{journal}{Communications in Computational Physics}
  \bibinfo{volume}{10}, \bibinfo{pages}{371--404}.
\newblock \DOIprefix\doi{https://doi.org/10.4208/cicp.220210.020710a}.
%Type = Article
\bibitem[{Bryson et~al.(2011)Bryson, Epshteyn, Kurganov and
  Petrova}]{M2AN_2011__45_3_423_0}
\bibinfo{author}{Bryson, S.}, \bibinfo{author}{Epshteyn, Y.},
  \bibinfo{author}{Kurganov, A.}, \bibinfo{author}{Petrova, G.},
  \bibinfo{year}{2011}.
\newblock \bibinfo{title}{Well-balanced positivity preserving central-upwind
  scheme on triangular grids for the saint-venant system}.
\newblock \bibinfo{journal}{ESAIM: Mathematical Modelling and Numerical
  Analysis - Mod\'elisation Math\'ematique et Analyse Num\'erique}
  \bibinfo{volume}{45}, \bibinfo{pages}{423--446}.
\newblock \DOIprefix\doi{https://doi.org/10.1051/m2an/2010060}.
%Type = Article
\bibitem[{Chertock et~al.(2015a)Chertock, Cui, Kurganov and
  Wu}]{chertock2015well}
\bibinfo{author}{Chertock, A.}, \bibinfo{author}{Cui, S.},
  \bibinfo{author}{Kurganov, A.}, \bibinfo{author}{Wu, T.},
  \bibinfo{year}{2015}a.
\newblock \bibinfo{title}{{Well-balanced positivity preserving central-upwind
  scheme for the shallow water system with friction terms}}.
\newblock \bibinfo{journal}{International Journal for Numerical Methods in
  Fluids} \bibinfo{volume}{78}, \bibinfo{pages}{355--383}.
\newblock \DOIprefix\doi{https://doi.org/10.1002/fld.4023}.
%Type = Article
\bibitem[{Chertock et~al.(2015b)Chertock, Jin and
  Kurganov}]{chertock2015welluq}
\bibinfo{author}{Chertock, A.}, \bibinfo{author}{Jin, S.},
  \bibinfo{author}{Kurganov, A.}, \bibinfo{year}{2015}b.
\newblock \bibinfo{title}{{A well-balanced operator splitting based stochastic
  Galerkin method for the one-dimensional Saint-Venant system with
  uncertainty}}.
\newblock \bibinfo{journal}{preprint} \URLprefix
  \url{https://ki-net.umd.edu/pubs/files/SW.pdf}.
%Type = Article
\bibitem[{Chertock et~al.(2015c)Chertock, Jin and
  Kurganov}]{chertock2015operator}
\bibinfo{author}{Chertock, A.}, \bibinfo{author}{Jin, S.},
  \bibinfo{author}{Kurganov, A.}, \bibinfo{year}{2015}c.
\newblock \bibinfo{title}{{An operator splitting based stochastic Galerkin
  method for the one-dimensional compressible Euler equations with
  uncertainty}}.
\newblock \bibinfo{journal}{preprint} , \bibinfo{pages}{1--21}\URLprefix
  \url{http://www.math.wisc.edu/~jin/PS/Euler-UQ.pdf}.
%Type = Article
\bibitem[{Dai et~al.(2021)Dai, Epshteyn and Narayan}]{doi:10.1137/20M1360736}
\bibinfo{author}{Dai, D.}, \bibinfo{author}{Epshteyn, Y.},
  \bibinfo{author}{Narayan, A.}, \bibinfo{year}{2021}.
\newblock \bibinfo{title}{{Hyperbolicity-Preserving and Well-Balanced
  Stochastic Galerkin Method for Shallow Water Equations}}.
\newblock \bibinfo{journal}{SIAM Journal on Scientific Computing}
  \bibinfo{volume}{43}, \bibinfo{pages}{A929--A952}.
\newblock \DOIprefix\doi{https://doi.org/10.1137/20M1360736}.
%Type = Article
\bibitem[{Debusschere et~al.(2004)Debusschere, Najm, P{\'e}bay, Knio, Ghanem
  and Le~Ma{\^{\i}}tre}]{debusschere2004numerical}
\bibinfo{author}{Debusschere, B.J.}, \bibinfo{author}{Najm, H.N.},
  \bibinfo{author}{P{\'e}bay, P.P.}, \bibinfo{author}{Knio, O.M.},
  \bibinfo{author}{Ghanem, R.G.}, \bibinfo{author}{Le~Ma{\^{\i}}tre, O.P.},
  \bibinfo{year}{2004}.
\newblock \bibinfo{title}{{Numerical Challenges in the Use of Polynomial Chaos
  Representations for Stochastic Processes}}.
\newblock \bibinfo{journal}{SIAM journal on scientific computing}
  \bibinfo{volume}{26}, \bibinfo{pages}{698--719}.
\newblock \DOIprefix\doi{https://doi.org/10.1137/S1064827503427741}.
%Type = Incollection
\bibitem[{Despr{\'e}s et~al.(2013)Despr{\'e}s, Po{\"e}tte and
  Lucor}]{despres2013robust}
\bibinfo{author}{Despr{\'e}s, B.}, \bibinfo{author}{Po{\"e}tte, G.},
  \bibinfo{author}{Lucor, D.}, \bibinfo{year}{2013}.
\newblock \bibinfo{title}{{Robust Uncertainty Propagation in Systems of
  Conservation Laws with the Entropy Closure Method}}, in:
  \bibinfo{booktitle}{Uncertainty Quantification in Computational Fluid
  Dynamics}. \bibinfo{publisher}{Springer}, pp. \bibinfo{pages}{105--149}.
%Type = Article
\bibitem[{Eigel et~al.(2014)Eigel, Gittelson, Schwab and
  Zander}]{eigel2014adaptive}
\bibinfo{author}{Eigel, M.}, \bibinfo{author}{Gittelson, C.J.},
  \bibinfo{author}{Schwab, C.}, \bibinfo{author}{Zander, E.},
  \bibinfo{year}{2014}.
\newblock \bibinfo{title}{{Adaptive stochastic Galerkin FEM}}.
\newblock \bibinfo{journal}{Computer Methods in Applied Mechanics and
  Engineering} \bibinfo{volume}{270}, \bibinfo{pages}{247--269}.
\newblock \DOIprefix\doi{https://doi.org/10.1016/j.cma.2013.11.015}.
%Type = Article
\bibitem[{Ernst et~al.(2012)Ernst, Mugler, Starkloff and
  Ullmann}]{ernst_convergence_2012}
\bibinfo{author}{Ernst, O.G.}, \bibinfo{author}{Mugler, A.},
  \bibinfo{author}{Starkloff, H.J.}, \bibinfo{author}{Ullmann, E.},
  \bibinfo{year}{2012}.
\newblock \bibinfo{title}{{On the {Convergence} of {Generalized} {Polynomial}
  {Chaos} {Expansions}}}.
\newblock \bibinfo{journal}{ESAIM: Mathematical Modelling and Numerical
  Analysis} \bibinfo{volume}{46}, \bibinfo{pages}{317--339}.
\newblock \DOIprefix\doi{https://doi.org/10.1051/m2an/2011045}.
%Type = Article
\bibitem[{Gerster and Herty(2020)}]{gerster2020entropies}
\bibinfo{author}{Gerster, S.}, \bibinfo{author}{Herty, M.},
  \bibinfo{year}{2020}.
\newblock \bibinfo{title}{{Entropies and Symmetrization of Hyperbolic
  Stochastic Galerkin Formulations}}.
\newblock \bibinfo{journal}{Communications in Computational Physics}
  \bibinfo{volume}{27}, \bibinfo{pages}{639--671}.
\newblock \DOIprefix\doi{https://doi.org/10.4208/cicp.OA-2019-0047}.
%Type = Article
\bibitem[{Gerster et~al.(2019)Gerster, Herty and
  Sikstel}]{gerster2019hyperbolic}
\bibinfo{author}{Gerster, S.}, \bibinfo{author}{Herty, M.},
  \bibinfo{author}{Sikstel, A.}, \bibinfo{year}{2019}.
\newblock \bibinfo{title}{{Hyperbolic stochastic Galerkin formulation for the
  p-system}}.
\newblock \bibinfo{journal}{Journal of Computational Physics}
  \bibinfo{volume}{395}, \bibinfo{pages}{186--204}.
\newblock \DOIprefix\doi{https://doi.org/10.1016/j.jcp.2019.05.049}.
%Type = Book
\bibitem[{Ghanem and Spanos(1991)}]{ghanem1991stochastic}
\bibinfo{author}{Ghanem, R.G.}, \bibinfo{author}{Spanos, P.D.},
  \bibinfo{year}{1991}.
\newblock \bibinfo{title}{{Stochastic Finite Elements: A Spectral Approach}}.
\newblock \bibinfo{publisher}{Springer}.
%Type = Article
\bibitem[{Ghazizadeh et~al.(2020)Ghazizadeh, Mohammadian and
  Kurganov}]{ghazizadeh2020adaptive}
\bibinfo{author}{Ghazizadeh, M.A.}, \bibinfo{author}{Mohammadian, A.},
  \bibinfo{author}{Kurganov, A.}, \bibinfo{year}{2020}.
\newblock \bibinfo{title}{An adaptive well-balanced positivity preserving
  central-upwind scheme on quadtree grids for shallow water equations}.
\newblock \bibinfo{journal}{Computers \& Fluids} \bibinfo{volume}{208},
  \bibinfo{pages}{104633}.
%Type = Article
\bibitem[{Gottlieb et~al.(2001)Gottlieb, Shu and Tadmor}]{gottlieb2001strong}
\bibinfo{author}{Gottlieb, S.}, \bibinfo{author}{Shu, C.W.},
  \bibinfo{author}{Tadmor, E.}, \bibinfo{year}{2001}.
\newblock \bibinfo{title}{{Strong Stability-Preserving High-Order Time
  Discretization Methods}}.
\newblock \bibinfo{journal}{SIAM review} \bibinfo{volume}{43},
  \bibinfo{pages}{89--112}.
\newblock \DOIprefix\doi{https://doi.org/10.1137/S003614450036757X}.
%Type = Article
\bibitem[{Hu and Jin(2016)}]{hu2016stochastic}
\bibinfo{author}{Hu, J.}, \bibinfo{author}{Jin, S.}, \bibinfo{year}{2016}.
\newblock \bibinfo{title}{{A stochastic Galerkin method for the Boltzmann
  equation with uncertainty}}.
\newblock \bibinfo{journal}{Journal of Computational Physics}
  \bibinfo{volume}{315}, \bibinfo{pages}{150--168}.
\newblock \DOIprefix\doi{https://doi.org/10.1016/j.jcp.2016.03.047}.
%Type = Article
\bibitem[{Jin and Shu(2019)}]{jin2019study}
\bibinfo{author}{Jin, S.}, \bibinfo{author}{Shu, R.}, \bibinfo{year}{2019}.
\newblock \bibinfo{title}{{A Study of Hyperbolicity of Kinetic Stochastic
  Galerkin System for the Isentropic Euler Equations with Uncertainty}}.
\newblock \bibinfo{journal}{Chinese Annals of Mathematics, Series B}
  \bibinfo{volume}{40}, \bibinfo{pages}{765--780}.
\newblock \DOIprefix\doi{https://doi.org/10.1007/s11401-019-0159-z}.
%Type = Article
\bibitem[{Kurganov(2018)}]{kurganov2018finite}
\bibinfo{author}{Kurganov, A.}, \bibinfo{year}{2018}.
\newblock \bibinfo{title}{{Finite-volume schemes for shallow-water equations}}.
\newblock \bibinfo{journal}{Acta Numerica} \bibinfo{volume}{27},
  \bibinfo{pages}{289--351}.
\newblock \DOIprefix\doi{https://doi.org/10.1017/S0962492918000028}.
%Type = Article
\bibitem[{Kurganov and Levy(2002)}]{kurganov2002central}
\bibinfo{author}{Kurganov, A.}, \bibinfo{author}{Levy, D.},
  \bibinfo{year}{2002}.
\newblock \bibinfo{title}{{Central-upwind schemes for the Saint-Venant
  system}}.
\newblock \bibinfo{journal}{ESAIM: Mathematical Modelling and Numerical
  Analysis} \bibinfo{volume}{36}, \bibinfo{pages}{397--425}.
\newblock \DOIprefix\doi{https://doi.org/10.1051/m2an:2002019}.
%Type = Article
\bibitem[{Kurganov and Lin(2007)}]{kurganov2007reduction}
\bibinfo{author}{Kurganov, A.}, \bibinfo{author}{Lin, C.T.},
  \bibinfo{year}{2007}.
\newblock \bibinfo{title}{{On the Reduction of Numerical Dissipation in
  Central-Upwind Schemes}}.
\newblock \bibinfo{journal}{Communications in Computational Physics}
  \bibinfo{volume}{2}, \bibinfo{pages}{141--163}.
%Type = Article
\bibitem[{Kurganov et~al.(2001)Kurganov, Noelle and
  Petrova}]{kurganov2001semidiscrete}
\bibinfo{author}{Kurganov, A.}, \bibinfo{author}{Noelle, S.},
  \bibinfo{author}{Petrova, G.}, \bibinfo{year}{2001}.
\newblock \bibinfo{title}{{Semidiscrete Central-Upwind Schemes for Hyperbolic
  Conservation Laws and Hamilton--Jacobi Equations}}.
\newblock \bibinfo{journal}{SIAM Journal on Scientific Computing}
  \bibinfo{volume}{23}, \bibinfo{pages}{707--740}.
\newblock \DOIprefix\doi{https://doi.org/10.1137/S1064827500373413}.
%Type = Article
\bibitem[{Kurganov and Petrova(2007)}]{kurganov2007second}
\bibinfo{author}{Kurganov, A.}, \bibinfo{author}{Petrova, G.},
  \bibinfo{year}{2007}.
\newblock \bibinfo{title}{{A second-order well-balanced positivity preserving
  central-upwind scheme for the Saint-Venant system}}.
\newblock \bibinfo{journal}{Communications in Mathematical Sciences}
  \bibinfo{volume}{5}, \bibinfo{pages}{133--160}.
%Type = Article
\bibitem[{Kurganov et~al.(2007)Kurganov, Petrova and
  Popov}]{kurganov2007adaptive}
\bibinfo{author}{Kurganov, A.}, \bibinfo{author}{Petrova, G.},
  \bibinfo{author}{Popov, B.}, \bibinfo{year}{2007}.
\newblock \bibinfo{title}{{Adaptive Semidiscrete Central-Upwind Schemes for
  Nonconvex Hyperbolic Conservation Laws}}.
\newblock \bibinfo{journal}{SIAM Journal on Scientific Computing}
  \bibinfo{volume}{29}, \bibinfo{pages}{2381--2401}.
\newblock \DOIprefix\doi{https://doi.org/10.1137/040614189}.
%Type = Article
\bibitem[{Kurganov and Tadmor(2000)}]{kurganov2000new}
\bibinfo{author}{Kurganov, A.}, \bibinfo{author}{Tadmor, E.},
  \bibinfo{year}{2000}.
\newblock \bibinfo{title}{{New High-Resolution Central Schemes for Nonlinear
  Conservation Laws and Convection–Diffusion Equations}}.
\newblock \bibinfo{journal}{Journal of Computational Physics}
  \bibinfo{volume}{160}, \bibinfo{pages}{241--282}.
\newblock \DOIprefix\doi{https://doi.org/10.1006/jcph.2000.6459}.
%Type = Article
\bibitem[{Kusch et~al.(2020)Kusch, McClarren and Frank}]{kusch2020filtered}
\bibinfo{author}{Kusch, J.}, \bibinfo{author}{McClarren, R.G.},
  \bibinfo{author}{Frank, M.}, \bibinfo{year}{2020}.
\newblock \bibinfo{title}{{Filtered stochastic Galerkin methods for hyperbolic
  equations}}.
\newblock \bibinfo{journal}{Journal of Computational Physics}
  \bibinfo{volume}{403}, \bibinfo{pages}{109073}.
\newblock \DOIprefix\doi{https://doi.org/10.1016/j.jcp.2019.109073}.
%Type = Article
\bibitem[{Kusch and Schlachter(2020)}]{kusch2020oscillation}
\bibinfo{author}{Kusch, J.}, \bibinfo{author}{Schlachter, L.},
  \bibinfo{year}{2020}.
\newblock \bibinfo{title}{{Oscillation Mitigation of Hyperbolicity-Preserving
  Intrusive Uncertainty Quantification Methods for Systems of Conservation
  Laws}}.
\newblock \bibinfo{journal}{arXiv:2008.07845} .
%Type = Article
\bibitem[{Le~Ma{\^{\i}}tre et~al.(2004)Le~Ma{\^{\i}}tre, Knio, Najm and
  Ghanem}]{le2004uncertainty}
\bibinfo{author}{Le~Ma{\^{\i}}tre, O.}, \bibinfo{author}{Knio, O.},
  \bibinfo{author}{Najm, H.}, \bibinfo{author}{Ghanem, R.},
  \bibinfo{year}{2004}.
\newblock \bibinfo{title}{{Uncertainty propagation using Wiener--Haar
  expansions}}.
\newblock \bibinfo{journal}{Journal of computational Physics}
  \bibinfo{volume}{197}, \bibinfo{pages}{28--57}.
\newblock \DOIprefix\doi{https://doi.org/10.1016/j.jcp.2003.11.033}.
%Type = Book
\bibitem[{Le~Ma{\^{\i}}tre and Knio(2010)}]{le2010spectral}
\bibinfo{author}{Le~Ma{\^{\i}}tre, O.}, \bibinfo{author}{Knio, O.M.},
  \bibinfo{year}{2010}.
\newblock \bibinfo{title}{{Spectral Methods for Uncertainty Quantification:
  With Applications to Computational Fluid Dynamics}}.
\newblock \bibinfo{publisher}{Springer Science \& Business Media}.
%Type = Article
\bibitem[{Liu et~al.(2018)Liu, Albright, Epshteyn and Kurganov}]{liu2018well}
\bibinfo{author}{Liu, X.}, \bibinfo{author}{Albright, J.},
  \bibinfo{author}{Epshteyn, Y.}, \bibinfo{author}{Kurganov, A.},
  \bibinfo{year}{2018}.
\newblock \bibinfo{title}{{Well-balanced positivity preserving central-upwind
  scheme with a novel wet/dry reconstruction on triangular grids for the
  Saint-Venant system}}.
\newblock \bibinfo{journal}{Journal of Computational Physics}
  \bibinfo{volume}{374}, \bibinfo{pages}{213--236}.
\newblock \DOIprefix\doi{https://doi.org/10.1016/j.jcp.2018.07.038}.
%Type = Article
\bibitem[{Mishra et~al.(2012)Mishra, Schwab and Sukys}]{mishra2012multilevel}
\bibinfo{author}{Mishra, S.}, \bibinfo{author}{Schwab, C.},
  \bibinfo{author}{Sukys, J.}, \bibinfo{year}{2012}.
\newblock \bibinfo{title}{{Multilevel Monte Carlo Finite Volume Methods for
  Shallow Water Equations with Uncertain Topography in Multi-dimensions}}.
\newblock \bibinfo{journal}{SIAM Journal on Scientific Computing}
  \bibinfo{volume}{34}, \bibinfo{pages}{B761--B784}.
\newblock \DOIprefix\doi{https://doi.org/10.1137/110857295}.
%Type = Article
\bibitem[{Nessyahu and Tadmor(1990)}]{nessyahu1990non}
\bibinfo{author}{Nessyahu, H.}, \bibinfo{author}{Tadmor, E.},
  \bibinfo{year}{1990}.
\newblock \bibinfo{title}{{Non-oscillatory central differencing for hyperbolic
  conservation laws}}.
\newblock \bibinfo{journal}{Journal of computational physics}
  \bibinfo{volume}{87}, \bibinfo{pages}{408--463}.
\newblock \DOIprefix\doi{https://doi.org/10.1016/0021-9991(90)90260-8}.
%Type = Article
\bibitem[{Nobile et~al.(2008)Nobile, Tempone and Webster}]{nobile2008sparse}
\bibinfo{author}{Nobile, F.}, \bibinfo{author}{Tempone, R.},
  \bibinfo{author}{Webster, C.G.}, \bibinfo{year}{2008}.
\newblock \bibinfo{title}{{A Sparse Grid Stochastic Collocation Method for
  Partial Differential Equations with Random Input Data}}.
\newblock \bibinfo{journal}{SIAM Journal on Numerical Analysis}
  \bibinfo{volume}{46}, \bibinfo{pages}{2309--2345}.
\newblock \DOIprefix\doi{https://doi.org/10.1137/060663660}.
%Type = Article
\bibitem[{Pettersson et~al.(2014)Pettersson, Iaccarino and
  Nordstr{\"o}m}]{pettersson2014stochastic}
\bibinfo{author}{Pettersson, P.}, \bibinfo{author}{Iaccarino, G.},
  \bibinfo{author}{Nordstr{\"o}m, J.}, \bibinfo{year}{2014}.
\newblock \bibinfo{title}{{A stochastic Galerkin method for the Euler equations
  with Roe variable transformation}}.
\newblock \bibinfo{journal}{Journal of Computational Physics}
  \bibinfo{volume}{257}, \bibinfo{pages}{481--500}.
\newblock \DOIprefix\doi{https://doi.org/10.1016/j.jcp.2013.10.011}.
%Type = Phdthesis
\bibitem[{Po{\"e}tte(2019)}]{poette2019contribution}
\bibinfo{author}{Po{\"e}tte, G.}, \bibinfo{year}{2019}.
\newblock \bibinfo{title}{{Contribution to the mathematical and numerical
  analysis of uncertain systems of conservation laws and of the linear and
  nonlinear Boltzmann equation}}.
\newblock Ph.D. thesis. Universit\'{e} de Bordeaux.
\newblock \URLprefix
  \url{https://hal.archives-ouvertes.fr/tel-02288678/document}.
%Type = Article
\bibitem[{Po{\"e}tte et~al.(2009)Po{\"e}tte, Despr{\'e}s and
  Lucor}]{poette2009uncertainty}
\bibinfo{author}{Po{\"e}tte, G.}, \bibinfo{author}{Despr{\'e}s, B.},
  \bibinfo{author}{Lucor, D.}, \bibinfo{year}{2009}.
\newblock \bibinfo{title}{{Uncertainty quantification for systems of
  conservation laws}}.
\newblock \bibinfo{journal}{Journal of Computational Physics}
  \bibinfo{volume}{228}, \bibinfo{pages}{2443--2467}.
\newblock \DOIprefix\doi{https://doi.org/10.1016/j.jcp.2008.12.018}.
%Type = Article
\bibitem[{Schlachter and Schneider(2018)}]{schlachter2018hyperbolicity}
\bibinfo{author}{Schlachter, L.}, \bibinfo{author}{Schneider, F.},
  \bibinfo{year}{2018}.
\newblock \bibinfo{title}{{A hyperbolicity-preserving stochastic Galerkin
  approximation for uncertain hyperbolic systems of equations}}.
\newblock \bibinfo{journal}{Journal of Computational Physics}
  \bibinfo{volume}{375}, \bibinfo{pages}{80--98}.
\newblock \DOIprefix\doi{https://doi.org/10.1016/j.jcp.2018.07.026}.
%Type = Article
\bibitem[{Schlachter et~al.(2020)Schlachter, Schneider and
  Kolb}]{schlachter2020weighted}
\bibinfo{author}{Schlachter, L.}, \bibinfo{author}{Schneider, F.},
  \bibinfo{author}{Kolb, O.}, \bibinfo{year}{2020}.
\newblock \bibinfo{title}{{Weighted Essentially Non-Oscillatory stochastic
  Galerkin approximation for hyperbolic conservation laws}}.
\newblock \bibinfo{journal}{Journal of Computational Physics}
  \bibinfo{volume}{419}, \bibinfo{pages}{109663}.
\newblock \DOIprefix\doi{https://doi.org/10.1016/j.jcp.2020.109663}.
%Type = Article
\bibitem[{Shirkhani et~al.(2016)Shirkhani, Mohammadian, Seidou and
  Kurganov}]{SHIRKHANI201625}
\bibinfo{author}{Shirkhani, H.}, \bibinfo{author}{Mohammadian, A.},
  \bibinfo{author}{Seidou, O.}, \bibinfo{author}{Kurganov, A.},
  \bibinfo{year}{2016}.
\newblock \bibinfo{title}{{A well-balanced positivity-preserving central-upwind
  scheme for shallow water equations on unstructured quadrilateral grids}}.
\newblock \bibinfo{journal}{Computers \& Fluids} \bibinfo{volume}{126},
  \bibinfo{pages}{25--40}.
\newblock \URLprefix
  \url{https://www.sciencedirect.com/science/article/pii/S0045793015003886},
  \DOIprefix\doi{https://doi.org/10.1016/j.compfluid.2015.11.017}.
%Type = Article
\bibitem[{Shu(2003)}]{shu2003high}
\bibinfo{author}{Shu, C.W.}, \bibinfo{year}{2003}.
\newblock \bibinfo{title}{{High-order Finite Difference and Finite Volume WENO
  Schemes and Discontinuous Galerkin Methods for CFD}}.
\newblock \bibinfo{journal}{International Journal of Computational Fluid
  Dynamics} \bibinfo{volume}{17}, \bibinfo{pages}{107--118}.
\newblock \DOIprefix\doi{https://doi.org/10.1080/1061856031000104851}.
%Type = Article
\bibitem[{Shu et~al.(2017)Shu, Hu and Jin}]{shu2017stochastic}
\bibinfo{author}{Shu, R.}, \bibinfo{author}{Hu, J.}, \bibinfo{author}{Jin, S.},
  \bibinfo{year}{2017}.
\newblock \bibinfo{title}{{A Stochastic Galerkin Method for the Boltzmann
  Equation with Multi-Dimensional Random Inputs Using Sparse Wavelet Bases}}.
\newblock \bibinfo{journal}{Numerical Mathematics: Theory, Methods and
  Applications} \bibinfo{volume}{10}, \bibinfo{pages}{465--488}.
\newblock \DOIprefix\doi{https://doi.org/10.4208/nmtma.2017.s12}.
%Type = Book
\bibitem[{Sullivan(2015)}]{sullivan2015introduction}
\bibinfo{author}{Sullivan, T.J.}, \bibinfo{year}{2015}.
\newblock \bibinfo{title}{{Introduction to Uncertainty Quantification}}.
  volume~\bibinfo{volume}{63}.
\newblock \bibinfo{publisher}{Springer}.
\newblock \DOIprefix\doi{https://doi.org/10.1007/978-3-319-23395-6}.
%Type = Article
\bibitem[{Tryoen et~al.(2010)Tryoen, Le~Ma{\^{\i}}tre, Ndjinga and
  Ern}]{tryoen2010intrusive}
\bibinfo{author}{Tryoen, J.}, \bibinfo{author}{Le~Ma{\^{\i}}tre, O.},
  \bibinfo{author}{Ndjinga, M.}, \bibinfo{author}{Ern, A.},
  \bibinfo{year}{2010}.
\newblock \bibinfo{title}{{Intrusive Galerkin methods with upwinding for
  uncertain nonlinear hyperbolic systems}}.
\newblock \bibinfo{journal}{Journal of Computational Physics}
  \bibinfo{volume}{229}, \bibinfo{pages}{6485--6511}.
\newblock \DOIprefix\doi{https://doi.org/10.1016/j.jcp.2010.05.007}.
%Type = Article
\bibitem[{Van~Leer(1974)}]{van1974towards}
\bibinfo{author}{Van~Leer, B.}, \bibinfo{year}{1974}.
\newblock \bibinfo{title}{{Towards the ultimate conservative difference scheme.
  II. Monotonicity and conservation combined in a second-order scheme}}.
\newblock \bibinfo{journal}{Journal of computational physics}
  \bibinfo{volume}{14}, \bibinfo{pages}{361--370}.
%Type = Article
\bibitem[{Wan and Karniadakis(2005)}]{wan2005adaptive}
\bibinfo{author}{Wan, X.}, \bibinfo{author}{Karniadakis, G.E.},
  \bibinfo{year}{2005}.
\newblock \bibinfo{title}{{An adaptive multi-element generalized polynomial
  chaos method for stochastic differential equations}}.
\newblock \bibinfo{journal}{Journal of Computational Physics}
  \bibinfo{volume}{209}, \bibinfo{pages}{617--642}.
\newblock \DOIprefix\doi{https://doi.org/10.1016/j.jcp.2005.03.023}.
%Type = Article
\bibitem[{Wiener(1938)}]{wiener1938homogeneous}
\bibinfo{author}{Wiener, N.}, \bibinfo{year}{1938}.
\newblock \bibinfo{title}{{The Homogeneous Chaos}}.
\newblock \bibinfo{journal}{American Journal of Mathematics}
  \bibinfo{volume}{60}, \bibinfo{pages}{897--936}.
\newblock \DOIprefix\doi{https://doi.org/10.2307/2371268}.
%Type = Article
\bibitem[{Wu et~al.(2017)Wu, Tang and Xiu}]{wu2017stochastic}
\bibinfo{author}{Wu, K.}, \bibinfo{author}{Tang, H.}, \bibinfo{author}{Xiu,
  D.}, \bibinfo{year}{2017}.
\newblock \bibinfo{title}{{A stochastic Galerkin method for first-order
  quasilinear hyperbolic systems with uncertainty}}.
\newblock \bibinfo{journal}{Journal of Computational Physics}
  \bibinfo{volume}{345}, \bibinfo{pages}{224--244}.
\newblock \DOIprefix\doi{https://doi.org/10.1016/j.jcp.2017.05.027}.
%Type = Article
\bibitem[{Xing et~al.(2010)Xing, Zhang and Shu}]{xing2010positivity}
\bibinfo{author}{Xing, Y.}, \bibinfo{author}{Zhang, X.}, \bibinfo{author}{Shu,
  C.W.}, \bibinfo{year}{2010}.
\newblock \bibinfo{title}{{Positivity-preserving high order well-balanced
  discontinuous Galerkin methods for the shallow water equations}}.
\newblock \bibinfo{journal}{Advances in Water Resources} \bibinfo{volume}{33},
  \bibinfo{pages}{1476--1493}.
\newblock \DOIprefix\doi{https://doi.org/10.1016/j.advwatres.2010.08.005}.
%Type = Book
\bibitem[{Xiu(2010)}]{xiu2010numerical}
\bibinfo{author}{Xiu, D.}, \bibinfo{year}{2010}.
\newblock \bibinfo{title}{Numerical Methods for Stochastic Computations: A
  Spectral Method Approach}.
\newblock \bibinfo{publisher}{Princeton university press}.
%Type = Article
\bibitem[{Xiu and Hesthaven(2005)}]{xiu2005high}
\bibinfo{author}{Xiu, D.}, \bibinfo{author}{Hesthaven, J.S.},
  \bibinfo{year}{2005}.
\newblock \bibinfo{title}{{High-Order Collocation Methods for Differential
  Equations with Random Inputs}}.
\newblock \bibinfo{journal}{SIAM Journal on Scientific Computing}
  \bibinfo{volume}{27}, \bibinfo{pages}{1118--1139}.
\newblock \DOIprefix\doi{https://doi.org/10.1137/040615201}.
%Type = Article
\bibitem[{Xiu and Karniadakis(2002)}]{xiu2002wiener}
\bibinfo{author}{Xiu, D.}, \bibinfo{author}{Karniadakis, G.E.},
  \bibinfo{year}{2002}.
\newblock \bibinfo{title}{{The Wiener--Askey Polynomial Chaos for Stochastic
  Differential Equations}}.
\newblock \bibinfo{journal}{SIAM journal on scientific computing}
  \bibinfo{volume}{24}, \bibinfo{pages}{619--644}.
\newblock \DOIprefix\doi{https://doi.org/10.1137/S1064827501387826}.
%Type = Article
\bibitem[{Xiu and Shen(2009)}]{xiu2009efficient}
\bibinfo{author}{Xiu, D.}, \bibinfo{author}{Shen, J.}, \bibinfo{year}{2009}.
\newblock \bibinfo{title}{{Efficient stochastic Galerkin methods for random
  diffusion equations}}.
\newblock \bibinfo{journal}{Journal of Computational Physics}
  \bibinfo{volume}{228}, \bibinfo{pages}{266--281}.
\newblock \DOIprefix\doi{https://doi.org/10.1016/j.jcp.2008.09.008}.
%Type = Article
\bibitem[{Zala et~al.(2020)Zala, Kirby and Narayan}]{zala2020structure}
\bibinfo{author}{Zala, V.}, \bibinfo{author}{Kirby, M.},
  \bibinfo{author}{Narayan, A.}, \bibinfo{year}{2020}.
\newblock \bibinfo{title}{{Structure-Preserving Function Approximation via
  Convex Optimization}}.
\newblock \bibinfo{journal}{SIAM Journal on Scientific Computing}
  \bibinfo{volume}{42}, \bibinfo{pages}{A3006--A3029}.
\newblock \DOIprefix\doi{10.1137/19M130128X}.

\end{thebibliography}
\end{document}